\UseRawInputEncoding
\documentclass[12pt,reqno]{amsart}
\usepackage{amsmath,amsfonts,amsthm,amsopn,amssymb}
\usepackage{cite,marginnote}
\usepackage{color,enumitem,graphicx}
\usepackage[colorlinks=true,urlcolor=blue,
citecolor=red,linkcolor=blue,linktocpage,pdfpagelabels,
bookmarksnumbered,bookmarksopen]{hyperref}
\usepackage[english]{babel}

\usepackage[left=2.9cm,right=2.9cm,top=2.8cm,bottom=2.8cm]{geometry}
\usepackage[hyperpageref]{backref}

\allowdisplaybreaks
\numberwithin{equation}{section}

\makeindex
\makeindex
\newenvironment{Abstract}{\textbf{Abstract}\mbox{  }}{ }
\newenvironment{key words}{\textbf{Keywords}\mbox{  }}{ }

\newtheorem{theorem}{Theorem}[section]
\newtheorem{definition}[theorem]{Definition}
\newtheorem{lemma}[theorem]{Lemma}
\newtheorem{corollary}[theorem]{Corollary}
\newtheorem{proposition}[theorem]{Proposition}
\renewenvironment{proof}{\noindent{\textbf{Proof.}}}{\hfill$\Box$}
\newtheorem{remark}[theorem]{Remark}

\newcommand{\s}{\section}

\newcommand{\bt}{\begin{theorem}}
\newcommand{\et}{\end{theorem}}
\newcommand{\bl}{\begin{lemma}}
\newcommand{\el}{\end{lemma}}
\newcommand{\bd}{\begin{definition}}
\newcommand{\ed}{\end{definition}}
\newcommand{\bc}{\begin{corollary}}
\newcommand{\ec}{\end{corollary}}
\newcommand{\bp}{\begin{proof}}
\newcommand{\ep}{\end{proof}}
\newcommand{\bx}{\begin{example}}
\newcommand{\ex}{\end{example}}
\newcommand{\bi}{\begin{exercise}}
\newcommand{\ei}{\end{exercise}}
\newcommand{\bo}{\begin{proposition}}
\newcommand{\eo}{\end{proposition}}
\newcommand{\br}{\begin{remark}}
\newcommand{\er}{\end{remark}}
\newcommand{\beq}{\begin{equation}}
\newcommand{\eeq}{\end{equation}}
\newcommand{\ba}{\begin{align}}
\newcommand{\ea}{\end{align}}
\newcommand{\bn}{\begin{enumerate}}
\newcommand{\en}{\end{enumerate}}
\newcommand{\bg}{\begin{align*}}
\newcommand{\bcs}{\begin{cases}}
\newcommand{\ecs}{\end{cases}}

\newcommand{\bean}{\begin{eqnarray*}}
\newcommand{\eean}{\end{eqnarray*}}

\def\bd{\mathrm{bd}\,}

\begin{document}

\title[Normalized ground state]{Sharp interaction estimates and their application: existence of normalized ground states to coupled Schr\"odinger systems with potentials}
\author[Y.~B.~Deng]{Yinbin Deng}
\author[Q.~H.~He]{Qihan He}
\author[X.~X.~Zhong]{Xuexiu Zhong}
\address[Y.~B.~Deng]{\newline\indent College of Mathematics and Statistics,
Xinyang Normal University,
Xinyang 464000, P. R. China.
\newline\indent
School of Mathematics and Statistics \& Hubei Key Laboratory of Mathematical Sciences,
Central China Normal University,
Wuhan 430079, P. R. China.}
\email{\href{mailto:ybdeng@mail.ccnu.edu.cn}{ybdeng@mail.ccnu.edu.cn}}
\address[Q.~H.~He]{\newline\indent ~College of Mathematics and Information Science, Guangxi Center for Mathematical Research,
Guangxi University,
Nanning, 530003, P. R. China.}
\email{\href{mailto:heqihan277@163.com}{heqihan277@163.com}}
\address[X.~X.~Zhong]{\newline\indent South China Research Center for Applied Mathematics and Interdisciplinary Studies,
South China Normal University,
Guangzhou 510631, P. R. China.}
\email{\href{mailto:zhongxuexiu1989@163.com}{zhongxuexiu1989@163.com},
\href{zhongxuexiu1989@hotmail.com}{zhongxuexiu1989@hotmail.com}}

\thanks{The research was supported  by the Natural
Science Foundation of China  (No. 11801581, 12271184, 12271196,11931012, 12061012), Young Top-notch Talent Project of Guangdong Province (2024TQ08A725), Guangdong Basic and Applied Basic Research Foundation (2021A1515010034),Guangzhou Basic and Applied Basic Research Foundation(2024A04J10001).  The second author also want to express his acknowledgement to the support of the special foundation for Guangxi Ba Gui Scholars and Guangxi Province Talent Project (Financial support from Science and Technology Department).}

\subjclass[2000]{}

\maketitle
\noindent
{\small
\begin{Abstract}
In this paper, our aim is to prove the existence of normalized ground state for  the following Schr\"odinger systems with potentials
$$\begin{cases}
-\Delta u_1+V_1(x)u_1+\lambda_1 u_1=\partial_1 G(u_1,u_2)\;\quad&\hbox{in}\;\mathbb{R}^N,\\
-\Delta u_2+V_2(x)u_2+\lambda_2 u_2=\partial_2G(u_1,u_2)\;\quad&\hbox{in}\;\mathbb{R}^N,\\
0<u_1,u_2\in H^1(\mathbb{R}^N), N\geq 1,\\
\int_{\mathbb{R}^N}u_1^2 \mathrm{d} x=a_1, \int_{\mathbb{R}^N}u_2^2 \mathrm{d} x=a_2.
\end{cases}$$
 The potentials $V_1(x),V_2(x)$ are general such that $\inf \text{ess}~\sigma(-\Delta+V_\iota)>-\infty$, which are allowed to be singular at some points. And the nonlinearities $G(u_1,u_2)$ are considered of the form
$$
\begin{cases}
G(u_1, u_2):=\sum_{i=1}^{\ell}\frac{\mu_i}{p_i}|u_1|^{p_i}+\sum_{j=1}^{m}\frac{\nu_j}{q_j}|u_2|^{q_j}+\sum_{k=1}^{n}\beta_k |u_1|^{r_{1,k}}|u_2|^{r_{2,k}},~~\ell,m,n\in \mathbb{N}^+_0,\\
 \mu_i, \nu_j,\beta_k>0, ~2<r_{1,k}+r_{2,k},  p_i, q_j<2+\frac{4}{N}, ~r_{1,k}, r_{2,k}>1, \\
 i=1,2,\cdots, \ell; j=1,2,\cdots, m; k=1,2,\cdots, n.
 \end{cases}
 $$
 Under the mass sub-critical assumption, the normalized ground states are obtained as the minimum of the functional $J$ on the manifold $S_{a_1,a_2}$. Since the functional is not weak lower semi-continuous, to prove the minimizing problem is achievable, the key step is establishing the strict sub-additive inequality.
 Among its main ingredients is the study of the sharp decay of the positive solutions and the interaction estimates.
\end{Abstract}\\
\begin{key words}
Nonlinear Schr\"odinger system; Ground state normalized solution; Potentials; Interaction estimates.
\end{key words}
}

\indent


\s{Introduction}
\renewcommand{\theequation}{1.\arabic{equation}}
The time-dependent system of coupled nonlinear  Schr\"odinger equations
\begin{equation}\label{eq:system-t}
\begin{cases}
-i\frac{\partial}{\partial t}\Phi_1=\Delta \Phi_1-V_1(x)\Phi_1+\partial_1F(|\Phi_1|^2,|\Phi_2|^2)\Phi_1,\\
-i\frac{\partial}{\partial t}\Phi_2=\Delta \Phi_2-V_2(x)\Phi_1+\partial_2F(|\Phi_1|^2,|\Phi_2|^2)\Phi_2,\\
\Phi_\iota=\Phi_\iota(x, t)\in \mathbb{C}, \iota=1,2,
\end{cases}
(x,t)\in \mathbb{R}^N\times \mathbb{R}
\end{equation}
appear  in mean-field models for binary mixtures of Bose-Einstein condensates (BEC) or
models for binary mixtures of ultracold quantum gases of fermion atoms, such as Bose-Fermi mixtures and Fermi-Fermi mixtures, see \cite{Adhikari2007,Bagnato2015,Esry1997,Malomed2008} and the references therein. Topological objects are interesting topics in many fields of physics. Most of the current work and the existence of topological structures are qualitative discussions, but rare are constructing corresponding solutions, see \cite{Harland-2011,Zheng-2010,Marcone-2007}. These solutions can help people analyze the structure and properties of topological objects in detail. The current cold atom system provides an ideal platform to study complex topological structures. As a representative, binary mixtures of BEC has been widely used in experiments, see \cite{Kawaguchi-2010}.

Consider the solitary wave solutions of \eqref{eq:system-t}, i.e., $\Phi_\iota(x,t)=e^{i\lambda_\iota t}u_\iota(x), \iota=1,2$, it equivalents to the following elliptic system,
\begin{equation}\label{eq:system}
\begin{cases}
-\Delta u_1+V_1(x)u_1+\lambda_1 u_1=\partial_1 G(u_1,u_2)\;\quad&\hbox{in}\;\mathbb{R}^N,\\
-\Delta u_2+V_2(x)u_2+\lambda_2 u_2=\partial_2G(u_1,u_2)\;\quad&\hbox{in}\;\mathbb{R}^N,\\
u_1,u_2\in H^1(\mathbb{R}^N), N\geq 1,
\end{cases}
\end{equation}
with $G(s,t)=\frac{1}{2}F(s^2,t^2)$.
An important and  well-known  feature of \eqref{eq:system-t} is conservation of masses: the $L^2$-norms $\|\Phi_1(\cdot,t)\|_2,$ $\|\Phi_2(\cdot,t)\|_2$ of solutions are independent of $t\in \mathbb{R}$. Then it is reasonable to require the solutions satisfying the constraint
\begin{equation}\label{eq:mass-constraint}
(u_1,u_2)\in \widetilde S_{a_1,a_2}:=\left\{(u_1,u_2)\in \widetilde{\mathcal{H}} : \int_{\mathbb{R}^N} u_\iota^2 \mathrm{d} x=a_\iota, \iota=1,2\right\},
\end{equation}
where
$$ \widetilde{\mathcal{H}}:=\left\{(u_1,u_2)\in H^1(\mathbb{R}^N)\times H^1(\mathbb{R}^N): \int_{\mathbb{R}^N} V_\iota (x) u_\iota^2 \mathrm{d} x <\infty, \iota=1,2\right\}.$$
To get a solution having the prescribed masses, it is natural to find critical points of the energy functional $J$ constrained on $\widetilde S_{a_1,a_2}$, where
 \begin{equation}\label{eq:def-J}
J[u_1,u_2]=\frac{1}{2}\sum_{\iota=1}^{2}\int_{\mathbb{R}^N} \left[|\nabla u_\iota|^2 +V_\iota(x)u_\iota^2 \right]\mathrm{d} x-\int_{\mathbb{R}^N}G(u_1,u_2)\mathrm{d} x,
\end{equation}
In such a approach, the parameters $\lambda_1$ and $\lambda_2$, which are the frequencies of the waves, are unknown and will appear as Lagrange multipliers.

 Consider the case of $V(x)\equiv 0$ and $G(u,v)=\frac{\mu_1}{p}|u|^p +\beta |u|^{r_1} |v|^{r_2} +\frac{\mu_2}{q}|v|^q$, that is,
\begin{equation}\label{101.114}
\begin{cases}
-\Delta u +\lambda_1 u= \mu_1|u|^{p-2}u+ \beta r_1|u|^{r_1-2}|v|^{r_2} u \quad &\hbox{in}\;\mathbb{R}^N, \\
-\Delta v + \lambda_2 v= \mu_2|v|^{q-2}v+\beta r_2|u|^{r_1}|v|^{r_2-2} v \quad &\hbox{in}\;\mathbb{R}^N, \\
\int_{\mathbb{R}^N}u^2 =a, \int_{\mathbb{R}^N}v^2 =b,
\end{cases}
\end{equation}
For the mass sub-critical case, i.e., $p,q,r_1+r_2<2+\frac{4}{N}$, the functional $J$ is coercive and thus bounded from below when constrained on $\widetilde S_{a,b}$.
So, it is possible to find a global minimizer after solving some compactness issues. Furthermore, the set of global minimizers is orbitally stable. In that direction, there had been a good amount of works, directly on \eqref{101.114}, or on related problems, see \cite{Bhattarai2016,Cao2011,Nguyen2016,Ohta1996,gou2016existence} and references therein.

For the pure mass super-critical case, i.e., $p,q,r_1+r_2>2+\frac{4}{N}$, the functional $J$ is not coercive and unbounded from below when constrained on $\widetilde S_{a,b}$. So, there is no global minimizer and one has to search for another kind critical points such as mountain pass type. The case of $N=3, p=q=4, r_1=r_2=2$ is the most important model in physical applications:
\begin{equation}\label{101.113}
\begin{cases}
-\Delta u +\lambda_1 u= \mu_1u^3+\beta uv^2 \quad &\hbox{in}\;\mathbb{R}^3, \\
-\Delta v + \lambda_2 v= \mu_2v^3+\beta u^2v \quad &\hbox{in}\;\mathbb{R}^3, \\
\int_{\mathbb{R}^N}u^2 =a, \int_{\mathbb{R}^N}v^2 =b.
\end{cases}
\end{equation}
Bartsch, Soave \cite{bartsch2017a}  studied the repulsive case, i.e., $\beta<0$, and obtained a positive normalized solutions for all $a,b>0$. They also proved the phase-separation as $\beta\rightarrow -\infty$ there. They also studied in \cite{bartsch2019multiple}the multiplicity of normalized solutions by considering $a=b$.
The case of $\beta>0$ is tough which is studied by Bartsch, Jeanjean, Soave \cite{bartsch2016normalized}. They managed to find solutions to \eqref{101.113} provided $\beta\in (0,\beta_1(a, b))\cup (\beta_2(a,b),+\infty)$. Here, for the technical reason, the values $\beta_1(a,b)$ and $\beta_2(a,b)$ depend on $a,b$ heavily that $\beta_1(a,b)\rightarrow 0$ as $\frac{a}{b}\rightarrow 0$ or $\frac{a}{b}\rightarrow +\infty$, while $\beta_{2}(a,b)\rightarrow +\infty$ as $\frac{a}{b}\rightarrow 0$ or $\frac{a}{b}\rightarrow +\infty$. In particular, there is no value of $\beta$ so that the results from \cite{bartsch2016normalized} yield a solution for all masses. Hence, like to the Sirakov's open problem for the fixed frequency problem, see \cite{Sirakov2007,WeiZhongZou-2022}, it is natural to ask what is the best range of the parameters for the existence of positive normalized (ground state) solution ? In particular, if there is some $\beta>0$ such that the existence result holds for all masses? We call this kind of problem Bartsch-Jeanjean-Soave's open problem, see \cite[Remark 1.3-(a) and (d)]{bartsch2016normalized}. Later in \cite{bartschzhongzou2020}, Bartsch, Zhong and Zou  presented  a new approach based on the fixed point index in cones, bifurcation theory and the continuation method, and  obtained the existence of normalized solution for $\beta >0$ belonging  to much  better ranges, which is independent of the prescribed masses $a$ and $b$. This kind of global branch approach is also developed by Jeanjean, Zhang and Zhong\cite{JeanZhangZhong2024} to study the scale equation for the mass sub-critical, critical, super-critical in a unified way. And in the recent work\cite{JeanZhangZhong2024b}, Jeanjean et al. have made some contributions to the general Bartsch-Jeanjean-Soave's open problem but Sobolev subcritical case for $1\leq N\leq 4$. And Zhong\cite{Zhong2024} extends the results to the Sobolev critical case.

The mass mixed case, i.e., $\min\{p,q,r_1+r_2\}<2+\frac{4}{N}<\max\{p,q,r_1+r_2\}$, are considered in \cite{GouJeanjean2018,BartschJeanjean2018}. When the parameters are involving the Sobolev critical exponent, we refer to \cite{Bartsch2023,Radulescu2024a,GuoHeShuaiZhong2024,Zhong2024}.

It is worth mentioning that for the study of mass sub-critical case, Shibata \cite{Shibata2017} introduced a new rearrangement inequality, which is useful to derive the compactness of a minimizing sequences for a minimizing problem of the form
\begin{equation}
\tilde{E}_{a_1,a_2}:=\inf_{(u_1,u_2)\in S_{a_1,a_2}} \tilde{J}[u_1,u_2]
\end{equation}
where
$$
 S_{a_1,a_2}:=\left\{(u_1,u_2)\in {\mathcal{H}=H^1(\mathbb{R}^N)\times H^1(\mathbb{R}^N)} : \int_{\mathbb{R}^N} u_\iota^2 \mathrm{d} x=a_\iota, \iota=1,2\right\}
$$
and
$$\tilde{J}[u_1,u_2]=\frac{1}{2}\sum_{\iota=1}^{2}\|\nabla u_\iota\|_2^2-\int_{\mathbb{R}^N}G(|u_1|^2, |u_2|^2)\mathrm{d} x.$$

When the system involves potentials, in particular for the non-radial case, it will become very complicated in mathematical processing. Let us consider  the following system
\begin{equation}\label{eq:20210414-system}
\begin{cases}
-\Delta u_1+V_1(x)u_1+\lambda_1 u_1=\mu_1|u_1|^{p_1-2}u_1+ \beta r_1|u_1|^{r_1-2}|u_2|^{r_2} u_1\;\quad&\hbox{in}\;\mathbb{R}^N,\\
-\Delta u_2+V_2(x)u_2+\lambda_2 u_2=\nu_1|u_2|^{q_1-2}u_2+\beta r_2 |u_1|^{r_1}|u_2|^{r_2-2} u_2 \;\quad&\hbox{in}\;\mathbb{R}^N,\\
\int_{\mathbb{R}^N}u_\iota^2 =a_\iota>0, ~\iota=1,2, N\geq 1,
\end{cases}
\end{equation}
as we know, the first result is given by Guo, Li, Wei and Zeng   \cite{glwz}, where the authors considered the mass critical case that $N=2,p_1=q_1=4,r_1=r_2=2$. They assumed that $0\leq V_i\in C^\alpha_{\text{loc}}(\mathbb{R}^2)~(i=1,2)$ satisfying
$$\lim\limits_{|x|\to \infty}V_i(x)=\infty,~\hbox{both}~\inf\limits_{x\in \mathbb{R}^2}V_i(x)=0~\hbox{and}~\inf\limits_{x\in \mathbb{R}^2}(V_1(x)+V_2(x))~\hbox{are~attained},$$
and showed the existence of normalized solutions to \eqref{eq:20210414-system} provided  $\mu_1, \nu_1, \beta>0$ are appropriately small. In their argument, the coercive condition plays an crucial role in the $L^2$-compactness.
Ikoma and Miyamoto \cite{Ikoma2021} made a  progress for mass sub-critical case with  $r_1=r_2$. They discussed  the compactness of minimizing sequence for energy functional $J$ by assuming that $V_i\in C(\mathbb{R}^N), V_i(x)\leq \lim\limits_{|x|\to \infty}V_i(x)=:V_{i,\infty}$ $(i=1, 2)$ and some other suitable conditions.

\begin{remark}\label{remark:20230301-r1}
\begin{itemize}
\item[(i)]
From the physical view, $V_\iota(x)$ stands for the external potential well in BEC, which is usually radial function, i.e., $V_\iota=V_{\iota}(|x|)$. Furthermore, $V_\iota(x)$ can be written by $V_\iota=V_{\iota,opt}+V_{\iota,mag}$, where $V_{\iota,opt}$ is the optical trap and $V_{\iota,mag}$ represents the energy shift of particle spin in the external magnetic field.
\item[(ii)]In present paper, from a mathematical point of view, the potentials considered here are much more general, which are allowed to be singular at some points and not necessary of radial.
\end{itemize}
\end{remark}

Comparing to \cite{Ikoma2021}, we also aim to make contributions to get rid of the  constraint that $r_1=r_2$.
Precisely, we are concerned with the existence of real numbers $\lambda_1,\lambda_2\in \mathbb{R}$ and $u_1,u_2\in H^1(\mathbb{R}^N)$ that solve the following system:
\begin{equation}\label{eq:system-general}
\begin{cases}
-\Delta u_1+V_1(x)u_1+\lambda_1 u_1=\partial_1 G(u_1,u_2)\;\quad&\hbox{in}\;\mathbb{R}^N,\\
-\Delta u_2+V_2(x)u_2+\lambda_2 u_2=\partial_2G(u_1,u_2)\;\quad&\hbox{in}\;\mathbb{R}^N,\\
\int_{\mathbb{R}^N}u_1^2 \mathrm{d} x=a_1>0, \int_{\mathbb{R}^N}u_2^2 \mathrm{d} x=a_2>0,\\
u_1,u_2\in H^1(\mathbb{R}^N), N\geq 1,
\end{cases}
\end{equation}
where
\begin{equation}\label{eq:G-definition}
\begin{cases}
\displaystyle G(u_1, u_2):=\sum_{i=1}^{\ell}\frac{\mu_i}{p_i}|u_1|^{p_i}+\sum_{j=1}^{m}\frac{\nu_j}{q_j}|u_2|^{q_j}+\sum_{k=1}^{n}\beta_k |u_1|^{r_{1,k}}|u_2|^{r_{2,k}},~~\ell,m,n\in \mathbb{N}^+_0,\\
 \mu_i, \nu_j,\beta_k>0, ~2+\frac{4}{N}>r_{1,k}+r_{2,k},  p_i, q_j>2, ~r_{1,k}, r_{2,k}>1, \\
 i=1,2,\cdots, \ell; j=1,2,\cdots, m; k=1,2,\cdots, n.
 \end{cases}
\end{equation}
and $V_1(x), \ V_2(x)$ satisfy the following assumption:
\begin{itemize}
\item[(VH1)] For $\iota=1,2$, $V_{\iota}(x)\leq 0$ in $\mathbb{R}^N$ and
$\displaystyle\lim_{|x|\rightarrow \infty}V_\iota(x)=0$.
Furthermore, there exist  some $\sigma_\iota\in [0,1)$ and $\tau_\iota\geq 0$ such that
$$\|u\|_{V_\iota}^{2}\leq \sigma_\iota \|\nabla u\|_2^2 +\tau_\iota \|u\|_2^2, \ \  \forall u\in H^1(\mathbb{R}^N)\cap L^2\left(\mathbb{R}^N, -V_\iota(x)\mathrm{d} x\right), \iota=1,2,$$
\end{itemize}
where
\begin{equation}\label{eq:norm-uV}
\|u\|_{V}:=\left(\int_{\mathbb{R}^N}-V(x)|u|^2\mathrm{d} x\right)^{\frac{1}{2}}.
\end{equation}

Under the assumption of (VH1), we remark that
$$\langle u,v\rangle_{V_\iota}:=\int_{\mathbb{R}^N} -V_\iota uv \mathrm{d} x,~ \forall u,v\in H^1(\mathbb{R}^N),$$
define an inner product. And that the norm $\|u\|_{V_\iota}$ is induced by $\sqrt{\langle u,u\rangle_{V_\iota}}$.

\begin{remark}\label{remark:potential}
\begin{itemize}
\item[(i)]
We note that for the bounded potentials $V_\iota(x)\leq 0$ with $-\tau_\iota:=\inf_{x\in\mathbb{R}^N}V_\iota (x)>-\infty$, then
$\displaystyle\|u\|_{V_\iota}^{2}\leq \tau_\iota \|u\|_2^2$, and thus (VH1) holds naturally.
Another important application for unbounded potential is that $V_\iota (x)=-\frac{\kappa_\iota }{|x|^{s_\iota}}$ with  $s_\iota\in (0,2), \kappa_\iota >0$ or $s_\iota=2, \kappa_\iota\in (0, \frac{(N-2)^2}{4})$ for $N\geq 3$. In particular, the Coulomb potential is within our consideration, which is very important in physics. On the other hand,  we emphasize that our argument is indeed valid for general $G$ of the form
$$G(u_1,u_2)=\frac{\mu_1}{p_1}|u_1|^{p_1}+\frac{\nu_1}{q_1}|u_2|^{q_1}+\beta_1 |u_1|^{r_{1,1}}|u_2|^{r_{2,1}}+H(u_1,u_2)$$
with some suitable assumptions on $H(u_1,u_2)$ and $r_{1, 1}\neq r_{2, 1}$ is allowed.
\item[(ii)]In present paper, the nonlinearities are of mass sub-critical, and thus are of Sobolev sub-critical, so under the assumption (VH1), it is standard to see that a PS sequence of $J\big|_{\tilde{S}_{a_1,a_2}}$ can not concentrate at any point, even for the singular case of $V_\iota(x)=-\frac{\kappa_\iota}{|x|^2}$ with $\kappa_\iota>0$ small.  Indeed, suppose that $u_n\rightharpoonup u$ in $H^1(\mathbb{R}^N)$ and let $|\nabla (u_n-u)|^2\rightharpoonup \mu, -V_\iota (u_n-u)^2 \rightharpoonup \bar{\mu}$ in $\mathcal{M}(\mathbb{R}^N)$. Then it is easy to prove that $\mu=\bar{\mu}$ due to the sub-critical growth of the nonlinearities. Under the assumption of (VH1), we also have that $\bar{\mu}\leq \sigma_{\iota} \mu$. So by $\sigma_\iota<1$, we obtain that $\bar{\mu}=\mu=0$. Combining with $V_\iota(x)\rightarrow 0$ as $|x|\rightarrow +\infty$, one can see that $\|u_n-u\|_{V_\iota}^2=o(1)$. This observation will be used in the proof of compactness, see for example \eqref{eq:20221216-e1}.
\item[(iii)]The normalized solution to the scalar case
\begin{equation}\label{eq:V-pontential}
-\Delta u+(V(x)+\lambda)u=g(u)\;\hbox{in}~\mathbb{R}^N, \int_{\mathbb{R}^N}u^2 \mathrm{d} x=a>0,
\end{equation}
also attracts a lot of researchers, see e.g.
\cite{bartsch2003multi, Bartsch2021,  Ikoma-Miyamoto2020,  ikomatanaka2019,   Molle2022, Jeanjean1997, sslu2020cvpde, Mederski2020, PellacciPistoiaVairaVerzini, soave2020,soave2020cri, zhongzou2020} and the references therein. Verzini and Yu \cite{VerziniYu-2026} investigated the case with a weakly attractive potential $V$ and small mass, and proved the existence and multiplicity of positive normalized solutions for a Sobolev critical power nonlinearity, under suitable assumptions concerning $V$. Very recently, Li, Yu, Wang and Liu \cite{Li-2026} also studied normalized solutions for Schr\"odinger equations with potential in the Sobolev supercritical case.
\end{itemize}
\end{remark}

Before stating our results, we give a definition.
\begin{definition}\label{def:ground-state}
For given $a_1>0,a_2>0$, a solution $(\lambda_1, \lambda_2, u_1,u_2)$ to \eqref{eq:system-general} is called normalized ground state, or least energy normalized solution, if
$$J[u_1,u_2]=\min\left\{J[v_1,v_2]: (\hat{\lambda}_1,\hat{\lambda}_2, v_1,v_2 )~\hbox{ solves \eqref{eq:system-general} }\right\}.$$
\end{definition}

Firstly, we are concerned with the radial potential case that
\begin{itemize}
\item[(VH2)] $V_\iota(x)=V_\iota(|x|), \iota=1,2$ are non-decreasing functions with respect to $r=|x|$.
\end{itemize}
Our first main result is stated as following:
\begin{theorem}\label{thm:main-radial}
Suppose that $(VH1)-(VH2)$ hold and $G(s,t)$ is given by \eqref{eq:G-definition}.
Then problem \eqref{eq:system-general} has a ground state normalized solution $(\lambda_1,\lambda_2,u_1,u_2)$ with $\lambda_1\geq 0,\lambda_2\geq 0$ and $(u_1, u_2)\in \mathcal{H}^{rad}(\mathbb{R}^N)$, where $\mathcal{H}^{rad}\subset \mathcal{H}$ denotes the radial subspace.
\end{theorem}

When the potential $V_\iota(x)$ is not of radial, we can not work in the radial subspace as \cite{BartschJeanjean2018,gou2016existence}, etc.  In order to state our next main theorem, we suppose that
\begin{equation}\label{eq:parameters-paixu}
p_1\leq p_i, q_1\leq q_j, r_{1,1}\leq r_{1,k}, r_{2,1}\leq r_{2,k}\;\hbox{for} ~i\geq 1, j\geq 1,k\geq 1.
\end{equation}
For the technical reason (see Remark \ref{remark:lambda=0} below),
we require the following assumption in addition.
\begin{equation}\label{eq:parameters-L}
\hbox{When} ~N\geq 5, \hbox{we suppose further that}~ p_1,q_1\in \left(2, 2+\frac{2}{N-2}\right].
\end{equation}

Here comes our another main theorem focused on the case that $V_\iota(x)$ is not necessary radial.
\begin{theorem}\label{thm:main-non-radial}
Suppose that $(VH1)$ holds and $G(s,t)$ satisfies \eqref{eq:G-definition},\eqref{eq:parameters-paixu} and \eqref{eq:parameters-L}. Then for any $a_1>0,a_2>0$, Problem \eqref{eq:system-general} has a ground state normalized solution $(\lambda_1,\lambda_2,u_1,u_2)$ with $\lambda_1> 0,\lambda_2> 0$, $(u_1, u_2)\in \mathcal{H}$ and $u_1 >0,  \ u_2>0$ in $\mathbb{R}^N$.
\end{theorem}

\begin{remark}\label{remark:lambda=0}
For these Lagrange multipliers $\lambda_1$ and $\lambda_2$, from a physical point of view it represents the chemical potentials of standing waves. We point out that there are situations in the Bose-Einstein condensate theory that requires the chemical potentials are positive, see \cite{LiebSeiringerSolovejYngvason2005,Pitaevskii2003}. By a well-known result own to Frank H. Clarke(see \cite[Theorem 1]{Clarke1976}), we always have that $\lambda_1\geq 0$ and $\lambda_2\geq 0$, see also Proposition \ref{prop:minimizer} in Section \ref{sec:properties-Ca}. When $V_\iota(x)$ is not of radial, to prove the strict sub-additive inequality, we need a sharp decay estimation of the positive solution. And that the positiveness of the Lagrange multipliers plays an crucial role in proving the $L^2$-compactness. Hence, we need to exclude the case of $\lambda_\iota=0$ $(\iota =1,\ 2)$. Just because of this, we need to use the  known Liouville type result.
Observing  that $\partial_\iota G(u_1,u_2)\geq 0$ in $\mathbb{R}^N$, $\iota=1,2$. So by the known Liouville type results, one can see that $\lambda_\iota>0$ provided $1\leq N\leq 4$. From a mathematical point of view, it seems that $\lambda_\iota=0$ may happen for high dimension case without any other condition. However, if suppose $p_1,q_1\in (2,2+\frac{2}{N-2}]$ in addition to the assumptions, one can also exclude the case of $\lambda_\iota=0$, see for instance \cite[Lemma A.2]{Ikoma2014}. For our case, see Lemma \ref{lemma:20210414-l1} in Section \ref{sec:compactness-proof-Th2}.
\end{remark}

For the mass sub-critical problem \eqref{eq:system-general}, the normalized ground state will be obtained as the minimum of the functional $J$ on the manifold $S_{a_1,a_2}$.
When $V(x)$ is not radial, we can not work in the radial subspace. We emphasize that the functional $J$ is not weakly lower semi-continuous. To prove the minimum is achievable, the key step is proving the strict sub-additive inequality, due to the well known concentration compactness principle of Lions\cite{Lions1984a,Lions1984b}.
 Suppose that $(u_1,u_2)\neq (0,0)$ is a nonnegative  solution to
\eqref{eq:system-20210329-1}
and $(\omega_1, \omega_2)\neq (0,0)$ is a nonnegative  solution to
\eqref{eq:system-20210329-2},
 we shall construct an element $(\phi_1,\phi_2)\in\mathcal{H}$ (in Proposition \ref{prop:binding-inequality}) such that
\begin{equation}\label{09131}
(\|\phi_1\|_2^2,\|\phi_2\|_2^2)=(\|u_1\|_2^2+\|\omega_1\|_2^2,\|u_2\|_2^2+\|\omega_2\|_2^2)~\hbox{and}
~J[\phi_1,\phi_2]<J[u_1,u_2]+I[\omega_1,\omega_2].
\end{equation}

In the process of this construction, sharp decay estimations of the positive solutions to \eqref{eq:system-20210329-1}  and \eqref{eq:system-20210329-2} will be helpful (See Proposition \ref{prop:sharp-estimation}) . And the key step is the interaction estimation for the couples terms, see Lemma \ref{lemma:intersection-estimation}.
Comparing with \cite{Ikoma2021}, our nonlinearities  are more general and  both $V_1(x)$ and $V_2(x)$  can be allowed to be  singular at some points, which brings
much more difficulties to us in mathematical processing. Since the proofs are too long and complicated,  we prefer to give these proofs of  Proposition \ref{prop:sharp-estimation}, Lemma \ref{lemma:intersection-estimation} and Proposition \ref{prop:binding-inequality}  in Appendixes-A,B,C respectively.

The paper is organized as follows. We collect and prove a few basic facts about the minimizing sequence related to \eqref{eq:system-general} in Section \ref{sec:prelim}.
In Section \ref{sec:properties-Ca}, we shall prove some properties of $C_{a_1,a_2}$, such as the monotonicity, continuity and sub-additive inequality. In Section \ref{radial-case}, we shall focus on the radial case and prove Theorem  \ref{thm:main-radial} via the rearrangement argument.
For the non-radial case, to establish the strict sub-additive inequality, a better decay estimation of the positive solutions  is necessary. So in Section \ref{sec:decay-estimation}, we shall give a sharp decay estimations of the positive solutions for the {\it fixed frequency problem}.  Basing on these results, we can prove the strict sub-additive inequality in Section  \ref{sec:interaction-estimate} (see Proposition \ref{prop:binding-inequality}). Finally in Section \ref{sec:compactness-proof-Th2}, we apply concentration compactness argument to prove Theorem \ref{thm:main-non-radial}.

Throughout the paper we use the notation $\|u\|_p$ to denote  the $L^p$-norm, and simply write $\mathcal{H}=H^1(\mathbb{R}^N)\times H^1(\mathbb{R}^N)$. $\mathcal{H}^{rad}$ denotes the radial subspace. The notation $\rightharpoonup$ denotes weak convergence in $H^1(\mathbb{R}^N)$ or $\mathcal{H}$. Capital latter $C$ stands for positive constant which may depend on some parameters, whose precise value can change from line to line.
\medskip


\s{Some Preliminaries }\label{sec:prelim}
\renewcommand{\theequation}{2.\arabic{equation}}
In this section, we summarize several results that will be used in the rest discussion.
\begin{lemma}\label{lemma:20210331-l1}
Suppose that $G(s,t)$ is given by \eqref{eq:G-definition}.
Then for any $\varepsilon>0$, there exists some $C_\varepsilon>0$ such that
\begin{equation}\label{eq:20210331-e1}
G(s,t)\leq C_\varepsilon (s^2+t^2)+\varepsilon (s^{2+\frac{4}{N}}+t^{2+\frac{4}{N}}),\forall (s,t)\in \mathbb{R}_+^2.
\end{equation}
\end{lemma}
\begin{proof}
It follows from Young inequality that
$$s^{r_{1,k}}t^{r_{2,k}}\leq \frac{1}{r_{1,k}+r_{2,k}} r_{1,k}^{\frac{r_{1,k}}{r_{1,k}+r_{2,k}}} r_{2,k}^{\frac{r_{2,k}}{r_{1,k}+r_{2,k}}} \left(s^{r_{1,k}+r_{2,k}}+t^{r_{1,k}+r_{2,k}}\right).$$
On the other hand, let $p\in (2, 2+\frac{4}{N})$,  for any $\varepsilon>0$, there exists some $C_\varepsilon>0$  such that
\begin{equation}\label{eq:Young-inequality}
y^p\leq C_\varepsilon y^2+\varepsilon y^{2+\frac{4}{N}}, \ \  \forall y\geq 0.
\end{equation}
It is easy to check that \eqref{eq:20210331-e1} holds by using \eqref{eq:G-definition}.
\end{proof}

\begin{lemma}\label{lemma:20210310-xl1}
Suppose that $(VH1)$ holds and $G(s,t)$ is given by \eqref{eq:G-definition}. Let $\sigma:=\frac{1-\max\{\sigma_1,\sigma_2\}}{2}>0$  determined by $(VH1)$. Then for any $\varepsilon>0$, there exists some $C_\varepsilon>0$  such that
\begin{align*}
&J[u_1,u_2]\geq \sigma (\|\nabla u_1\|_2^2+\|\nabla u_2\|_2^2)-C_\varepsilon (\|u_1\|_2^2+\|u_2\|_2^2)\\
&-\varepsilon \left(\|\nabla u_{1}\|_2^2 \|u_{1}\|_{2}^{\frac{4}{N}}+\|\nabla u_{2}\|_2^2 \|u_{2}\|_{2}^{\frac{4}{N}}\right),\quad \forall ~(u_1,u_2)\in \mathcal{H}.
\end{align*}
\end{lemma}
\begin{proof}
By Lemma \ref{lemma:20210331-l1} and a direct computation.
\end{proof}

\begin{lemma}\label{lemma:J-bounded-blow}
Suppose that $(VH1)$ holds and $G(s,t)$ is given by \eqref{eq:G-definition}. For any $a_1\geq 0,a_2\geq 0$, $J\big|_{S_{a_1,a_2}}$ is bounded from below, i.e.,
$$C_{a_1,a_2}:=\inf_{(u_1,u_2)\in S_{a_1,a_2}}J[u_1,u_2]>-\infty.$$
\end{lemma}
\begin{proof}
This is a direct conclusion from Lemma \ref{lemma:20210310-xl1}.
\end{proof}

\begin{lemma}\label{lemma:bounded-bounded}
Suppose that $(VH1)$ holds and $G(s,t)$ is given by \eqref{eq:G-definition}.
Then any $L^2$-bounded sequence $\{(u_{1,\vartheta}, u_{2,\vartheta})\}$ with $\sup_{\vartheta\in \mathbb{N}}J[u_{1,\vartheta}, u_{2,\vartheta}]<\infty$ is also bounded in $\mathcal{H}$.
\end{lemma}
\begin{proof}
Let $\{(u_{1,\vartheta}, u_{2,\vartheta})\}$ be a $L^2$-bounded sequence, i.e., $\|u_{\iota,\vartheta}\|_2^2\leq M_\iota, \iota=1,2$. Suppose that
$J[u_{1,\vartheta}, u_{2,\vartheta}]\leq M_3,\;\forall ~\vartheta\in \mathbb{N}$.
 Recalling Lemma \ref{lemma:20210310-xl1}, for any $\varepsilon>0$, there exists some $C_\varepsilon>0$  such that
\begin{align*}
M_3\geq J[u_{1,\vartheta},u_{2,\vartheta}]\geq& \sigma (\|\nabla u_{1,\vartheta}\|_2^2+\|\nabla u_{2,\vartheta}\|_2^2)-C_\varepsilon (\|u_{1,\vartheta}\|_2^2+\|u_{2,\vartheta}\|_2^2)\\
&-\varepsilon \left(\|\nabla u_{1,\vartheta}\|_2^2 \|u_{1,\vartheta}\|_{2}^{\frac{4}{N}}+\|\nabla u_{2,\vartheta}\|_2^2 \|u_{2,\vartheta}\|_{2}^{\frac{4}{N}}\right)\\
\geq&\sigma (\|\nabla u_{1,\vartheta}\|_2^2+\|\nabla u_{2,\vartheta}\|_2^2) - C_\varepsilon (M_1+M_2)\\
&-\varepsilon \left(M_{1}^{\frac{2}{N}}\|\nabla u_{1,\vartheta}\|_2^2 +M_{2}^{\frac{2}{N}}\|\nabla u_{2,\vartheta}\|_2^2 \right).
\end{align*}
We can take $\varepsilon>0$ small enough such that
$\varepsilon   M_{\iota}^{\frac{2}{N}}<\frac{1}{2}\sigma, \iota=1,2.$
Then it follows that
\begin{equation}\label{eq:bounded}
\sup_{\vartheta\in \mathbb{N}}\sum_{\iota=1}^{2}\|\nabla u_{\iota,\vartheta}\|_2^2<\infty,
\end{equation}
and this implies that the sequence $\{(u_{1,\vartheta}, u_{2,\vartheta})\}_{\vartheta=1}^{\infty}$ is  bounded in $\mathcal{H}$.
\end{proof}

\begin{remark}
\label{remark:existence-bdd-minizing} \ \
\begin{itemize}
\item[(i)]Suppose that $(VH1)$ holds and $G(s,t)$ is given by \eqref{eq:G-definition}. For any given $a_1>0,a_2>0$, by Lemma \ref{lemma:J-bounded-blow} and Ekeland's variational principle,  the existence of Palais-Smale minimizing sequence is trivial. Furthermore, by Lemma \ref{lemma:bounded-bounded}, this Palais-Smale minimizing sequence is also bounded in $\mathcal{H}$. So the main difficulty of mass sub-critical problem is the convergence of minimizing sequences.
\item[(ii)]Define
\begin{equation}\label{eq:functional-infinity}
I[u_1,u_2]:=\frac{1}{2}\sum_{\iota=1}^{2} \|\nabla u_{\iota}\|_2^2-\int_{\mathbb{R}^N}G(u_1,u_2)\mathrm{d} x, \forall (u_1, u_2)\in \mathcal{H},
\end{equation}
and
\begin{equation}
E_{a_1,a_2}:=\inf_{(u_1,u_2)\in S_{a_1,a_2}}I[u_1,u_2].
\end{equation}
It is equivalent to take $V_1(x)\equiv V_2(x)\equiv 0$ in $J[u_1,u_2]$ (which means that $\sigma_\iota=0,\tau_\iota=0$ and $\sigma:=\frac{1-\max\{\sigma_1,\sigma_2\}}{2}=\frac{1}{2}$). So we have a sequence of conclusions for $I[u_1,u_2]$ and $E_{a_1,a_2}$, as that in Lemmas \ref{lemma:20210310-xl1}-\ref{lemma:bounded-bounded}.
\end{itemize}
\end{remark}


\s{Some properties of $C_{a_1,a_2}$ and $E_{a_1,a_2}$ }\label{sec:properties-Ca}
\renewcommand{\theequation}{3.\arabic{equation}}
In this section, we shall establish some basic properties of $C_{a_1,a_2}$ and $E_{a_1,a_2}$ including monotonicity, continuity and sub-additive inequality.

Observing that $J[u_1,u_2]=I[u_1,u_2]-\frac{1}{2}\sum_{\iota=1}^{2}\|u_\iota\|_{V_\iota}^{2}$, so we have
$C_{a_1,a_2}\leq E_{a_1,a_2}$.
Recalling the fiber map
\begin{equation}\label{eq:ut-invariant}
u(x)\mapsto (t\star u)(x):=t^{\frac{N}{2}}u(tx),
\end{equation}
for $(t, u)\in \mathbb{R}^+\times H^1(\mathbb{R}^N)$, which preserves the $L^2$-norm, we have the following result.

\begin{lemma}\label{lemma:min-negative}
Suppose that $(VH1)$ holds and $G(s,t)$ is given by \eqref{eq:G-definition}. For any $a_1\geq 0,a_2\geq 0$ and $(a_1,a_2)\neq (0,0)$, then $C_{a_1,a_2}\leq E_{a_1,a_2}<0$. In particular, if $E_{a_1,a_2}$ is attained, and $V_\iota(x)\not\equiv 0, \iota=1,2$, then $C_{a_1,a_2}<E_{a_1,a_2}<0$.
\end{lemma}
\begin{proof}
Letting  $(u_1,u_2)\in S_{a_1,a_2}$, we remark that $t\star (u_1,u_2)\in S_{a_1,a_2}$ for any $t>0$. A direct computation shows that
\begin{equation}
I[t\star (u_1,u_2)]=\frac{1}{2}\sum_{\iota=1}^{2} \|\nabla u_{\iota}\|_2^2 t^2-t^{-N}\int_{\mathbb{R}^N}G(t^{\frac{N}{2}}u_1, t^{\frac{N}{2}}u_2)\mathrm{d} x.
\end{equation}
Noting that
$$G(t^{\frac{N}{2}}u_1, t^{\frac{N}{2}}u_2)=o(1)t^N(|u_1|^2+|u_2|^2),\;\hbox{as}\;t\rightarrow 0^+,$$
we have $\lim_{t\rightarrow 0^+}I[t\star (u_1,u_2)]=0$.

On the other hand, for $t>0$ small enough, by \eqref{eq:G-definition}, there exists some $\varepsilon_0>0$ such  that
\begin{align*}
\frac{I[t\star (u_1,u_2)]}{t^2}=&\frac{1}{2}\sum_{\iota=1}^{2} \|\nabla u_{\iota}\|_2^2 -t^{-N-2}\int_{\mathbb{R}^N}G(t^{\frac{N}{2}}u_1, t^{\frac{N}{2}}u_2)\mathrm{d} x\\
\leq & \frac{1}{2}\sum_{\iota=1}^{2} \|\nabla u_{\iota}\|_2^2 - t^{-N-2} \int_{\mathbb{R}^N}\frac{\mu_1}{p_1}|t^{\frac{N}{2}}u_{1}|^{p_{1}}+\frac{\nu_1}{q_1}|t^{\frac{N}{2}}u_{2}|^{q_{1}}\mathrm{d} x\\
=&\frac{1}{2}\sum_{\iota=1}^{2} \|\nabla u_{\iota}\|_2^2-\frac{\mu_1}{p_1}\|u_{1}\|_{p_{1}}^{p_{1}} t^{\frac{N}{2}p_{1}-N-2}- \frac{\nu_1}{q_1}\|u_{2}\|_{q_{1}}^{q_{1}} t^{\frac{N}{2}q_{1}-N-2}.
\end{align*}
 Noting that $p_1,q_1\in (2,2+\frac{4}{N})$, we have that $\frac{N}{2}p_{1}-N-2<0, \frac{N}{2}q_{1}-N-2<0$. By $(a_1,a_2)\neq (0,0)$, at least one of $\|u_1\|_{p_1}^{p_1}$ and $\|u_2\|_{q_1}^{q_1}$ is positive. Hence, there exists some $t_0>0$ such that
\begin{equation}
\frac{d}{dt}I[t\star (u_1,u_2)]<0,\;\forall t\in (0,t_0),
\end{equation}
which implies that $I[t\star (u_1,u_2)]<0,\;\forall t\in (0,t_0)$.
Recalling that $J[u_1,u_2]=I[u_1,u_2]-\frac{1}{2}\sum_{\iota=1}^{2}\|u_\iota\|_{V_\iota}^{2}\leq I[u_1,u_2]$, we have that
$$C_{a_1,a_2}:=\inf_{(u_1,u_2)\in S_{a_1,a_2}}J[u_1,u_2]\leq E_{a_1,a_2}:=\inf_{(u_1,u_2)\in S_{a_1,a_2}}I[u_1,u_2]<0.$$

In particular, if $E_{a_1,a_2}$ is attained, we can take $(0,0)\neq (u_1,u_2)\in S_{a_1,a_2}$ with $u_1\geq 0,u_2\geq 0$ such that $I[u_1,u_2]=E_{a_1,a_2}$. A standard elliptic estimation shows the solution is of classical. So by strong maximum principle, it is easy to see that at least one of $u_1$ and $u_2$ is  positive in $\mathbb{R}^N$.
Then by $V_\iota(x)\leq 0$ but $V_\iota(x)\not\equiv 0, \iota=1,2$, we obtain that
$$C_{a_1,a_2}\leq J[u_1,u_2]=I[u_1,u_2]-\frac{1}{2}\sum_{\iota=1}^{2}\|u_\iota\|_{V_\iota}^{2}<E_{a_1,a_2}<0.$$
\end{proof}

\begin{lemma}\label{lemma:sub-addictive-1}
For $(a_1,a_2), (b_1,b_2)\in \mathbb{R}_+^2$, we have
$$C_{a_1+b_1,a_2+b_2}\leq C_{a_1,a_2}+E_{b_1,b_2}\;\hbox{and}\;E_{a_1+b_1,a_2+b_2}\leq E_{a_1,a_2}+E_{b_1,b_2}.$$
\end{lemma}
\begin{proof}
We only prove $C_{a_1+b_1,a_2+b_2}\leq C_{a_1,a_2}+E_{b_1,b_2}$, and $E_{a_1+b_1,a_2+b_2}\leq E_{a_1,a_2}+E_{b_1,b_2}$ can be proved by taking $V_1(x)\equiv V_2(x)\equiv 0$.
If $(a_1,a_2)=(0,0)$ or $(b_1,b_2)=(0,0)$, it is trivial. Next, we suppose that $(a_1,a_2)\neq (0,0)$ and $(b_1,b_2)\neq (0,0)$.
For any $\varepsilon>0$, by density, we can take some   $(u_1, u_2)\in S_{a_1,a_2} , (v_1,v_2)\in S_{b_1,b_2}$ with $u_\iota,v_\iota\in C_0^\infty(\mathbb{R}^N), \iota=1,2$,   such that
$$C_{a_1,a_2}\leq J[u_1,u_2]<C_{a_1,a_2}+\frac{\varepsilon}{2}\;\hbox{and}\;E_{b_1,b_2}\leq I[v_1,v_2]<E_{b_1,b_2}+\frac{\varepsilon}{2}.$$
Let $(v_{1,R}, v_{2,R}):=(v_1(\cdot-R{\bf e_1}), v_2(\cdot-R{\bf e_1}))$.
Take some $R_0>0$ large enough such that
$$(\text{supp}~u_1\cup \text{supp}~u_2)\cap (\text{supp}~v_{1,R_0}\cup \text{supp}~v_{2,R_0})=\emptyset.$$
Then we have that
$(u_1+v_{1,R}, u_2+v_{2,R})\in S_{a_1+b_1,a_2+b_2}, \forall R\geq R_0,$
and
\begin{align*}
&J[u_1+v_{1,R}, u_2+v_{2,R}]=J[u_1,u_2]+J[v_{1,R},v_{2,R}]\leq J[u_1,u_2]+I[v_{1,R},v_{2,R}]\\
=&J[u_1,u_2]+I[v_1,v_2]\leq C_{a_1,a_2}+E_{b_1,b_2}+\varepsilon.
\end{align*}
By the arbitrary of $\varepsilon$, we obtain that
$C_{a_1+b_1,a_2+b_2}\leq C_{a_1,a_2}+E_{b_1,b_2}.$
\end{proof}

\begin{corollary}\label{cro:monotonicity-Ca}
For $(a_1,a_2), (b_1,b_2)\in \mathbb{R}_+^2$ with $(b_1,b_2)\neq (0,0)$, then
$C_{a_1+b_1,a_2+b_2}<C_{a_1,a_2}$ and $E_{a_1+b_1,a_2+b_2}<E_{a_1,a_2}$.
\end{corollary}
\begin{proof}

Recalling Lemma \ref{lemma:min-negative}, we have $E_{b_1,b_2}<0$ for $(b_1,b_2)\neq (0,0)$. Then by Lemma \ref{lemma:sub-addictive-1}, one can see that $C_{a_1+b_1,a_2+b_2}<C_{a_1,a_2}$ and $E_{a_1+b_1,a_2+b_2}<E_{a_1,a_2}$.
\end{proof}

\begin{lemma}\label{lemma:CE-continuous}
Suppose that $(VH1)$ holds and $G(s,t)$ is given by \eqref{eq:G-definition}.
Let  $\{(u_{1,\vartheta}, u_{2,\vartheta})\}$ be bounded in $\mathcal{H}$ with
$(\|u_{1,\vartheta}\|_2^2, \|u_{2,\vartheta}\|_2^2)\rightarrow (a_1,a_2)\in \mathbb{R}_+^2$.
Then we have
\begin{equation}\label{eq:20210412-xle1}
\liminf_{\vartheta\rightarrow \infty}J[u_{1,\vartheta}, u_{2,\vartheta}]\geq C_{a_1,a_2}
\end{equation}
and
\begin{equation}\label{eq:20210412-xle2}
\liminf_{\vartheta\rightarrow \infty}I[u_{1,\vartheta}, u_{2,\vartheta}]\geq E_{a_1,a_2}.
\end{equation}
\end{lemma}
\begin{proof}
We only prove \eqref{eq:20210412-xle1}.
If $a_1=a_2=0$, then we have that
$(u_{1,\vartheta}, u_{2,\vartheta})\rightarrow (0,0)$ in $L^2(\mathbb{R}^N)\times L^2(\mathbb{R}^N)$.
Since $\{(u_{1,\vartheta}, u_{2,\vartheta})\}$ is bounded in $\mathcal{H}$, by Gagliardo-Nirenberg inequality, one can see that
$u_{\iota,\vartheta}\rightarrow 0$ in $L^p(\mathbb{R}^N)$ for $2\leq p<2^*$ and $\iota=1,2$. Hence,
$$\int_{\mathbb{R}^N}G(u_{1,\vartheta}, u_{2,\vartheta})\mathrm{d} x\rightarrow 0~\hbox{as}~\vartheta\rightarrow \infty.$$
By $(VH1)$, we obtain that
$$\liminf_{\vartheta\rightarrow \infty}J[u_{1,\vartheta}, u_{2,\vartheta}]\geq 0=C_{0,0}.$$

If $a_1=0<a_2$, we have that
$u_{1,\vartheta}\rightarrow 0$ in $L^2(\mathbb{R}^N)$. Then by Gagliardo-Nirenberg inequality again, $u_{1,\vartheta}\rightarrow 0$ in $L^p(\mathbb{R}^N)$ for $2\leq p<2^*$. By mean value theorem, there exists $v_{1,\vartheta}=\theta u_{1,\vartheta}$ for
some $\theta(x)\in (0,1)$ such that
\begin{equation}\label{eq:20210413-e1}
\int_{\mathbb{R}^N}G(u_{1,\vartheta}, u_{2,\vartheta})\mathrm{d} x-\int_{\mathbb{R}^N}G(0, u_{2,\vartheta})\mathrm{d} x
=\int_{\mathbb{R}^N}\partial_1 G(\theta u_{1,\vartheta},u_{2,\vartheta}) u_{1,\vartheta} \mathrm{d} x\rightarrow 0.
\end{equation}
Hence,
$$\liminf_{\vartheta\rightarrow \infty}J[u_{1,\vartheta}, u_{2,\vartheta}]\geq \liminf_{\vartheta\rightarrow \infty}J[0, u_{2,\vartheta}].$$
By $\|u_{2,\vartheta}\|_2^2\rightarrow a_2>0$, without loss of generality, we may assume that $\|u_{2,\vartheta}\|_2>0$ for all $\vartheta\in \mathbb{N}$. Let $v_{2,\vartheta}:= \frac{\sqrt{a_2}}{\|u_{2,\vartheta}\|_2} u_{2,\vartheta}, \forall \vartheta\in \mathbb{N}$. Then we have that
\begin{equation}\label{eq:20210413-e2}
\|v_{2,\vartheta}\|_2^2=a_2, u_{2,\vartheta}=t_\vartheta v_{2,\vartheta}~\hbox{with}~t_\vartheta\rightarrow 1.
\end{equation}
Then we have
\begin{equation}\label{eq:20210413-e2}
\|\nabla u_{2,\vartheta}\|_2^2-\|u_{2,\vartheta}\|_{V_2}^{2}=t_\vartheta^2 [\|\nabla v_{2,\vartheta}\|_2^2-\|v_{2,\vartheta}\|_{V_2}^{2}]=(1+o(1))[\|\nabla v_{2,\vartheta}\|_2^2-\|v_{2,\vartheta}\|_{V_2}^{2}],
\end{equation}
and by mean value theorem, there exists some $\eta_\vartheta(x)\in [\min\{t_\vartheta, 1\}, \max \{t_\vartheta, 1\}]$ such that
\begin{equation}\label{eq:20210413-e3}
\int_{\mathbb{R}^N}G(0,u_{2,\vartheta})\mathrm{d} x-\int_{\mathbb{R}^N}G(0,v_{2,\vartheta})\mathrm{d} x
=(t_\vartheta -1)\int_{\mathbb{R}^N} \partial_2G(0,\eta_\vartheta v_{2,\vartheta}) v_{2,\vartheta} \mathrm{d} x\rightarrow 0.
\end{equation}
Hence,
\begin{align*}
\liminf_{\vartheta\rightarrow \infty} J[u_{1,\vartheta}, u_{2,\vartheta}]\geq &\liminf_{\vartheta\rightarrow \infty} J[0, v_{2,\vartheta}]\geq C_{0,a_2}.
\end{align*}

We remark that the case of $a_1>0=a_2$ can be proved similarly as above.

Finally, we consider the case of $a_1>0,a_2>0$.
We let
$$v_{\iota,\vartheta}=\frac{\sqrt{a_\iota}}{\|u_{\iota,\vartheta}\|_2} u_{\iota,\vartheta}, \ \ \text{and} \ \        t_{\iota,\vartheta}:=\frac{\|u_{\iota,\vartheta}\|_2}{\sqrt{a_\iota}}, \quad\iota=1,2.$$
Then we have that
$u_{\iota,\vartheta}=t_{\iota,\vartheta} v_{\iota,\vartheta}$ with
$\|v_{\iota,\vartheta}\|_2^2=a_\iota, \forall \vartheta\in \mathbb{N}N~\hbox{and}~t_{\iota,\vartheta}=1+o(1)~\hbox{as}~\vartheta\rightarrow \infty.$
Hence, similar to \eqref{eq:20210413-e2}, for $\iota=1,2$,  we also have
\begin{equation}\label{eq:20210413-e4}
\|\nabla u_{\iota,\vartheta}\|_2^2-\|u_{\iota,\vartheta}\|_{V_\iota}^{2}=(1+o(1))[\|\nabla v_{\iota,\vartheta}\|_2^2-\|v_{\iota,\vartheta}\|_{V_\iota}^{2}].
\end{equation}
On the other hand, there exist $\eta_{\iota,\vartheta}(x)\in [\min\{t_{\iota,\vartheta}, 1\}, \max \{t_{\iota,\vartheta}, 1\}]$ such that
\begin{align*}
\int_{\mathbb{R}^N}G(u_{1,\vartheta}, u_{2,\vartheta})\mathrm{d} x- \int_{\mathbb{R}^N}G(v_{1,\vartheta}, v_{2,\vartheta})\mathrm{d} x
=&\int_{\mathbb{R}^N}G(u_{1,\vartheta}, u_{2,\vartheta})\mathrm{d} x-\int_{\mathbb{R}^N}G(u_{1,\vartheta}, v_{2,\vartheta})\mathrm{d} x\\
&+\int_{\mathbb{R}^N}G(u_{1,\vartheta}, v_{2,\vartheta})\mathrm{d} x-\int_{\mathbb{R}^N}G(v_{1,\vartheta}, v_{2,\vartheta})\mathrm{d} x\\
=&(t_{2,\vartheta}-1)\int_{\mathbb{R}^N} \partial_2G(u_{1,\vartheta}, \eta_{2,\vartheta} v_{2,\vartheta}) v_{2,\vartheta} \mathrm{d} x\\
&+(t_{1,\vartheta}-1)\int_{\mathbb{R}^N} \partial_2G(\eta_{1,\vartheta}v_{1,\vartheta},  v_{2,\vartheta}) v_{1,\vartheta} \mathrm{d} x\\
=&o(1).
\end{align*}
We obtain that
\begin{align*}
\liminf_{\vartheta\rightarrow \infty} J[u_{1,\vartheta}, u_{2,\vartheta}]\geq &\liminf_{\vartheta\rightarrow \infty} J[v_{1,\vartheta}, v_{2,\vartheta}]\geq C_{a_1,a_2}.
\end{align*}
\end{proof}

\begin{corollary}\label{cro:20210413-c1}
Suppose that $(VH1)$ holds and $G(s,t)$ is given by \eqref{eq:G-definition}. Then $C_{a_1,a_2}$ and $E_{a_1,a_2}$ is continuous with respect to $(a_1,a_2)\in \mathbb{R}_+^2$.
\end{corollary}
\begin{proof}
We only prove the continuity of $C_{a_1,a_2}$. Let $(a_1,a_2)\in \mathbb{R}_+^2$.  For any sequence $(a_{1,\vartheta}, a_{2,\vartheta})\in (\mathbb{R}_+^2)$ with $(a_{1,\vartheta}, a_{2,\vartheta})\rightarrow (a_1,a_2)$, we shall prove that
\begin{equation}\label{eq:20210413-e5}
\lim_{\vartheta\rightarrow \infty} C_{a_{1,\vartheta}, a_{2,\vartheta}}=C_{a_1,a_2}.
\end{equation}
Firstly, $\forall ~\varepsilon>0, \exists~(u_{1,\vartheta},u_{2,\vartheta})\in S_{a_{1,\vartheta}, a_{2,\vartheta}}$ such that
$C_{a_{1,\vartheta}, a_{2,\vartheta}}\leq J[u_{1,\vartheta},u_{2,\vartheta}]\leq C_{a_{1,\vartheta}, a_{2,\vartheta}}+\varepsilon.$
Noting that $C_{a_{1,\vartheta}, a_{2,\vartheta}}\leq 0$, by Lemma \ref{lemma:bounded-bounded},  $(u_{1,\vartheta},u_{2,\vartheta})$ is bounded in $\mathcal{H}$. Then by Lemma \ref{lemma:CE-continuous}, we have that
\begin{align*}
\liminf_{\vartheta\rightarrow \infty} C_{a_{1,\vartheta}, a_{2,\vartheta}}\geq
\liminf_{\vartheta\rightarrow \infty} J[u_{1,\vartheta},u_{2,\vartheta}] -\varepsilon
\geq C_{a_1,a_2}-\varepsilon.
\end{align*}
By the arbitrary of $\varepsilon$, we obtain that
\begin{equation}\label{eq:20210413-e6}
\liminf_{\vartheta\rightarrow \infty} C_{a_{1,\vartheta}, a_{2,\vartheta}}\geq C_{a_1,a_2}.
\end{equation}

If $a_1>0,a_2>0$, we put
$t_{\iota,\vartheta}=\left(\frac{a_{\iota,\vartheta}}{a_\iota}\right)^{\frac{1}{2}},\ \iota=1,2.$
$\forall \varepsilon>0, \exists~(u_1,u_2)\in S_{a_1,a_2}$ such that $J[u_1,u_2]\leq C_{a_1,a_2}+\varepsilon$.
Then $(t_{1,\vartheta}u_1, t_{2,\vartheta}u_2)\in S_{a_{1,\vartheta},a_{2,\vartheta}}$ and thus
$C_{a_{1,\vartheta},a_{2,\vartheta}}\leq J[t_{1,\vartheta}u_1, t_{2,\vartheta}u_2].$
Observing that $t_{\iota,\vartheta}=1+o(1), \iota=1,2$, we obtain that
\begin{align*}
\limsup_{\vartheta\rightarrow \infty}C_{a_{1,\vartheta},a_{2,\vartheta}}\leq
\limsup_{\vartheta\rightarrow \infty}J[t_{1,\vartheta}u_1, t_{2,\vartheta}u_2]
=J[u_1,u_2]\leq C_{a_1,a_2}+\varepsilon.
\end{align*}
By the arbitrary of $\varepsilon$, we obtain that
\begin{equation}\label{eq:20210413-e7}
\limsup_{\vartheta\rightarrow \infty} C_{a_{1,\vartheta}, a_{2,\vartheta}}\leq C_{a_1,a_2}.
\end{equation}
Therefore, by \eqref{eq:20210413-e6} and \eqref{eq:20210413-e7}, we obtain that \eqref{eq:20210413-e5} holds.

If $a_1=a_2=0$, we have $C_{a_1,a_2}=0$. By Lemma \ref{lemma:min-negative}, it is trivial that
$\limsup_{\vartheta\rightarrow \infty} C_{a_{1,\vartheta}, a_{2,\vartheta}}\leq 0,$
combing with \eqref{eq:20210413-e6}, we see that the assertion \eqref{eq:20210413-e5} also holds.

If $a_1=0, a_2>0$. Let $t_\vartheta:=\left(\frac{a_{2,\vartheta}}{a_2}\right)^{\frac{1}{2}}$ and $\omega\in H^1(\mathbb{R}^N)$ such that $\|\omega\|_2^2=1$.
For $\forall ~\varepsilon>0, \exists u_2\in H^1(\mathbb{R}^N)$ such that $\|u_2\|_2^2=a_2$ and
$J[0,u_2]\leq C_{0,a_2}+\varepsilon.$
Then we see that $(\sqrt{a_{1,\vartheta}}\omega, t_\vartheta u_2)\in S_{a_{1,\vartheta}, a_{2,\vartheta}}$, thus
\begin{equation}\label{eq:20210413-e8}
C_{a_{1,\vartheta}, a_{2,\vartheta}}\leq J[\sqrt{a_{1,\vartheta}}\omega, t_\vartheta u_2].
\end{equation}
Noting that $a_{1,\vartheta}\rightarrow 0$, we can prove that
\begin{align*}
\limsup_{\vartheta\rightarrow \infty} C_{a_{1,\vartheta}, a_{2,\vartheta}}\leq&
\limsup_{\vartheta\rightarrow \infty} J[\sqrt{a_{1,\vartheta}}\omega, t_\vartheta u_2]
=\limsup_{\vartheta\rightarrow \infty} J[0, t_\vartheta u_2]\\
=&J[0, u_2]\leq  C_{0,a_2}+\varepsilon.
\end{align*}
Then by the arbitrary of $\varepsilon$ and  \eqref{eq:20210413-e6}, we also prove \eqref{eq:20210413-e5}.

And the case of $a_1>0, a_2=0$ can be proved similarly.

The proof is completed.
\end{proof}

Noting that $S_{a_1,a_2}$ is not weak closed, for this reason we introduce a new  set
\begin{equation}\label{def:Ba}
B_{a_1,a_2}:=\left\{(u_1,u_2)\in \mathcal{H}: \int_{\mathbb{R}^N}u_\iota^2 \mathrm{d} x\leq a_\iota, \iota=1,2\right\}, \forall (a_1,a_2)\in \mathbb{R}_+^2.
\end{equation}
Then it is clear that $B_{a_1,a_2}$ is weak closed (compact). Define
\begin{equation}\label{def:Da}
D_{a_1,a_2}:=\inf_{u\in B_{a_1,a_2}}J[u_1,u_2]\;\hbox{and}\;F_{a_1,a_2}:=\inf_{(u_1,u_2)\in B_{a_1,a_2}}I[u_1,u_2].
\end{equation}

\begin{lemma}\label{lemma:C=D}
For any $(a_1,a_2)\in \mathbb{R}_+^2$, we have that $C_{a_1,a_2}=D_{a_1,a_2}$ and $E_{a_1,a_2}=F_{a_1,a_2}$.
\end{lemma}
\begin{proof}
We only prove $C_{a_1,a_2}=D_{a_1,a_2}$.
The case of $(a_1,a_2)=(0,0)$ is trivial. For $(a_1,a_2)\neq (0,0)$, observing that $S_{a_1,a_2}\subset B_{a_1,a_2}$, we have that $D_{a_1,a_2}\leq C_{a_1,a_2}$. If $D_{a_1,a_2}\neq C_{a_1,a_2}$, then $D_{a_1,a_2}<C_{a_1,a_2}$. And thus there exists some $(u_1,u_2)\in B_{a_1,a_2}$ such that
$$D_{a_1,a_2}\leq J[u_1,u_2]<C_{a_1,a_2}.$$
It is clearly that $(u_1,u_2)\in B_{a_1,a_2}\backslash S_{a_1,a_2}$. Put
$b_\iota:=\int_{\mathbb{R}^N}u_\iota^2\mathrm{d} x, \iota=1,2.$
Then $(u_1,u_2)\in S_{b_1,b_2}$. Furthermore, It is easy to see that $a_\iota-b_\iota\geq 0, \iota=1,2$ and $(a_1-b_1, a_2-b_2)\neq (0,0)$. Then it follows that
$$C_{b_1,b_2}=\inf_{(v_1,v_2)\in S_{b_1,b_2}}J[v_1,v_2]\leq J[u_1,u_2]<C_{a_1,a_2},$$
a contradiction to Corollary \ref{cro:monotonicity-Ca}.
\end{proof}

For the convenience to estimate the Lagrange multipliers, we prefer to give the following result directly without proof, which essentially own to Frank H. Clarke(see \cite[Theorem 1]{Clarke1976}).
\begin{proposition}\label{prop:minimizer}
If $(u_1,u_2)$ attains $D_{a_1,a_2}$, then there exists $\lambda_\iota\geq 0, \iota=1,2$ such that
$$J'[u_1,u_2]+\lambda_1 (u_1,0)+\lambda_2 (0,u_2)=0\;\hbox{in}\;\mathcal{H}^*.$$
That is, $(u_1,u_2)$ solves equation \eqref{eq:system-general} with $\lambda_1,\lambda_2\geq 0$.
Furthermore, we have
$$\lambda_\iota(\|u_\iota\|_2^2-a_\iota)=0, \quad\iota=1,2.$$
\end{proposition}


\s{Radial potential case and Proof of Theorem \ref{thm:main-radial}}\label{radial-case}
\renewcommand{\theequation}{4.\arabic{equation}}


In this Section, we focus on the radial case and prove Theorem  \ref{thm:main-radial} via the rearrangement argument. Some properties relate to the rearrangement are provided first.
\begin{remark}\label{remark:monotone}
It is easy to see that for any fixed $s\in \mathbb{R}_+$, $\partial_1G(s,t)$ is increasing by $t\in \mathbb{R}_+$. Similarly, for any fixed $t\in \mathbb{R}_+$, $\partial_2G(s,t)$ is increasing by $s\in \mathbb{R}_+$ under the assumption \eqref{eq:G-definition}.
\end{remark}

\begin{lemma}\label{lemma:rearrange-system}
 Assume that $G(s,t)$ is given by \eqref{eq:G-definition}.
Let $u_1, u_2$ be nonnegative functions on $\mathbb{R}^N$, vanishing at infinity and let $u_1^*$ and $u_2^*$ be their symmetric-decreasing  rearrangements. Then
\begin{equation}\label{eq:rearrangement-G}
\int_{\mathbb{R}^N}G(u_1^*, u_2^*)\mathrm{d} x\geq \int_{\mathbb{R}^N}G(u_1, u_2)\mathrm{d} x.
\end{equation}
\end{lemma}
\begin{proof}
Noting that $G(0,0)=0$, we have
{\allowdisplaybreaks
\begin{align*}
&\int_{\mathbb{R}^N}G(u_1, u_2)\mathrm{d} x
=\int_{\mathbb{R}^N}\left[G(u_1, u_2)-G(0, u_2)+G(0, u_2)-G(0,0)\right]\mathrm{d} x\\
=&\int_{\mathbb{R}^N}\mathrm{d} x\int_{0}^{u_1(x)}\partial_1G(s, u_2(x))ds+\int_{\mathbb{R}^N}\mathrm{d} x\int_{0}^{u_2(x)}\partial_2G(0,t)dt\\
=&\int_{\mathbb{R}^N}\mathrm{d} x\int_0^\infty \partial_1G(s, u_2(x)) \chi_{\{x:u_1(x)>s\}}(x)ds\\
&+\int_{\mathbb{R}^N}\mathrm{d} x\int_0^\infty \partial_2G(0,t)\chi_{\{x:u_2(x)>t\}}(x)dt\\
=&\int_{\mathbb{R}^N}\mathrm{d} x\int_0^\infty\int_0^\infty\chi_{\{x:\partial_1G(s, u_2(x))>r\}}(x)\chi_{\{x:u_1(x)>s\}}(x)drds\\
&+\int_{\mathbb{R}^N}\mathrm{d} x\int_0^\infty \partial_2G(0,t)\chi_{\{x:u_2(x)>t\}}(x)dt\\
=&\int_0^\infty\int_0^\infty\left|\{x:\partial_1G(s, u_2(x))>r\}\cap \{x:u_1(x)>s\}\right|drds\\
& + \int_0^\infty \partial_2G(0,t)\left|\{x:u_2(x)>t\}\right|dt.
\end{align*}}
For $r, s\in \mathbb{R}_+$, we put
\begin{equation}
t(r,s):=\begin{cases}
\sup\{t: \partial_1G(s,t)\leq r\},\quad &\hbox{if}\; \{t: \partial_1G(s,t)\leq r\}\neq \emptyset,\\
-\infty&\hbox{if}\;\{t: \partial_1G(s,t)\leq r\}=\emptyset.
\end{cases}
\end{equation}
Then by Remark \ref{remark:monotone}, we have that $\partial_1G(s, u_2(x))>r$ if and only if $u_2(x)>t(r,s)$.
So
\begin{equation}
\left|\{x:\partial_1G(s, u_2(x))>r\}\cap \{x:u_1(x)>s\}\right|=\left|\{x:u_2(x)>t(r,s)\}\cap \{x:u_1(x)>s\}\right|,
\end{equation}
and thus
\begin{align}
\int_{\mathbb{R}^N}G(u_1, u_2)\mathrm{d} x=&\int_0^\infty\int_0^\infty\left|\{x:u_2(x)>t(r,s)\}\cap \{x:u_1(x)>s\}\right| drds\nonumber\\
&+ \int_0^\infty \partial_2G(0,t)\left|\{x:u_2(x)>t\}\right|dt.
\end{align}
Similarly, we have
\begin{align}
\int_{\mathbb{R}^N}G(u_1^*, u_2^*)\mathrm{d} x=&\int_0^\infty\int_0^\infty\left|\{x:u_2^*(x)>t(r,s)\}\cap \{x:u_1^*(x)>s\}\right| drds\nonumber\\
& +\int_0^\infty \partial_2G(0,t)\left|\{x:u_2^*(x)>t\}\right|dt.
\end{align}
Recalling the definition of symmetric-decreasing  rearrangement, we have that $\forall t, s\geq 0$,
$$\begin{cases}
\left|\{x:u_2^*(x)>t\}\right|=\left|\{x:u_2(x)>t\}\right|,\\
\left|\{x:u_1^*(x)>s\}\right|=\left|\{x:u_1(x)>s\}\right|,\\
\left|\{x:u_2^*(x)>t(r,s)\}\right|=\left|\{x:u_2(x)>t(r,s)\}\right|,
\end{cases}$$
and thus
\begin{align*}
&\left|\{x:u_2^*(x)>t(r,s)\}\cap \{x:u_1^*(x)>s\}\right|\\
=&\min\left\{\left|\{x:u_2^*(x)>t(r,s)\}\right|, \left|\{x:u_1^*(x)>s\}\right|\right\}\\
=&\min\left\{\left|\{x:u_2(x)>t(r,s)\}\right|, \left|\{x:u_1(x)>s\}\right|\right\}\\
\geq&\left|\{x:u_2(x)>t(r,s)\}\cap \{x:u_1(x)>s\}\right|.
\end{align*}
Hence, one can see that
$$\int_{\mathbb{R}^N}G(u_1^*, u_2^*)\mathrm{d} x\geq \int_{\mathbb{R}^N}G(u_1,u_2)\mathrm{d} x.$$
\end{proof}

We note that under the assumption of $(VH2)$, $|V_\iota(x)|, \iota=1,2$ are symmetric decreasing functions, i.e.,
$$-V_\iota(x)=|V_\iota(x)|=|V_\iota|^*(x), \iota=1,2.$$
Then by the simplest rearrangement inequality (see \cite[Theorem 3.4]{LiebLoss1997}), we have that
\begin{equation}\label{eq:rearrangement-V}
\|u^{*}\|_{V_\iota}^{2}=\int_{\mathbb{R}^N}|V_\iota(x)| {u^*}^2(x) \mathrm{d} x \geq \int_{\mathbb{R}^N}|V_\iota(x)| u^2(x)\mathrm{d} x=\|u\|_{V_\iota}^{2}, \iota=1,2, \forall u\in H^1(\mathbb{R}^N).
\end{equation}

Define
\begin{equation}\label{def:minimizion-rad-C}
C_{a_1,a_2}^{rad}:=\inf_{(u_1,u_2)\in S_{a_1,a_2}^{rad}} J[u_1,u_2],
\end{equation}
and
\begin{equation}\label{def:minimizion-rad-D}
D_{a_1,a_2}^{rad}:=\inf_{(u_1,u_2)\in B_{a_1,a_2}^{rad}} J[u_1,u_2],
\end{equation}
where $S_{a_1,a_2}^{rad}:=S_{a_1,a_2}\cap \mathcal{H}^{rad}, D_{a_1,a_2}^{rad}:=B_{a_1,a_2}\cap \mathcal{H}^{rad}$.
Then we have the following result.
\begin{lemma}\label{lemma:radial-minimizion}
Under the assumptions $(VH2)$ and \eqref{eq:G-definition}, we have that
$$C_{a_1,a_2}=C_{a_1,a_2}^{rad}=D_{a_1,a_2}=D_{a_1,a_2}^{rad}, \forall  a_1>0, a_2>0.$$
\end{lemma}
\begin{proof}
It is only need to prove $C_{a_1,a_2}=C_{a_1,a_2}^{rad}$ due to Lemma \ref{lemma:C=D}.
It is trivial that $C_{a_1,a_2}^{rad}\geq C_{a_1,a_2}$ since $S_{a_1,a_2}^{rad} \subset S_{a_1,a_2}$.
On the other hand, for any $(u_1,u_2)\in S_{a_1,a_2}$, we have $(u_1^*, u_2^*)\in S_{a_1,a_2}^{rad}$. Recalling the P\'olya-Szeg\"o inequality, we have
\begin{equation}\label{eq:Polya-Szego}
\int_{\mathbb{R}^N}|\nabla u_\iota^*|^2 \mathrm{d} x\leq \int_{\mathbb{R}^N}|\nabla u_\iota|^2 \mathrm{d} x, \iota=1,2.
\end{equation}
Under the assumption $(VH2)$, by Lemma \ref{lemma:rearrange-system}, combining with the formulas \eqref{eq:rearrangement-V} and \eqref{eq:Polya-Szego}, we have that
\begin{align*}
C_{a_1,a_2}^{rad}\leq J[u_1^*,u_2^*]=&\frac{1}{2}\sum_{\iota=1}^{2}\left(\|\nabla u_{\iota}^{*}\|_2^2-\|u_{\iota}^{*}\|_{V_\iota}^{2}\right)-\int_{\mathbb{R}^N}G(u_1^*,u_2^*)\mathrm{d} x\\
\leq&\frac{1}{2}\sum_{\iota=1}^{2}\left(\|\nabla u_{\iota}\|_2^2-\|u_{\iota}\|_{V_\iota}^{2}\right)-\int_{\mathbb{R}^N}G(u_1,u_2)\mathrm{d} x\\
=&J[u_1,u_2].
\end{align*}
By the arbitrary of $(u_1,u_2)\in S_{a_1,a_2}$, we obtain the opposite inequality
$C_{a_1,a_2}^{rad}\leq C_{a_1,a_2}$. Hence,
$$C_{a_1,a_2}^{rad}=C_{a_1,a_2}.$$
\end{proof}


{\bf Now we are ready to prove Theorem \ref{thm:main-radial}:}
Under the assumptions, by Lemma \ref{lemma:radial-minimizion} and Remark \ref{remark:existence-bdd-minizing}-(i), there exists  $\{(u_{1,\vartheta}, u_{2,\vartheta})\}\subset  \mathcal{H}^{rad}$ such that
$\|u_{\iota,\vartheta}\|_2^2=a_\iota, \iota=1,2, \forall \vartheta\in \mathbb{N}$
and
$$J[u_{1,\vartheta},u_{2,\vartheta}]\rightarrow C_{a_1,a_2}=D_{a_1,a_2}\;\hbox{as}\;\vartheta\rightarrow \infty.$$
Up to a subsequence, we assume that $(u_{1,\vartheta}, u_{2,\vartheta})\rightharpoonup (u_1,u_2)$ in $\mathcal{H}$ and
$$(u_{1,\vartheta}, u_{2,\vartheta})\rightarrow (u_1,u_2)\;\hbox{in}\;L_{loc}^{2}(\mathbb{R}^N)\times L_{loc}^{2}(\mathbb{R}^N), (u_{1,\vartheta}, u_{2,\vartheta})\rightarrow (u_1,u_2)\;a.e.\;\hbox{in}\;\mathbb{R}^N.$$
By the Br\'ezis-Lieb lemma \cite{BrezisLieb1983}, we have
\begin{equation}\label{eq:Brezis-Lieb}
\int_{\mathbb{R}^N}G(u_{1,\vartheta}, u_{2,\vartheta})\mathrm{d} x=\int_{\mathbb{R}^N}G(u_1,u_2)\mathrm{d} x+\int_{\mathbb{R}^N}G(u_{1,\vartheta}-u_1, u_{2,\vartheta}-u_2)\mathrm{d} x+o(1).
\end{equation}
Noting that  the embedding $H^{1,rad}(\mathbb{R}^N)\hookrightarrow\hookrightarrow L^q(\mathbb{R}^N)$ is compact for $2<q<2^*$, provided $N\geq 2$. In particular, if $N=1$, we can take $\{(u_{1,\vartheta}, u_{2,\vartheta})\}$ as non-increasing functions of $|x|$. That is,
$$u_{\iota,\vartheta}(x)=u_{\iota,\vartheta}^{*}(x), \iota=1,2,\forall \vartheta\in \mathbb{N}.$$
Then up to a subsequence, we can also have that $u_{\iota,\vartheta}\rightarrow u_{\iota}$ in $L^q(\mathbb{R})$, $\forall 2<q\leq \infty, \iota=1,2$ (see \cite[Proposition 1.7.1]{Cazenave2003}). Essentially, it is due to the fact that
$$u_{\iota,\vartheta}(x)=u_{\iota,\vartheta}(|x|)\leq C \|u_{\iota,\vartheta}\|_2 |x|^{-\frac{1}{2}}\rightarrow 0$$
as $|x|\rightarrow \infty$, uniformly in $\vartheta\in \mathbb{N}$ and $\iota=1,2$.
So by Lemma \ref{lemma:20210331-l1}, we have that
\begin{equation}\label{eq:20210319-e1}
\lim_{\vartheta\rightarrow \infty}\int_{\mathbb{R}^N}G(u_{1,\vartheta}-u_1, u_{2,\vartheta}-u_2)\mathrm{d} x=0.
\end{equation}

It follows from \eqref{eq:Brezis-Lieb} and \eqref{eq:20210319-e1} that
\begin{equation}\label{eq:20210319-e2}
\lim_{\vartheta\rightarrow \infty}\int_{\mathbb{R}^N}G(u_{1,\vartheta}, u_{2,\vartheta})\mathrm{d} x=\int_{\mathbb{R}^N}G(u_1, u_2)\mathrm{d} x.
\end{equation}
On the other hand, by the weak lower semi-continuity of norm and $\lim_{|x|\rightarrow \infty}V_\iota(x)=0$, one can prove that
\begin{equation}\label{eq:20210319-e3}
\int_{\mathbb{R}^N}|\nabla u_\iota|^2 +V_\iota(x)|u_\iota|^2 \mathrm{d} x \leq \liminf_{\vartheta\rightarrow \infty} \int_{\mathbb{R}^N}|\nabla u_{\iota,\vartheta}|^2 +V_\iota(x)|u_{\iota,\vartheta}|^2 \mathrm{d} x, \iota=1,2.
\end{equation}
Hence, by \eqref{eq:20210319-e2} and \eqref{eq:20210319-e3}, we have that
\begin{equation}\label{eq:20210319-e4}
J[u_1,u_2]\leq \lim_{\vartheta\rightarrow \infty}J[u_{1,\vartheta}, u_{2,\vartheta}]=C_{a_1,a_2}=D_{a_1,a_2}.
\end{equation}
It is clear that $(u_1,u_2)\in B_{a_1,a_2}$. Put
$b_\iota=\|u_\iota\|_2^2, \iota=1,2,$
then  $(u_1,u_2)\in S_{b_1,b_2}$ with $b_1\leq a_1,b_2\leq a_2$, and thus
 $$D_{a_1,a_2}\leq C_{b_1,b_2}\leq J[u_1,u_2]\leq C_{a_1,a_2}.$$
 If $(u_1,u_2)\in B_{a_1,a_2}\backslash S_{a_1,a_2}$, we have $a_\iota-b_\iota\geq 0, \iota=1,2$ and $(a_1-b_1, a_2-b_2)\neq (0,0)$, then by Corollary \ref{cro:monotonicity-Ca}, we have that
 $C_{b_1,b_2}>C_{a_1,a_2}$, a contradiction.
Hence, we have that $(u_1,u_2)\in S_{a_1,a_2}$ with $J[u_1,u_2]=C_{a_1,a_2}$. Finally by Proposition \ref{prop:minimizer}, $\lambda_\iota\geq 0~(\iota=1,2)$, and we  finish the proof of Theorem \ref{thm:main-radial}. We also note that $\lambda_\iota> 0~(\iota=1,2)$ if $1\leq N\leq 4$.
\hfill$\Box$

\begin{corollary}\label{cro:IE}
For any $(a_1,a_2)\in \mathbb{R}_+^2$, $E_{a_1,a_2}$ is attained.
\end{corollary}
\begin{proof}
By taking $V_\iota(x)\equiv 0, \iota=1,2$, then it follows Theorem \ref{thm:main-radial} that $E_{a_1,a_2}$ is attained by some nonnegative symmetric decreasing function of $|x|$.
\end{proof}


\s{ The sharp decay estimation of positive solution}\label{sec:decay-estimation}
\renewcommand{\theequation}{5.\arabic{equation}}

In this section, we shall give the sharp decay estimation of  the positive solution for the {\it fixed frequency problem}
\begin{equation}\label{eq:system-20210329-1}
\begin{cases}
-\Delta u_1+V_1(x)u_1+\lambda_1 u_1=\partial_1G(u_1,u_2)\;\quad&\hbox{in}\;\mathbb{R}^N,\\
-\Delta u_2+V_2(x)u_2+\lambda_2 u_2=\partial_2G(u_1,u_2)\;\quad&\hbox{in}\;\mathbb{R}^N,
\end{cases}
\end{equation}
where $G(s,t)$ is given by \eqref{eq:G-definition} with the parameters satisfying \eqref{eq:parameters-paixu} and $V_1,\ V_2$ are given potentials.

The following Lemma is due to Ikoma and Miyamoto\cite[Lemma 3.3]{Ikoma2021}, which is a variant of \cite[Proposition 1.2]{BahriLi1990}
\begin{lemma}\label{lemma:Ikoma-Miyamoto}(cf.\cite{Ikoma2021})
Let $0\leq f(x), g(x)$ satisfy $f,g\in C(\mathbb{R}^N\backslash\{0\})\cap L^1(\mathbb{R}^N)$ and
$$\lim_{|x|\rightarrow \infty}(1+|x|)^\alpha e^{\beta|x|}g(x)=c\in [0,\infty), f(x)\leq C e^{-\gamma|x|}\;\hbox{for each}\;|x|\geq 1$$
for some $\alpha\geq 0,0\leq \beta<\gamma$. Then
$$\lim_{r\rightarrow \infty} (1+r)^\alpha e^{\beta r}\int_{\mathbb{R}^N}g(r\omega -y)f(y)dy=c \int_{\mathbb{R}^N}f(y)e^{\beta \omega\cdot y}dy\;\hbox{unifromly with respect to}\;\omega\in \mathcal{S}^{N-1}.$$
\end{lemma}

\begin{corollary}\label{cro:20210329-zcro1}
Let $\lambda>0$ and $G_\lambda$ be the Green function of $-\Delta +\lambda$ in $\mathbb{R}^N$. Then for any $0<\nu<\lambda$ and $t>0$,
\begin{equation}\label{eq:20210329-be1}
\lim_{r\rightarrow \infty} (1+r)^t e^{\nu r}\int_{\mathbb{R}^N}G_\lambda(r\omega-y)(1+|y|)^{-t}e^{-\nu|y|}dy=\int_{\mathbb{R}^N}G_\lambda (y) e^{\nu \omega\cdot y}dy>0
\end{equation}
uniformly with respect to $\omega\in \mathcal{S}^{N-1}$.
\end{corollary}
\begin{proof}
For $\eta=\frac{\lambda-\nu}{2}>0$, we have $\nu+\eta=\frac{\lambda+\nu}{2}<\lambda$. By comparison principle, one can find some $C>0$ such that
$$G_\lambda(x)\leq C e^{-(\nu+\eta)|x|}\;\hbox{for each}~|x|\geq 1.$$
So we can put $g(x)=(1+|x|)^{-t}e^{-\nu |x|}, f(x)=G_\lambda(x)$. And then apply Lemma \ref{lemma:Ikoma-Miyamoto} with $\alpha=t, \beta=\nu$, one can see that \eqref{eq:20210329-be1} holds uniformly with respect to $\omega\in \mathcal{S}^{N-1}$.
Here we use the property of convolution $[f\star g](r\omega)=[g\star f](r\omega)$ and note that $c=1$ in this application.
\end{proof}

\begin{lemma}\label{lemma:partial-G}
Let $G(s,t)$ be given by \eqref{eq:G-definition} with the parameters satisfying \eqref{eq:parameters-paixu}. We simply write  $(r_{1,1},r_{2,1})$ by $(r_1, r_2)$. Then there exists some $C>0$ such that
\begin{equation}\label{eq:partial-G-1}
\begin{cases}
 \partial_1G(s,t)\geq C(s^{p_1-1}+s^{r_1-1}t^{r_2})\\
 \partial_2G(s,t)\geq C(t^{q_1-1}+s^{r_1}t^{r_2-1})
 \end{cases}\;\hbox{for any}~(s,t)\in \mathbb{R}_+^2.
 \end{equation}
Furthermore, for any $M>0$,  there exists $C_{M}>0$ such that
\begin{equation}\label{eq:partial-G-2}
\begin{cases}
\partial_1G(s,t)\leq C_M\left(s^{p_1-1}+s^{r_1-1}(s+t)^{r_2}\right)\\
\partial_2G(s,t)\leq C_M\left(t^{q_1-1}+t^{r_2-1}(s+t)^{r_1}\right)
\end{cases} \;\hbox{for}~(s,t)~\in [0,M]\times [0,M].
\end{equation}
\end{lemma}
\begin{proof}
It is trivial.
\end{proof}

Basing on \eqref{eq:partial-G-1} and \eqref{eq:partial-G-2} with $p_1,q_1>2, r_1,r_2>1$ , we can obtain the sharp decay estimation of  the positive solution for the {\it fixed frequency problem  (\ref {eq:system-20210329-1})}. Precisely, we have the following result.
\begin{proposition}\label{prop:sharp-estimation}
Suppose $G(s,t)$ satisfies \eqref{eq:partial-G-1} and \eqref{eq:partial-G-2} with $p_1,q_1>2, r_1,r_2>1$.  Assume  $V_1(x),V_2(x)\leq 0$ with
$\lim_{|x|\rightarrow \infty}V_\iota(x)=0, \iota=1,2.$
Let $(u_1,u_2)\in \mathcal{H}$ be a nonnegative  solution of (\ref {eq:system-20210329-1}) with $\lambda_1>0, \lambda_2>0$.
Then we have the following sharp decay estimations.
\begin{itemize}
\item[(i)] If $u_1\not\equiv 0, u_2\equiv 0$, then for every $\nu\in (0,\lambda_1)$, there exists $c_0>0$, which is independent of $\nu$, and $C_\nu>0$ such that
$$c_0(1+|x|)^{-\frac{N-1}{2}} e^{-\sqrt{\lambda_1}|x|}\leq u_1(x)\leq C_\nu e^{-\sqrt{\nu}|x|}~\hbox{for all}~x\in \mathbb{R}^N.$$
\item[(ii)] If $u_1\equiv 0, u_2\not\equiv 0$, then for every $\nu\in (0,\lambda_2)$, there exists $c_0>0$, which is independent of $\nu$, and $C_\nu>0$ such that
$$c_0(1+|x|)^{-\frac{N-1}{2}} e^{-\sqrt{\lambda_2}|x|}\leq u_2(x)\leq C_\nu e^{-\sqrt{\nu}|x|}~\hbox{for all}~x\in \mathbb{R}^N.$$
\item[(iii)] Assume $u_1\not\equiv 0,u_2\not\equiv 0$.
\begin{itemize}
\item[(iii-1)]Put $$\bar{\lambda}_1:=\min\left\{\lambda_1, \frac{r_2^2}{(2-r_1)_+^2}\lambda_2\right\}.$$ Then for any $0<\nu_2<\bar{\lambda}_1< \nu_1$, there exists $C_{\nu_2}\geq C_{\nu_1}>0$ such that
    \begin{equation}\label{eq:20210328-estimation1}
    C_{\nu_1}e^{-\sqrt{\nu_1}|x|}\leq u_1(x)\leq C_{\nu_2}e^{-\sqrt{\nu_2}|x|}\;\hbox{in}\;\mathbb{R}^N.
    \end{equation}
    In particular, if $\bar{\lambda}_1=\lambda_1$, the lower bound can be improved as that there exists some $c_1>0$ such that
    \begin{equation}\label{eq:20210329-e1}
    c_1(1+|x|)^{-\frac{N-1}{2}} e^{-\sqrt{\lambda_1}|x|}\leq u_1(x)\;\hbox{in}~\mathbb{R}^N.
    \end{equation}
\item[(iii-2)] Assume $u_1\not\equiv 0,u_2\not\equiv 0$. Put $$\bar{\lambda}_2:=\min\left\{\lambda_2, \frac{r_1^2}{(2-r_2)_+^2}\lambda_1\right\}.$$ Then for any $0<\nu_2<\bar{\lambda}_2< \nu_1$, there exists $C_{\nu_2}\geq C_{\nu_1}>0$ such that
    \begin{equation}\label{eq:20210328-estimation2}
    C_{\nu_1}e^{-\sqrt{\nu_1}|x|}\leq u_2(x)\leq C_{\nu_2}e^{-\sqrt{\nu_2}|x|}\;\hbox{in}\;\mathbb{R}^N.
    \end{equation}
    In particular, if $\bar{\lambda}_2=\lambda_2$, the lower bound can be improved as that there exists some $c_2>0$ such that
    \begin{equation}\label{eq:20210329-e2}
    c_2(1+|x|)^{-\frac{N-1}{2}} e^{-\sqrt{\lambda_2}|x|}\leq u_2(x)\;\hbox{in}~\mathbb{R}^N.
    \end{equation}
\end{itemize}
\end{itemize}
\end{proposition}

The proof of this proposition is long but necessary. In order to avoid readers getting lost in such a complex detail of proof, we tend to give its proof in the Appendix A.

\begin{remark}\label{remark:20210407-r1}
We emphasize that $V_\iota(x)\equiv 0$ is within our consideration in Proposition \ref{prop:sharp-estimation}.
\end{remark}

\begin{corollary}\label{cro:20210329-crow1}
Suppose $G(s,t)$ satisfying \eqref{eq:partial-G-1} and \eqref{eq:partial-G-2} with $p_1,q_1>2, r_1,r_2>1$.  Assume  $V_1(x),V_2(x)\leq 0$ with
$\lim_{|x|\rightarrow \infty}V_\iota(x)=0, \iota=1,2.$
Let $0<\lambda_1\leq \lambda_2$. Suppose that $(u_1,u_2)$ is a nonnegative nontrivial  solution to (\ref {eq:system-20210329-1})
and $(\omega_1, \omega_2)$ is a nonnegative nontrivial solution to
\begin{equation}\label{eq:system-20210329-2}
\begin{cases}
-\Delta \omega_1+\lambda_1 \omega_1=\partial_1G(\omega_1,\omega_2)\;\quad&\hbox{in}\;\mathbb{R}^N,\\
-\Delta \omega_2+\lambda_2 \omega_2=\partial_2G(\omega_1,\omega_2)\;\quad&\hbox{in}\;\mathbb{R}^N.
\end{cases}
\end{equation}
Define
\begin{equation}\label{eq:20210329-xhhe2}
\omega_{\iota,R}(x):=\omega_\iota(x-R{\bf e_1}),~\iota=1,2.
\end{equation}
We have the following conclusions:
\begin{itemize}
\item[(i)]If $u_1\neq 0$ and $\omega_1\neq 0$, then for any $\lambda\in (0,\lambda_1)$, we can find some $c>0$ independent of $\lambda$ and another $C_\lambda>0$ such that
    \begin{equation}\label{20210329-we1}
    c(1+R)^{-\frac{N-1}{2}} e^{-\sqrt{\lambda_1}R}\leq \int_{\mathbb{R}^N}u_1\omega_{1,R}\mathrm{d} x\leq C_\lambda e^{-\sqrt{\lambda}R}.
    \end{equation}
    Furthermore, if $\lambda_1=\lambda_2$, and $u_2\neq 0$,$\omega_2\neq 0$, then \eqref{20210329-we1} also holds for $\int_{\mathbb{R}^N}u_2\omega_{2,R}\mathrm{d} x$. And if $\lambda_1<\lambda_2$, then we can find some $\theta_1>1$ and $C>0$ such that
    \begin{equation}\label{eq:20210408-bue1}
    \int_{\mathbb{R}^N}u_2\omega_{2,R}\mathrm{d} x\leq C \left(\int_{\mathbb{R}^N}u_1\omega_{1,R}\mathrm{d} x\right)^{\theta_1}.
    \end{equation}
    \item[(ii)] Suppose $\bar{\lambda}_2=\lambda_2$ or $\bar{\lambda}_2<\lambda_2$ with $u_1=\omega_1=0$. If $u_2\neq 0$ and $\omega_2\neq 0$, then for any $\lambda\in (0,\lambda_2)$,
        we can find some $c>0$ independent of $\lambda$ and another $C_\lambda>0$ such that
    \begin{equation}\label{20210329-we2}
    c(1+R)^{-\frac{N-1}{2}} e^{-\sqrt{\lambda_2}R}\leq \int_{\mathbb{R}^N}u_2\omega_{2,R}\mathrm{d} x\leq C_\lambda e^{-\sqrt{\lambda}R}.
    \end{equation}
    \item[(iii)] Suppose $\bar{\lambda}_2<\lambda_2$ and $(u_1,\omega_1)\neq (0,0)$. If $u_2\neq 0$ and $\omega_2\neq 0$, then for any $0<\nu_1<\bar{\lambda}_2<\nu_2$, there exists $C_{\nu_1},C_{\nu_2}>0$ such that
        \begin{equation}\label{20210329-we3}
        C_{\nu_2}e^{-\sqrt{\nu_2}R}\leq \int_{\mathbb{R}^N}u_2\omega_{2,R}\mathrm{d} x\leq C_{\nu_1}e^{-\sqrt{\nu_1}R}.
        \end{equation}
    \item[(iv)]If $u_1\neq 0,  u_2\neq 0, \omega_1\neq 0$ and $\omega_2\neq 0$, then for any $\gamma_1\geq r_1, \gamma_2\geq r_2$ with $\eta\in (0, \min\{r_1-1,r_2-1\})$, there exists some $\theta=\theta_{\gamma_1,\gamma_2,\eta}\in (1,2)$ such that
        \begin{equation}\label{20210409-e4}
        \begin{cases}
        \displaystyle \int_{\mathbb{R}^N}u_{1}^{1+\eta}\omega_{1,R}^{\gamma_1-1-\eta}\omega_{2,R}^{\gamma_2}=o\left((\int_{\mathbb{R}^N}u_1\omega_{1,R}\mathrm{d} x)^{\theta}\right),\\
        \displaystyle \int_{\mathbb{R}^N}u_{2}^{1+\eta}\omega_{1,R}^{\gamma_1}\omega_{2,R}^{\gamma_2-1-\eta}=o\left((\int_{\mathbb{R}^N}u_1\omega_{1,R}\mathrm{d} x)^{\theta}\right).
        \end{cases}
        \end{equation}
\end{itemize}
\end{corollary}
\begin{proof}
(i) When $\lambda_1\leq \lambda_2$, by Proposition \ref{prop:sharp-estimation}, we see that $\bar{\lambda}_1=\lambda_1$. And then for any $\lambda\in (0,\lambda_1)$, we have
$$\lim_{|x|\rightarrow \infty}e^{\sqrt{\lambda}|x|} u_1(x)=0, \omega_1(x)\leq C_\lambda e^{-\sqrt{\frac{\lambda+\lambda_1}{2}}|x|}\;\hbox{in}\;\mathbb{R}^N.$$
Then by Lemma \ref{lemma:Ikoma-Miyamoto}, we have that
$$\lim_{R\rightarrow \infty} e^{\sqrt{\lambda}R}\int_{\mathbb{R}^N} u_1(R\omega+y)\omega_1(y)dy=0\;\hbox{uniformly with respect to}~\omega\in \mathcal{S}^{N-1}.$$
It follows that $$\int_{\mathbb{R}^N}u_1\omega_{1,R}\mathrm{d} x=o(e^{-\sqrt{\lambda}R})\;\hbox{as}\;R\rightarrow \infty.$$
Hence, we can find some $C_\lambda>0$ such that
$\int_{\mathbb{R}^N}u_1\omega_{1,R}\mathrm{d} x\leq C_\lambda e^{-\sqrt{\lambda}R}.$
On the other hand, by Proposition \ref{prop:sharp-estimation} again, there exists some $c>0$ such that
$$c(1+|x|)^{-\frac{N-1}{2}} e^{-\sqrt{\lambda_1}|x|}\leq u_1(x), c e^{-\sqrt{\lambda_1+1}|x|}\leq \omega_1(x).$$
Then by Lemma \ref{lemma:Ikoma-Miyamoto} again, we have that
\begin{align*}
&(1+R)^{\frac{N-1}{2}}e^{\sqrt{\lambda_1}R} \int_{\mathbb{R}^N}u_1\omega_{1,R}\mathrm{d} x\\
=&(1+R)^{\frac{N-1}{2}}e^{\sqrt{\lambda_1}R} \int_{\mathbb{R}^N}u_1(x+R{\bf e_1})\omega_{1}(x)\mathrm{d} x\\
\geq&c^2(1+R)^{\frac{N-1}{2}}e^{\sqrt{\lambda_1}R} \\
&\times \int_{\mathbb{R}^N} (1+|x+R{\bf e_1}|)^{-\frac{N-1}{2}}e^{-\sqrt{\lambda_1}|x+R{\bf e_1}|}e^{-\sqrt{\lambda_1+1}|x|}\mathrm{d} x\\
=&c^2(1+o_R(1))e^{\sqrt{\lambda_1}R}\int_{\mathbb{R}^N}e^{-\sqrt{\lambda_1}|R{\bf (-e_1)}-x|}e^{-\sqrt{\lambda_1+1}|x|}\mathrm{d} x\\
\rightarrow& c^2 \int_{\mathbb{R}^N}e^{-\sqrt{\lambda_1+1}|y|}e^{-\sqrt{\lambda_1}{\bf e_1}\cdot y}dy>0.
\end{align*}
Hence, the assertion of \eqref{20210329-we1} holds.

Finally, if $\lambda_1<\lambda_2$, we remark that $\lambda_1<\bar{\lambda}_2$, so we can take $\varepsilon>0$ small such that $\lambda_1+3\varepsilon<\bar{\lambda}_2$. Then, we have that
$$0\leq u_2(x)\leq Ce^{-\sqrt{\lambda_1+2\varepsilon}|x|}, \omega_2(x)\leq C e^{-\sqrt{\lambda_1+3\varepsilon}|x|}.$$
Hence, Lemma \ref{lemma:Ikoma-Miyamoto} yields
$\int_{\mathbb{R}^N}u_2 \omega_{2,R}\mathrm{d} x=o( e^{-\sqrt{\lambda_1+\varepsilon}R})~\hbox{as}~R\rightarrow +\infty.$
Thus, we can find some $\theta_1>1$ and $C>0$ such that \eqref{eq:20210408-bue1} holds.

(ii) When $\bar{\lambda}_2=\lambda_2$ or $\bar{\lambda}_2<\lambda_2$ with $u_1=\omega_1=0$, by Proposition \ref{prop:sharp-estimation},
\eqref{20210329-we2} can be proved as \eqref{20210329-we1}.

(iii) Suppose $\bar{\lambda}_2<\lambda_2$ and $(u_1,\omega_1)\neq (0,0)$. If $u_1\neq 0, \omega_1\neq 0$, then for any $0<\nu_1<\bar{\lambda}_2<\nu_2$, we can find some $C, c>0$ such that
\begin{equation}\label{eq:20231004-e1}
\lim_{|x|\rightarrow \infty}u_2(x)e^{\sqrt{\nu_1}|x|}=0, \omega_2(x)\leq C e^{-\sqrt{\frac{\nu_1+\bar{\lambda}_2}{2}}|x|}\;\hbox{in}\;\mathbb{R}^N,
\end{equation}
and
\begin{equation}\label{eq:20231004-e2}
c e^{-\sqrt{\nu_2}|x|}\leq u_2(x), ~ c e^{-\sqrt{\nu_2+1}|x|}\leq \omega_2(x).
\end{equation}
We have that
$$\lim_{R\rightarrow +\infty}e^{\sqrt{\nu_1}R}\int_{\mathbb{R}^N}u_2\omega_{2,R}\mathrm{d} x=0$$
and
\begin{align*}
e^{\sqrt{\nu_2}R} \int_{\mathbb{R}^N}u_2\omega_{2,R}\mathrm{d} x
\geq&c^2e^{\sqrt{\nu_2}R} \int_{\mathbb{R}^N}e^{-\sqrt{\nu_2}|x|} e^{-\sqrt{\nu_2+1}|x-R{\bf e_1}|}\mathrm{d} x\\
=&c^2e^{\sqrt{\nu_2}R} \int_{\mathbb{R}^N}e^{-\sqrt{\nu_2}|x+R{\bf e_1}|} e^{-\sqrt{\nu_2+1}|x|}\mathrm{d} x\\
\rightarrow &c^2 \int_{\mathbb{R}^N} e^{-\sqrt{\nu_2+1}|x|} e^{-\sqrt{\nu_2}{\bf e_1}\cdot x}\mathrm{d} x>0~\hbox{as}~R\rightarrow +\infty.
\end{align*}
Hence, there exists $C_{\nu_1}>0,C_{\nu_2}>0$  such that \eqref{20210329-we3} holds.

If $u_1=0, \omega_1\neq 0$, then we can choose some suitable positive constants such that
\begin{equation}\label{eq:20231004-e3}
\begin{cases}
C_{\varepsilon_0}e^{-\sqrt{\lambda_2+\varepsilon_0}|x|}\leq c_0(1+|x|)^{-\frac{N-1}{2}}e^{-\sqrt{\lambda_2}|x|}\leq u_2(x)\leq C_{\varepsilon_1} e^{-\sqrt{\lambda_2-\varepsilon_1}|x|},\\
C_{\varepsilon_2}e^{-\sqrt{\bar{\lambda}_2+\varepsilon_2}|x|}\leq \omega_2(x)\leq C_{\varepsilon_3}e^{-\sqrt{\bar{\lambda}_2-\varepsilon_3}|x|}.
\end{cases}
\end{equation}
By $\bar{\lambda}_2<\lambda_2$, we can take $\varepsilon_1,\varepsilon_3$ small enough that
\begin{equation}\label{eq:20231004-e4}
\int_{\mathbb{R}^N}u_2\omega_{2,R}\mathrm{d} x\leq C_{\varepsilon_1,\varepsilon_3} e^{-\sqrt{\bar{\lambda}_2-\varepsilon_3}R}.
\end{equation}
In particular, by taking $\varepsilon_3<\bar{\lambda}_2-\nu_1$, we conclude the right hand side of \eqref{20210329-we3}.
Also by $\bar{\lambda}_2<\lambda_2$, one has that
\begin{equation}\label{eq:20231004-e5}
\int_{\mathbb{R}^N}u_2\omega_{2,R}\mathrm{d} x\geq C_{\varepsilon_0,\varepsilon_2} e^{-\sqrt{\bar{\lambda}_2+\varepsilon_2}R}.
\end{equation}
In particular, by taking $\varepsilon_2<\nu_2-\bar{\lambda}_2$, we conclude that left hand side of \eqref{20210329-we3}.

The case of $u_1\neq 0, \omega_1=0$ can be proved in a similar way.

(iv)
Thanks to $1+\eta>1,\gamma_1-1-\eta +r_2>r_2>1$ and Proposition \ref{prop:sharp-estimation}, we may find an $\eta_1>0$ such that
\begin{equation}\label{eq:20210409-e5}
\omega_{1}^{\gamma_1-1-\eta}\omega_{2}^{\gamma_2}(x)+u_{1}^{1+\eta}(x)\leq C e^{-\sqrt{\lambda_1+\eta_1}|x|}.
\end{equation}
Therefore, exploiting Lemma \ref{lemma:Ikoma-Miyamoto}, combining with the conclusion (i) above (see formula \eqref{20210329-we1}), we can find some $\tilde{\theta}_1\in (1,2)$ such that
$$\int_{\mathbb{R}^N} u_{1}^{1+\eta}\omega_{1,R}^{\gamma_1-1-\eta} \omega_{2,R}^{\gamma_2}=o\left((\int_{\mathbb{R}^N}u_1\omega_{1,R}\mathrm{d} x)^{\tilde{\theta}_1}\right).$$
Similarly, we may find some $\tilde{\theta}_2\in (1,2)$ such that
$$\int_{\mathbb{R}^N}u_{2}^{1+\eta}\omega_{1,R}^{\gamma_1}\omega_{2,R}^{\gamma_2-1-\eta}=o\left((\int_{\mathbb{R}^N}u_1\omega_{1,R}\mathrm{d} x)^{\tilde{\theta}_2}\right).$$
Then we can take $\theta:=\min\{\tilde{\theta}_1,\tilde{\theta}_2\}$.
\end{proof}


\s{ The strict sub-additive inequality}\label{sec:interaction-estimate}
\renewcommand{\theequation}{6.\arabic{equation}}
In this Section, we shall prove the strict sub-additive inequality associated with the  functionals  $J$ and $I$ basing on sharp decay estimations of the positive solutions for the {\it fixed frequency problem  (\ref {eq:system-20210329-1})} obtained in Section \ref{sec:decay-estimation}.

The following is prepared to give the interaction estimation for the general couples terms $(a_1+b_1)^{\gamma_1} (a_2+b_2)^{\gamma_2}$.

\begin{lemma}\label{lemma:20210401-l1}
Let $\gamma_1>1, \gamma_2>1$. For any $0<\eta<\min\{\gamma_1,\gamma_2\}-1$, There exists some $C_{\gamma_1,\gamma_2,\eta}>0$ such that for any $(x, y)\in \mathbb{R}_+^2$,
\begin{align}\label{20210401-xe1}
(1+x)^{\gamma_1}(1+y)^{\gamma_2}\geq &1+x^{\gamma_1}y^{\gamma_2}+\gamma_1x+\gamma_2 y \nonumber\\
&-C_{\gamma_1,\gamma_2,\eta}\left(x^{\gamma_1-1-\eta}y^{\gamma_2}+x^{\gamma_1}y^{\gamma_2-1-\eta}\right).
\end{align}
\end{lemma}
\begin{proof}
Observing that for $\gamma>1$, one can see that
$$\begin{cases}
x^{\gamma-1}<\gamma, &\forall ~0<x<\gamma^{\frac{1}{\gamma-1}}, \\
x^{\gamma-1}>\gamma , &\forall~ x>\gamma^{\frac{1}{\gamma-1}}.
\end{cases}$$
Then for any $\displaystyle(x, y)\in \left[0,\gamma_{1}^{\frac{1}{\gamma_1-1}}\right]\times \left[0,\gamma_{2}^{\frac{1}{\gamma_2-1}}\right]$, by Bernoulli inequality,
\begin{align*}
(1+x)^{\gamma_1}(1+y)^{\gamma_2}\geq &(1+\gamma_1 x) (1+\gamma_2 y)\\
=&1+\gamma_1 x+\gamma_2 y+(\gamma_1 x)\cdot (\gamma_2 y)\\
\geq&1+\gamma_1 x+\gamma_2 y+x^{\gamma_1}y^{\gamma_2},
\end{align*}
the assertion \eqref{20210401-xe1} holds.

For $(x, y)\in (\gamma_{1}^{\frac{1}{\gamma_1-1}},\infty)\times (\gamma_{2}^{\frac{1}{\gamma_2-1}},\infty)$, we can take $C=C_{\gamma_1,\gamma_2,\eta}$ large enough such that
$$\frac{2}{C}<\min\left\{\gamma_{1}^{\frac{\gamma_1-1-\eta}{\gamma_1-1}},\gamma_{2}^{\frac{\gamma_2-1-\eta}{\gamma_2-1}}\right\}.$$
And then
$$C\left(x^{\gamma_1-1-\eta}y^{\gamma_2}+x^{\gamma_1}y^{\gamma_2-1-\eta}\right)\geq 2 y^{\gamma_2}+2x^{\gamma_1}>1+\gamma_1x+\gamma_2y.$$
Since
$(1+x)^{\gamma_1}(1+y)^{\gamma_2}>x^{\gamma_1}y^{\gamma_2},$
we have that
\begin{align*}
(1+x)^{\gamma_1}(1+y)^{\gamma_2}>1+\gamma_1x+\gamma_2y+x^{\gamma_1}y^{\gamma_2}-C\left(x^{\gamma_1-1-\eta}y^{\gamma_2} +x^{\gamma_1}y^{\gamma_2-1-\eta}\right).
\end{align*}
Hence the assertion \eqref{20210401-xe1} also holds.

For $y\in (\gamma_{2}^{\frac{1}{\gamma_2-1}},+\infty)$, by the mean value theorem, there  exists some $z\in [y, y+1]$ such that
$$(1+y)^{\gamma_2}-y^{\gamma_2}=\gamma_2 z^{\gamma_2-1}\geq \gamma_2 y^{\gamma_2-1}>\gamma_2^2>1,$$
which implies that $(1+y)^{\gamma_2}>1+y^{\gamma_2}.$
 So if $\displaystyle(x, y)\in [0,\gamma_{1}^{\frac{1}{\gamma_1-1}}]\times (\gamma_{2}^{\frac{1}{\gamma_2-1}},+\infty)$,
we have
\begin{align*}
(1+x)^{\gamma_1}(1+y)^{\gamma_2}>&(1+\gamma_1 x)\cdot (1+y^{\gamma_2})\\
=&1+\gamma_1 x+\gamma_1 xy^{\gamma_2}+y^{\gamma_2}\\
\geq&1+\gamma_1 x+x^{\gamma_1}y^{\gamma_2}+\gamma_2 y.
\end{align*}
Hence, the assertion \eqref{20210401-xe1} also holds.
The case of $\displaystyle(x, y)\in (\gamma_{1}^{\frac{1}{\gamma_1-1}},+\infty)\times [0,\gamma_{2}^{\frac{1}{\gamma_2-1}}]$ can be proved similarly.
\end{proof}

\begin{corollary}\label{cro:20210409-c1}
Let $\gamma_1>1, \gamma_2>1$. For any $0<\eta<\min\{\gamma_1,\gamma_2\}-1$, There exists some $C_{\gamma_1,\gamma_2,\eta}>0$ such that for any $a_\iota,b_\iota\geq 0$, $\iota=1,2$,
\begin{equation}\label{091531}
\begin{array}{ll}
(a_1+b_1)^{\gamma_1}(a_2+b_2)^{\gamma_2}\geq &a_{1}^{\gamma_1}a_{2}^{\gamma_2}+\gamma_1 a_{1}^{\gamma_1-1}a_{2}^{\gamma_2}b_1+\gamma_2a_{1}^{\gamma_1}a_{2}^{\gamma_2-1}b_2
+b_{1}^{\gamma_1}b_{2}^{\gamma_2}\\
&-C_{\gamma_1,\gamma_2,\eta}\left\{a_{1}^{1+\eta}b_{1}^{\gamma_1-1-\eta}b_{2}^{\gamma_2} +a_{2}^{1+\eta}b_{2}^{\gamma_2-1-\eta}b_{1}^{\gamma_1}\right\}.
\end{array}
\end{equation}
Furthermore,  if $a_1=0$, we have that
\begin{equation}\label{eq:20210409-xe1}
(a_1+b_1)^{\gamma_1}(a_2+b_2)^{\gamma_2}\geq b_{1}^{\gamma_1}\max\left\{a_{2}^{\gamma_2}+\gamma_2 a_{2}^{\gamma_2-1}b_2,b_{2}^{\gamma_2}+\gamma_2 b_{2}^{\gamma_2-1}a_2\right\}.
\end{equation}
Similarly, if $b_1=0$, we have that
\begin{equation}\label{eq:20240911-e1}
(a_1+b_1)^{\gamma_1}(a_2+b_2)^{\gamma_2}\geq a_{1}^{\gamma_1}\max\left\{a_{2}^{\gamma_2}+\gamma_2 a_{2}^{\gamma_2-1}b_2,b_{2}^{\gamma_2}+\gamma_2 b_{2}^{\gamma_2-1}a_2\right\}.
\end{equation}
\end{corollary}
\begin{proof} If $(a_1+b_1)=0$ or $ (a_2+b_2)=0$, it is trivial. So we only need to consider the case of $(a_1+b_1)>0$ and  $ (a_2+b_2)>0$. Therefore, without loss of generality, we can
 assume that $a_1>0$ and $a_2>0$, since $a_\iota,b_\iota\geq 0$, $\iota=1,2$. By putting $x=\frac{b_1}{a_1}, y=\frac{b_2}{a_2}$ in \eqref{20210401-xe1}, we
 can get the inequalities \eqref{091531}.

 In particular, \eqref{eq:20210409-xe1} and \eqref{eq:20240911-e1} follow by the well known Bernoulli's inequality.
\end{proof}

Put
$a_\iota:=\|u_\iota\|_2^2,\ \  b_{\iota}:=\|\omega_\iota\|_2^2.$
Define
\begin{equation}\label{eq:20210401-e1}
\omega_{\iota, R}:=\omega_{\iota}(\cdot-R{\bf e_1}), \ \  \iota=1,2,
\end{equation}
\begin{equation}\label{eq:20210401-e2}
\sigma_{\iota,R}:=\int_{\mathbb{R}^N}u_\iota \omega_{\iota, R} \mathrm{d} x, \ \ \iota=1,2.
\end{equation}
We remark that if $u_\iota\equiv 0$ or $\omega_{\iota}\equiv 0$, we have that $\sigma_{\iota,R}\equiv 0, \forall R\in \mathbb{R}_+$. If $u_\iota\not\equiv 0$ and $\omega_{\iota}\not\equiv 0$, by the strong maximum principle, it holds that
$u_\iota(x)>0, \ \ \omega_{\iota}(x)>0\; \ \hbox{in}\; \ \mathbb{R}^N.$
And thus
\begin{equation}\label{eq:20210401-e3}
\sigma_{\iota,R}>0\;\hbox{and}\;\sigma_{\iota,R}\rightarrow 0~\hbox{as}~R\rightarrow \infty.
\end{equation}
Let
\begin{equation}\label{eq:20210401-be1}
\tau_{\iota,R}^{2}:=\frac{a_\iota+b_\iota}{a_\iota+b_\iota+2\sigma_{\iota,R}}=1-\frac{2}{a_\iota+b_\iota} \sigma_{\iota,R}+O(\sigma_{\iota,R}^{2})
\end{equation}
and
$\tau_{\iota,R}=1$ if $\sigma_{\iota,R}=0$.
Then one can see that
\begin{equation}\label{eq:20210401-be2}
\tau_{\iota,R}=1-\frac{1}{a_\iota+b_\iota} \sigma_{\iota,R}+O(\sigma_{\iota,R}^{2})
\end{equation}
and
\begin{equation}\label{eq:20210401-be3}
\|\tau_{\iota,R}(u_\iota+\omega_{\iota,R})\|_2^2=a_\iota+b_\iota, \ \  \iota=1,2.
\end{equation}

\begin{lemma}\label{lemma:intersection-estimation}
Suppose that $(VH1)$ holds and $G(s,t)$ is given by \eqref{eq:G-definition} with the parameters satisfying \eqref{eq:parameters-paixu}.
Let $0<\lambda_1\leq \lambda_2$. Assume that $(u_1,u_2)\neq (0,0)$ is a nonnegative solution to \eqref{eq:system-20210329-1} and $(\omega_{1},\omega_2)\neq (0,0)$ is a nonnegative solution to \eqref{eq:system-20210329-2}. We have the following interaction estimation.
\begin{itemize}
\item[(i)] For  $i=1,2,\cdots,\ell$, there exists some $\theta_i(p_i)\in (1,\frac{p_i}{2})\cap (1,2)$ such that
{\allowdisplaybreaks
\begin{align*}
&-\frac{\mu_i}{p_i}\tau_{1,R}^{p_i}\int_{\mathbb{R}^N}(u_1+\omega_{1,R})^{p_i}\mathrm{d} x\\
\leq &-\frac{\mu_i}{p_i}\left[\int_{\mathbb{R}^N}u_{1}^{p_i} +\omega_{1,R}^{p_i}\mathrm{d} x\right]
-\mu_i\int_{\mathbb{R}^N} \left[u_{1}^{p_i-1}\omega_{1,R}+u_{1}\omega_{1,R}^{p_i-1}\right]\mathrm{d} x\\
&+\frac{\mu_i\sigma_{1,R}}{a_1+b_1} \left[\int_{\mathbb{R}^N}u_{1}^{p_i} +\omega_{1,R}^{p_i}\mathrm{d} x\right]+o(\sigma_{1,R}^{\theta_i(p_i)}).
\end{align*}}
\item[(ii)] For  $j=1,2,\cdots,m$, there exists some $\theta_j(q_j)\in (1,\frac{q_j}{2})\cap(1,2)$ such that
{\allowdisplaybreaks
\begin{align*}
&-\frac{\nu_j}{q_j}\tau_{2,R}^{q_j}\int_{\mathbb{R}^N}(u_2+\omega_{2,R})^{q_j}\mathrm{d} x\\
\leq &-\frac{\nu_j}{q_j}\left[\int_{\mathbb{R}^N}u_{2}^{q_j} +\omega_{2,R}^{q_j}\mathrm{d} x\right]
-\nu_j\int_{\mathbb{R}^N} \left[u_{2}^{q_j-1}\omega_{2,R}+u_{2}\omega_{2,R}^{q_j-1}\right]\mathrm{d} x\\
&+\frac{\nu_j\sigma_{2,R}}{a_2+b_2} \left[\int_{\mathbb{R}^N}u_{2}^{q_j} +\omega_{2,R}^{q_j}\mathrm{d} x\right]+o(\sigma_{2,R}^{\theta_j(q_j)}).
\end{align*}}
\item[(iii)]If $\sigma_{1,R}>0$,
for $k=1,2,\cdots,n$, there exists some $\theta=\theta(r_{1,k},r_{2,k})\in (1,2)$ such that
{\allowdisplaybreaks
\begin{align*}
&-\beta_k\tau_{1,R}^{r_{1,k}}\tau_{2,R}^{r_{2,k}}\int_{\mathbb{R}^N}(u_1+\omega_{1,R})^{r_{1,k}} (u_2+\omega_{2,R})^{r_{2,k}}\mathrm{d} x\\
\leq&-\beta_k \int_{\mathbb{R}^N}\left[ u_{1}^{r_{1,k}}u_{2}^{r_{2,k}}+r_{1,k}u_{1}^{r_{1,k}-1}u_{2}^{r_{2,k}}\omega_{1,R}
+r_{2,k}u_{1}^{r_{1,k}}u_{2}^{r_{2,k}-1}\omega_{2,R} +\omega_{1,R}^{r_{1,k}}\omega_{2,R}^{r_{2,k}}\right]\mathrm{d} x\\
&+\frac{\beta_k r_{1,k}}{a_1+b_1}\sigma_{1,R}\int_{\mathbb{R}^N}\left[u_{1}^{r_{1,k}}u_{2}^{r_{2,k}}
+\omega_{1,R}^{r_{1,k}}\omega_{2,R}^{r_{2,k}}\right]\mathrm{d} x\\
&+\frac{\beta_k r_{2,k}}{a_2+b_2}\sigma_{2,R}\int_{\mathbb{R}^N}\left[u_{1}^{r_{1,k}}u_{2}^{r_{2,k}}
+\omega_{1,R}^{r_{1,k}}\omega_{2,R}^{r_{2,k}}\right]\mathrm{d} x
+o(\sigma_{1,R}^{\theta})
\end{align*}}
And if $\sigma_{1,R}=0$, for $k=1,2,\cdots,n$, there exists some $\theta=\theta(r_{1,k},r_{2,k})\in (1,2)$ such that
{\allowdisplaybreaks
\begin{align*}
&-\beta_k\tau_{1,R}^{r_{1,k}}\tau_{2,R}^{r_{2,k}}\int_{\mathbb{R}^N}(u_1+\omega_{1,R})^{r_{1,k}} (u_2+\omega_{2,R})^{r_{2,k}}\mathrm{d} x\\
\leq&-\beta_k \int_{\mathbb{R}^N}\left[\omega_{1,R}^{r_{1,k}}\omega_{2,R}^{r_{2,k}}+r_{2,k}\omega_{1,R}^{r_{1,k}}\omega_{2,R}^{r_{2,k}-1}u_2\right]\mathrm{d} x\\
&+\frac{\beta_k r_{2,k}}{a_2+b_2}\sigma_{2,R}\int_{\mathbb{R}^N}\left[\omega_{1,R}^{r_{1,k}}\omega_{2,R}^{r_{2,k}}+r_{2,k}\omega_{1,R}^{r_{1,k}}\omega_{2,R}^{r_{2,k}-1}u_2\right]\mathrm{d} x+o(\sigma_{2,R}^{\theta})\\
=&-\beta_k \int_{\mathbb{R}^N}\left[\omega_{1,R}^{r_{1,k}}\omega_{2,R}^{r_{2,k}}+r_{2,k}\omega_{1,R}^{r_{1,k}}\omega_{2,R}^{r_{2,k}-1}u_2\right]\mathrm{d} x\\
&+\frac{\beta_k r_{2,k}}{a_2+b_2}\sigma_{2,R}\int_{\mathbb{R}^N}\omega_{1,R}^{r_{1,k}}\omega_{2,R}^{r_{2,k}}\mathrm{d} x+o(\sigma_{2,R}^{\theta})~\hbox{if $u_1\equiv 0$}
\end{align*}}
and
{\allowdisplaybreaks
\begin{align*}
&-\beta_k\tau_{1,R}^{r_{1,k}}\tau_{2,R}^{r_{2,k}}\int_{\mathbb{R}^N}(u_1+\omega_{1,R})^{r_{1,k}} (u_2+\omega_{2,R})^{r_{2,k}}\mathrm{d} x\\
\leq&-\beta_k \int_{\mathbb{R}^N}\left[ u_{1}^{r_{1,k}}u_{2}^{r_{2,k}}+r_{2,k}u_{1}^{r_{1,k}}u_{2}^{r_{2,k}-2}\omega_{2,R}\right]\mathrm{d} x\\
&+\frac{\beta_k r_{2,k}}{a_2+b_2}\sigma_{2,R}\int_{\mathbb{R}^N}\left[u_{1}^{r_{1,k}}u_{2}^{r_{2,k}}
\right]\mathrm{d} x+o(\sigma_{2,R}^{\theta})~\hbox{if}~\omega_1\equiv 0.
\end{align*}}
In particular, we can rewrite them unifiedly as follows:
{\allowdisplaybreaks
\begin{align*}
&-\beta_k\tau_{1,R}^{r_{1,k}}\tau_{2,R}^{r_{2,k}}\int_{\mathbb{R}^N}(u_1+\omega_{1,R})^{r_{1,k}} (u_2+\omega_{2,R})^{r_{2,k}}\mathrm{d} x\\
\leq&-\beta_k \int_{\mathbb{R}^N}\left[\omega_{1,R}^{r_{1,k}}\omega_{2,R}^{r_{2,k}}+r_{2,k}\omega_{1,R}^{r_{1,k}}\omega_{2,R}^{r_{2,k}-1}u_2+ u_{1}^{r_{1,k}}u_{2}^{r_{2,k}}+r_{2,k}u_{1}^{r_{1,k}}u_{2}^{r_{2,k}-2}\omega_{2,R}\right]\mathrm{d} x\\
&+\frac{\beta_k r_{2,k}}{a_2+b_2}\sigma_{2,R}\int_{\mathbb{R}^N}\left[u_{1}^{r_{1,k}}u_{2}^{r_{2,k}}+\omega_{1,R}^{r_{1,k}}\omega_{2,R}^{r_{2,k}}
\right]\mathrm{d} x+o(\sigma_{2,R}^{\theta})~\hbox{if}~\sigma_{1,R}= 0.
\end{align*}}
\end{itemize}
\end{lemma}

This  lemma gives a good intersection estimation,  which plays a crucial role in the proof of strict sub-additive inequality. Its proof required a sequence of complicated calculations. We shall give the details in Appendix B.
Basing on these estimations, we obtain the following strict sub-additive inequality. The proof of the following proposition will be given in Appendix C.

\begin{proposition}\label{prop:binding-inequality}
Suppose that $(VH1)$ holds and $G(s,t)$ is given by \eqref{eq:G-definition} with the parameters satisfying \eqref{eq:parameters-paixu}.
Let $0<\lambda_1\leq \lambda_2$.
Suppose that $(u_1,u_2)\neq (0,0)$ is a nonnegative  solution to (\ref {eq:system-20210329-1})
and $(\omega_1, \omega_2)\neq (0,0)$ is a nonnegative  solution to (\ref {eq:system-20210329-2}).
 Then there exists some $(\phi_1,\phi_2)\in\mathcal{H}$ such that
\begin{equation}\label{eq:20210410-e1}
(\|\phi_1\|_2^2,\|\phi_2\|_2^2)=(\|u_1\|_2^2+\|\omega_1\|_2^2,\|u_2\|_2^2+\|\omega_2\|_2^2)~\hbox{and}
~J[\phi_1,\phi_2]<J[u_1,u_2]+I[\omega_1,\omega_2].
\end{equation}
\end{proposition}


\s{ Compactness analysis and Proof of Theorem \ref{thm:main-non-radial}}\label{sec:compactness-proof-Th2}
\renewcommand{\theequation}{7.\arabic{equation}}
In  this section, we shall apply the standard concentration compactness argument to prove Theorem \ref{thm:main-non-radial}.

\begin{remark}\label{remark:20210411-r1}
We remark that for $a_1>0,a_2>0$, $S_{a_1,a_2}$ is a Hilbert manifold with codimension $2$. And for any minimizing sequence $\{(u_{1,\vartheta}, u_{2,\vartheta})\}_{\vartheta=1}^{\infty}\subset S_{a_1,a_2}$, by Ekeland's variational principle for $J$ and $\{(u_{1,\vartheta}, u_{2,\vartheta})\}$ on $S_{a_1,a_2}$, there exist $\{(v_{1,\vartheta}, v_{2,\vartheta})\}\subset S_{a_1,a_2}$ and $\{\lambda_{1,\vartheta}\}, \{\lambda_{2,\vartheta}\}\subset \mathbb{R}$ such that
$$\|(u_{1,\vartheta}-v_{1,\vartheta},u_{2,\vartheta}-v_{2,\vartheta})\|_{\mathcal{H}}\rightarrow 0$$
and
$$J'[v_{1,\vartheta}, v_{2,\vartheta}]+\lambda_{1,\vartheta} (v_{1,\vartheta},0)+\lambda_{2,\vartheta}(0, v_{2,\vartheta})\rightarrow 0\;\hbox{strongly in}~\mathcal{H}^*~\hbox{as}~\vartheta\rightarrow \infty.$$
See also \cite[Lemma 2.3]{Ikoma2021}.
\end{remark}

For  $a_1>0,a_2>0$. Let us take any minimizing sequence $\{(u_{1,\vartheta}, u_{2,\vartheta})\}_{\vartheta=1}^{\infty}\subset S_{a_1,a_2}$ for $C_{a_1,a_2}$.
Since $J[u_1,u_2]=J[|u_1|, |u_2|]$,  $\{(|u_{1,\vartheta}|, |u_{2,\vartheta})|\}_{\vartheta=1}^{\infty}$ is also a minimizing sequence. Hence, we may suppose that $\{u_{\iota,\vartheta}\}~ (\iota=1,2, \vartheta\in \mathbb{N})$ are nonnegative functions.
By Remark \ref{remark:20210411-r1}, without loss of generality, we may suppose further that
\begin{equation}\label{eq:20210412-e1}
J'[u_{1,\vartheta}, u_{2,\vartheta}]+(\lambda_{1,\vartheta} u_{1,\vartheta}, \lambda_{2,\vartheta} u_{2,\vartheta})\rightarrow 0\;\hbox{in}~\mathcal{H}^*.
\end{equation}

\begin{lemma}\label{lemma:Lagrange-multiplier-bounded}
Suppose that $(VH1)$ holds and $G(s,t)$ is given by \eqref{eq:G-definition}. For any given $a_1>0,a_2>0$, let $\{(u_{1,\vartheta}, u_{2,\vartheta})\}_{\vartheta=1}^{\infty}\subset S_{a_1,a_2}$ be a minimizing sequence such that $u_{\iota,\vartheta}, \iota=1,2, \vartheta\in \mathbb{N}$ are nonnegative functions and \eqref{eq:20210412-e1} holds, then $\{(\lambda_{1,\vartheta}, \lambda_{2,\vartheta})\}$ is bounded.
\end{lemma}
\begin{proof}
Recalling that $\{(u_{1,\vartheta}, u_{2,\vartheta})\}$ are bounded in $\mathcal{H}$ (see Remark \ref{remark:existence-bdd-minizing}), by
$$\lambda_{1,\vartheta}=\frac{\|u_{1,\vartheta}\|_{H^1}\times o(1)-J'[u_{1,\vartheta}, u_{2,\vartheta}](u_{1,\vartheta},0)}{a_1},$$
we see that $\{\lambda_{1,\vartheta}\}_{\vartheta=1}^{\infty}$ is bounded.
Similarly, by
$$\lambda_{2,\vartheta}=\frac{\|u_{2,\vartheta}\|_{H^1}\times o(1)-J'[u_{1,\vartheta}, u_{2,\vartheta}](0,u_{2,\vartheta})}{a_2},$$
we obtain that $\{\lambda_{2,\vartheta}\}_{\vartheta=1}^{\infty}$ is also bounded.
\end{proof}

\begin{lemma}\label{lemma:no-vanish}
Under the assumptions of Lemma \ref{lemma:Lagrange-multiplier-bounded}, let $Q^N:=[0,1]^N$. Then
\begin{equation}\label{eq:20210412-e2}
\limsup_{\vartheta\rightarrow \infty} \sup_{y\in \mathbb{R}^N} \|u_{\iota,\vartheta}\|_{L^2(Q^N+y)}>0, \quad \iota=1,2.
\end{equation}
\end{lemma}
\begin{proof}
We only prove the case of $\iota=1$. And the case of $\iota=2$ can be proved in a same way.
We argue by a way of negation and suppose that
$$\sup_{y\in \mathbb{R}^N} \|u_{1,\vartheta}\|_{L^2(Q^N+y)}\rightarrow 0\;\hbox{as}\;\vartheta\rightarrow \infty.$$
Then it follows that $u_{1,\vartheta}\rightarrow 0$ strongly in $L^p(\mathbb{R}^N)$ for any $2<p<2^*$, see \cite[Lemma 1.21]{Willem1996}.
Therefore, similar to \eqref{eq:20210413-e1}, we can prove that
\begin{equation}\label{eq:20210412-e3}
\liminf_{\vartheta\rightarrow \infty}\int_{\mathbb{R}^N}G(u_{1,\vartheta}, u_{2,\vartheta})\mathrm{d} x
=\liminf_{\vartheta\rightarrow \infty} G(0,u_{2,\vartheta})\mathrm{d} x,
\end{equation}
which implies that
\begin{equation}\label{eq:20210412-e4}
C_{a_1,a_2}=\lim_{\vartheta\rightarrow \infty}J[u_{1,\vartheta},u_{2,\vartheta}]\geq \liminf_{\vartheta\rightarrow \infty}J[0,u_{2,\vartheta}]\geq C_{0,a_2}.
\end{equation}
On the other hand, by Corollary \ref{cro:monotonicity-Ca} and $a_1>0$, we have that
$C_{a_1,a_2}<C_{0,a_2}$, a contradiction.
\end{proof}

\begin{lemma}\label{lemma:20210414-l1}
Under the assumptions of Lemma \ref{lemma:Lagrange-multiplier-bounded}, up to a subsequence, we assume that $(\lambda_{1,\vartheta}, \lambda_{2,\vartheta})\rightarrow (\lambda_1,\lambda_2)$. If $N\geq 5$, we suppose further  \eqref{eq:parameters-paixu} and \eqref{eq:parameters-L}. Then
$\lambda_1>0, \lambda_2>0.$
\end{lemma}
\begin{proof}
We only prove $\lambda_1>0$. From Lemma \ref{lemma:no-vanish},
we deduce that there exists $y_\vartheta\in \mathbb{R}^N$ such that
$\displaystyle \limsup_{\vartheta\rightarrow \infty}\|\phi_\vartheta(x)\|_{L^2(Q^N)}=:\delta>0$, where
$\phi_\vartheta(x):=u_{1,\vartheta}(x+y_\vartheta).$
Up to a subsequence, there exists some $0\not\equiv \phi_0(x)\geq 0$ in $\mathbb{R}^N$  such that
$\phi_\vartheta(x)\rightharpoonup \phi_0(x)$ in $H^1(\mathbb{R}^N)$.
If $\lambda_1\leq 0$, we have that
\begin{equation}\label{eq:20210412-e5}
-\Delta \phi_0\geq \mu_1 \phi_{0}^{p_1-1}~\hbox{in}~\mathbb{R}^N.
\end{equation}
If $1\leq N\leq 4$, \eqref{eq:20210412-e5} implies that
\begin{equation}\label{eq:20210412-e6}
-\Delta \phi_0\geqslant 0\;\hbox{in}\;\mathbb{R}^N, \phi_0\in L^2(\mathbb{R}^N),
\end{equation}
a contradiction to \cite[Lemma A.2]{Ikoma2014}.
If $N\geq 5$, by \eqref{eq:parameters-L}, we see that $p_1-1\leq \frac{N}{N-2}$. Then follows by \cite[Theorem 8.4]{QuittnerSouplet2007} that  \eqref{eq:20210412-e5} has only trivial non-negative solution, which is also a contradiction.
Hence, we conclude that $\lambda_1>0$.
\end{proof}

In the following, without loss of generality, we may assume that $0<\lambda_1\leq \lambda_2$, $(u_{1,\vartheta}, u_{2,\vartheta})\rightharpoonup (u_1,u_2)$ weakly in $\mathcal{H}$ and $\lambda_{\iota,\vartheta}\rightarrow \lambda_{\iota}, \iota=1,2$.
Remark that
$\|u_\iota\|_2^2\leq a_\iota, \iota=1,2,$ and
notice that
$$(u_{1,\vartheta}, u_{2,\vartheta})\rightarrow (u_1,u_2)\;\hbox{strongly in}~L^2(\mathbb{R}^N)\times L^2(\mathbb{R}^N)~\hbox{if}~(\|u_1\|_2^2,\|u_2\|_2^2)=(a_1,a_2).$$
Recalling $\{(u_{1,\vartheta}, u_{2,\vartheta})\}$ is bounded in $\mathcal{H}$, by the well known Brezis-Lieb Lemma and Gagliardo-Nirenberg inequality, one can prove that
\begin{equation}\label{eq:20210412-e7}
\int_{\mathbb{R}^N}G(u_{1,\vartheta}, u_{2,\vartheta})\mathrm{d} x\rightarrow \int_{\mathbb{R}^N}G(u_1,u_2)\mathrm{d} x.
\end{equation}
On the other hand, by $V_\iota(x)\rightarrow 0$ as $|x|\rightarrow \infty$ and the weak lower semicontinuity of norm, it follows that
\begin{align*}
C_{a_1,a_2}=&\lim_{\vartheta\rightarrow \infty} J[u_{1,\vartheta}, u_{2,\vartheta}]\\
=&\lim_{\vartheta\rightarrow \infty} \left\{\sum_{\iota=1}^{2}\left[\|\nabla u_{\iota,\vartheta}\|_2^2-\|u_{\iota,\vartheta}\|_{V_\iota}^{2} \right] -\int_{\mathbb{R}^N}G(u_{1,\vartheta}, u_{2,\vartheta})\mathrm{d} x\right\}\\
\geq& \sum_{\iota=1}^{2}\left[\|\nabla u_{\iota}\|_2^2-\|u_{\iota}\|_{V_\iota}^{2} \right] -\int_{\mathbb{R}^N}G(u_{1}, u_{2})\mathrm{d} x\\
=&J[u_1,u_2]\geq C_{a_1,a_2}.
\end{align*}
Therefore, $(u_{1,\vartheta}, u_{2,\vartheta})\rightarrow (u_1,u_2)$ strongly in $\mathcal{H}$ and $J[u_1,u_2]=C_{a_1,a_2}$.

\

\noindent
{\bf Proof of Theorem \ref{thm:main-non-radial}:}
It is sufficient to prove that  $(\|u_1\|_2^2,\|u_2\|_2^2)=(a_1,a_2)$.
We argue by contradiction and suppose that
$$(\|u_1\|_2^2,\|u_2\|_2^2)=:(b_1,b_2)\neq (a_1,a_2).$$
Firstly we prove
\begin{lemma}\label{lemma:20210412-wl1}
 $(u_1,u_2)\neq (0,0)$ and $\displaystyle(b_1,b_2)\neq (0,0)$.
\end{lemma}
\begin{proof}
Suppose that $(u_1,u_2)=(0,0)$, then we also have $(b_1,b_2)= (0,0)$. Remark that $u_{\iota,\vartheta}\rightarrow 0$ in $L_{loc}^{p}(\mathbb{R}^N)$ for $2\leq p<2^*$ and $\iota=1,2$.
So by $V_\iota(x)\rightarrow 0$ as $|x|\rightarrow \infty$, recalling the observation in Remark \ref{remark:potential}-(ii), one can prove that
\begin{equation}\label{eq:20221216-e1}
\|u_{\iota,\vartheta}\|_{V_\iota}^{2}\rightarrow 0~\hbox{as}~\vartheta\rightarrow \infty, \iota=1,2.
\end{equation}
Hence,
\begin{align*}
C_{a_1,a_2}=&\lim_{\vartheta\rightarrow \infty} J[u_{1,\vartheta},u_{2,\vartheta}]
=\lim_{\vartheta\rightarrow \infty} I[u_{1,\vartheta},u_{2,\vartheta}]
\geq E_{a_1,a_2},
\end{align*}
a contradiction to Lemma \ref{lemma:min-negative} and Corollary \ref{cro:IE}.
\end{proof}

Secondly we prove
\begin{lemma}\label{lemma:20210412-xl1}
Let $Q^N:=[0,1]^N$. Then
\begin{equation}\label{eq:20210412-e8}
\limsup_{\vartheta\rightarrow \infty} \sup_{y\in \mathbb{R}^N} \left[\|u_{1,\vartheta}-u_{1}\|_{L^2(Q^N+y)}+\|u_{2,\vartheta}-u_{2}\|_{L^2(Q^N+y)}\right]>0.
\end{equation}
\end{lemma}
\begin{proof}
Suppose that $\displaystyle \limsup_{\vartheta\rightarrow \infty} \sup_{y\in \mathbb{R}^N} \left[\|u_{1,\vartheta}-u_{1}\|_{L^2(Q^N+y)}+\|u_{2,\vartheta}-u_{2}\|_{L^2(Q^N+y)}\right]=0$. Then up to a subsequence, by \cite[Lemma 1.21]{Willem1996} again, we have that
$u_{\iota, \vartheta}\rightarrow u_\iota$ strongly in $L^p(\mathbb{R}^N)$ for any $2<p<2^*$ and $\iota=1,2$. On the other hand, $V_\iota(x)\rightarrow 0$ as $|x|\rightarrow \infty$ and
$u_{\iota, \vartheta}\rightarrow u_\iota$ in $L_{loc}^{2}(\mathbb{R}^N)$ for $\iota=1,2$. Therefore,
$\int_{\mathbb{R}^N}G(u_{1,\vartheta},u_{2,\vartheta})\mathrm{d} x\rightarrow \int_{\mathbb{R}^N}G(u_{1},u_{2})\mathrm{d} x$ as $\vartheta\rightarrow \infty$ and then
\begin{equation}\label{eq:20210412-e9}
C_{b_1,b_2}\leq J[u_1,u_2]\leq \liminf_{\vartheta\rightarrow \infty} J[u_{1,\vartheta},u_{2,\vartheta}]=C_{a_1,a_2},
\end{equation}
a contradiction to Corollary \ref{cro:monotonicity-Ca}, due to the fact $(0,0)\neq (a_1-b_1, a_2-b_2)\in \mathbb{R}_+^2$.
\end{proof}

Then by Lemma \ref{lemma:20210412-xl1}, we may find $y_\vartheta\in \mathbb{R}^N$ such that
$|y_\vartheta|\rightarrow \infty$ and up to a subsequence,
$$\|u_{1,\vartheta}\|_{L^2(Q^N+y_\vartheta)}+\|u_{2,\vartheta}\|_{L^2(Q^N+y_\vartheta)}\rightarrow c_0>0.$$
So we may assume that
$\big(u_{1,\vartheta}(\cdot+y_\vartheta),u_{2,\vartheta}(\cdot+y_\vartheta)\big)\rightharpoonup (\omega_1,\omega_2)~\hbox{weakly in}~\mathcal{H}.$
We remark that
\begin{equation}\label{eq:20210412-e10}
(\omega_1,\omega_2)\not\equiv (0,0).
\end{equation}
By $|y_\vartheta|\rightarrow \infty$, one can prove that
\begin{equation}\label{eq:20210412-e11}
\|u_{\iota,\vartheta}-u_\iota-\omega_{\iota}(\cdot-y_\vartheta)\|_2^2 =\|u_{\iota,\vartheta}\|_2^2-\|u_\iota\|_2^2-\|\omega_\iota\|_2^2+o(1), \iota=1,2.
\end{equation}
In particular, put
$c_\iota:=\|\omega_\iota\|_2^2, \iota=1,2$. Then we have
\begin{equation}\label{eq:20210412-e12}
c_\iota:=\|\omega_\iota\|_2^2\leq \liminf_{\vartheta\rightarrow \infty} \left(\|u_{\iota,\vartheta}\|_2^2-\|u_\iota\|_2^2\right)=a_\iota-b_\iota, \iota=1,2.
\end{equation}
Next, we show
\begin{lemma}\label{lemma:20210412-xl2}
$(a_1,a_2)=(b_1,b_2)+(c_1,c_2), \ J[u_1,u_2]=C_{b_1,b_2},  \ I[\omega_1,\omega_2]=E_{c_1,c_2}$ and
\begin{equation}\label{eq:20210412-wbe2}
C_{a_1,a_2}=C_{b_1,b_2}+E_{c_1,c_2}.
\end{equation}
\end{lemma}
\begin{proof}
Suppose that
$$\lim_{\vartheta\rightarrow \infty} \sum_{\iota=1}^{2}\|u_{\iota,\vartheta} -u_\iota-\omega_{\iota}(\cdot-y_\vartheta)\|_2^2>0.$$
Then we see that
$$(b_1+c_1,b_2+c_2)\neq (a_1,a_2).$$
By Brezis-Lieb Lemma again, combing with \eqref{eq:20210412-e11} and $V_\iota(x)\rightarrow 0$ as $|x|\rightarrow \infty$, one can prove that
\begin{align}\label{eq:20210412-e13}
J[u_{1,\vartheta}, u_{2,\vartheta}]& -J[u_1,u_2]-J[\omega_{1}(\cdot-y_\vartheta), \omega_{2}(\cdot-y_\vartheta)]\nonumber\\
\quad &-J[u_{1,\vartheta} -u_1-\omega_{1}(\cdot-y_\vartheta),u_{2,\vartheta} -u_2-\omega_{2}(\cdot-y_\vartheta)]\rightarrow 0,
\end{align}
and
\begin{equation}\label{eq:20210412-e14}
\lim_{\vartheta\rightarrow \infty}J[\omega_{1}(\cdot-y_\vartheta), \omega_{2}(\cdot-y_\vartheta)]=I[\omega_1,\omega_2]\geq E_{c_1,c_2}.
\end{equation}
Put
$\delta_{\iota}:=\lim_{\vartheta\rightarrow \infty}\|u_{\iota,\vartheta} -u_\iota-\omega_{\iota}(\cdot-y_\vartheta)\|_2^2,$
then by \eqref{eq:20210412-e11}, we have
\begin{equation}\label{eq:20210412-e15}
\delta_{\iota}=a_\iota-b_\iota-c_\iota\geq 0, \iota=1,2.
\end{equation}
Since $\displaystyle u_{\iota,\vartheta} -u_\iota-\omega_{\iota}(\cdot-y_\vartheta)\rightarrow 0$ in $L_{loc}^{p}(\mathbb{R}^N)$ for $2\leq p<2^*$ and $V_\iota(x)\rightarrow 0$ as $|x|\rightarrow \infty$, we have that
\begin{align*}
&J[u_{1,\vartheta} -u_1-\omega_{1}(\cdot-y_\vartheta),u_{2,\vartheta} -u_2-\omega_{2}(\cdot-y_\vartheta)]\\
&-I[u_{1,\vartheta} -u_1-\omega_{1}(\cdot-y_\vartheta),u_{2,\vartheta} -u_2-\omega_{2}(\cdot-y_\vartheta)]\\
=&-\frac{1}{2}\sum_{\iota=1}^{2} \|u_{\iota,\vartheta} -u_\iota-\omega_{\iota}(\cdot-y_\vartheta)\|_{V_\iota}^{2}\rightarrow 0.
\end{align*}
And then it follows Lemma \ref{lemma:CE-continuous} that
\begin{align}\label{eq:20210412-e16}
&\liminf_{\vartheta\rightarrow \infty} J[u_{1,\vartheta} -u_1-\omega_{1}(\cdot-y_\vartheta),u_{2,\vartheta} -u_2-\omega_{2}(\cdot-y_\vartheta)]\nonumber\\
=&\liminf_{\vartheta\rightarrow \infty} I[u_{1,\vartheta} -u_1-\omega_{1}(\cdot-y_\vartheta),u_{2,\vartheta} -u_2-\omega_{2}(\cdot-y_\vartheta)]
\geq E_{\delta_1,\delta_2}.
\end{align}

By the formulas  \eqref{eq:20210412-e13}- \eqref{eq:20210412-e16},
we have that
\begin{align}\label{eq:20210412-be1}
C_{a_1,a_2}=&\lim_{\vartheta\rightarrow \infty}J[u_{1,\vartheta}, u_{2,\vartheta}] \nonumber\\
=&\lim_{\vartheta\rightarrow \infty}\left\{J[u_1,u_2]+J[\omega_{1}(\cdot-y_\vartheta), \omega_{2}(\cdot-y_\vartheta)]\right.\nonumber\\
&\left.\quad\quad+J[u_{1,\vartheta} -u_1-\omega_{1}(\cdot-y_\vartheta),u_{2,\vartheta} -u_2-\omega_{2}(\cdot-y_\vartheta)]\right\}\\
\geq&C_{b_1,b_2} +E_{c_1,c_2}+E_{\delta_1,\delta_2}.\nonumber
\end{align}

However, by Corollary \ref{cro:IE}, we may assume that $E_{c_1,c_2}$ and $E_{\delta_1,\delta_2}$ are attained by some nonnegative symmetry decreasing functions $(\phi_1,\phi_2)$ and $(\psi_1,\psi_2)$ respectively.
We recall the coupled rearrangement of $u$ and $v$, which is introduced by  M. Shibata\cite{Shibata2017}.
Then one can see that
$$(\phi_1\star \psi_1, \phi_2\star \psi_2)\in S_{c_1+\delta_1,c_2+\delta_2}$$
and
$$E_{c_1+\delta_1,c_2+\delta_2}\leq I[\phi_1\star \psi_1, \phi_2\star \psi_2]<I[\phi_1,\phi_2]+I[\psi_1,\psi_2]=E_{c_1,c_2}+E_{\delta_1,\delta_2},$$
see \cite[Theorem 2.4]{Shibata2017}.

Therefore,
$$C_{b_1,b_2} +E_{c_1,c_2}+E_{\delta_1,\delta_2}>C_{b_1,b_2}+E_{c_1+\delta_1,c_2+\delta_2}
\geq C_{b_1+c_1+\delta_1, b_2+c_2+\delta_2}=C_{a_1,a_2},$$
a contradiction to \eqref{eq:20210412-be1}. Hence, we prove that $\delta_\iota=0, \iota=1,2$. That is,
$$(a_1,a_2)=(b_1,b_2)+(c_1,c_2).$$
Again by \eqref{eq:20210412-e14}, \eqref{eq:20210412-be1} and
$$ C_{b_1,b_2} +E_{c_1,c_2}\geq C_{b_1+c_1,b_2+c_2}=C_{a_1,a_2},$$
we see that
$$J[u_1,u_2]=C_{b_1,b_2}, \  I[\omega_1,\omega_2]=E_{c_1,c_2} \ \ and \ \ C_{a_1,a_2}=C_{b_1,b_2} +E_{c_1,c_2}. $$
\end{proof}

Now, we see that on the contrary hypothesis $(b_1,b_2)\neq (a_1,a_2)$,
we obtain that $(0,0)\neq (u_1, u_2), (0,0)\neq (\omega_1,\omega_2)$.
Recalling \eqref{eq:20210412-e1},  $|y_\vartheta|\rightarrow \infty$ and $V_\iota(x)\rightarrow 0$ as $|x|\rightarrow \infty$, we see that
$(u_1,u_2)\neq (0,0)$ is a nonnegative solution to \eqref{eq:system-20210329-1} and $(\omega_{1},\omega_2)\neq (0,0)$ is a nonnegative solution to \eqref{eq:system-20210329-2}. Then by Proposition \ref{prop:binding-inequality}, there exists some $(\phi_1,\phi_2)\in\mathcal{H}$ such that
$$
(\|\phi_1\|_2^2,\|\phi_2\|_2^2)=(\|u_1\|_2^2+\|\omega_1\|_2^2,\|u_2\|_2^2+\|\omega_2\|_2^2)~\hbox{and}
~J[\phi_1,\phi_2]<J[u_1,u_2]+I[\omega_1,\omega_2].
$$
We see that $(\phi_1,\phi_2)\in S_{a_1,a_2}$. On the other hand, by Lemma \ref{lemma:20210412-xl2}, we have $J[u_1,u_2]=C_{b_1,b_2}$ and $I[\omega_1,\omega_2]=E_{c_1,c_2}$. Hence,
$$C_{a_1,a_2}\leq J[\phi_1,\phi_2]<J[u_1,u_2]+I[\omega_1,\omega_2]=C_{b_1,b_2}+E_{c_1,c_2},$$
a contradiction to \eqref{eq:20210412-wbe2}.

To sum up, we can prove that $(b_1,b_2)=(a_1,a_2)$ and $J[u_1,u_2]=C_{a_1,a_2}$.


\s{Appendix A}\label{Appdeneix-A}
\renewcommand{\theequation}{A.\arabic{equation}}
In this appendix, we give the detailed  proof of Proposition \ref{prop:sharp-estimation}. We remark that the statements of (i) and (ii) is well known, and one can also refer to \cite[Lemma 3.1-(i) and (ii)]{Ikoma2021}. Hence, we only need to prove the statements of (iii) here.

Without loss of generality, we assume that $0<\lambda_1\leq \lambda_2$. By $r_1,r_2>1$, we have that $\bar{\lambda}_1=\lambda_1$ and $\bar{\lambda}_2=\min\left\{\lambda_2, \frac{r_1^2}{(2-r_2)_+^2}\lambda_1\right\}$.

(iii-1) Put $u(x):=u_1(x)+u_2(x)$, we have that
\begin{align*}
-\Delta u+\lambda_1 u=&\partial_1G(u_1,u_2)+\partial_2G(u_1,u_2)-V_1(x)u_1-V_2(x)u_2+(\lambda_1-\lambda_2)u_2\\
\leq&\partial_1G(u_1,u_2)+\partial_2G(u_1,u_2)-V_1(x)u_1-V_2(x)u_2.
\end{align*}
 It follows from  $V_1(x),V_2(x)\leq 0$ that
\begin{equation}\label{eq:20210328-xe1}
-\Delta u+[\lambda_1+V_1(x)+V_2(x)]u\leq \partial_1G(u_1,u_2)+\partial_2G(u_1,u_2).
\end{equation}
We can take some $R_1>0$ large enough such that
$$ u(x):=u_1(x)+u_2(x)<1\;\hbox{for}~|x|\geq R_1.$$
By \eqref{eq:partial-G-2}, we have that
\begin{align*}
\partial_1G(u_1,u_2)+\partial_2G(u_1,u_2)\leq & C_1\left(u_{1}^{p_1-1}+u_{1}^{r_1-1}(u_1+u_2)^{r_2}
+u_{2}^{q_1-1}+u_{2}^{r_2-1}(u_1+u_2)^{r_1}\right)\\
\leq& C_{R_1} u^{p}\;\hbox{for}~|x|\geq R_1,
\end{align*}
where $p:=\min\{p_1-1,q_1-1,r_1+r_2-1\}>1$.
Then by \eqref{eq:20210328-xe1}, we obtain that
\begin{equation}
-\Delta u+[\lambda_1+V_1(x)+V_2(x)-C_{R_1}u^{p-1}]u\leq 0~\hbox{for}\;|x|\geq R_1.
\end{equation}
So for any $0<\nu_2<\lambda_1$, by $\displaystyle\lim_{|x|\rightarrow \infty}V_\iota(x)=0, \iota=1,2$, we can take some $R_2>R_1$ such that
$$\lambda_1+V_1(x)+V_2(x)-C_{R_1}u^{p-1}>\frac{\nu_2+\lambda_1}{2}\;\hbox{for}\;|x|\geq R_2.$$
Hence, we have
\begin{equation}\label{eq:20210328-xe2}
-\Delta u+\nu_2 u\leq 0\;\hbox{in}\;B_{R_2}^{c}.
\end{equation}
Then by the comparison principle, we have that
$u(x)\leq C e^{-\sqrt{\nu_2}|x|}\;\hbox{in}\;B_{R_2}^{c}.$
Therefore, we can find some $C_{\nu_2}>0$ such that
$u_1(x)\leq u(x)\leq C_{\nu_2}e^{-\sqrt{\nu_2}|x|}\;\hbox{in}\;\mathbb{R}^N.$

On the other hand, by \eqref{eq:partial-G-1}, we see that $\partial_1G(u_1,u_2)\geq 0$, and thus
$$-\Delta u_1+\lambda_1u_1=\partial_1G(u_1,u_2)-V_1(x)u_1\geq 0=-\Delta G_{\lambda_1}+\lambda_1 G_{\lambda_1}\;\hbox{in}\;\mathbb{R}^N,$$
where $G_{\lambda_1}$ is the Green function for $-\Delta+\lambda_1$ in $\mathbb{R}^N$. It is known (see  \cite{GidasNiNirenberg1981,Ikoma2021} for instance) that
\begin{equation}\label{eq:Green-function}
\begin{cases}
G_{\lambda_1}\in C(\mathbb{R}^N\backslash\{0\})\cap L^1(\mathbb{R}^N), G_{\lambda_1}(x)>0\;\hbox{in}\;\mathbb{R}^N\backslash \{0\},\\
\lim_{|x|\rightarrow \infty}G_{\lambda_1}(x)|x|^{\frac{N-1}{2}}e^{\sqrt{\lambda_1}|x|}=\gamma_0\in(0,\infty).
\end{cases}
\end{equation}
 Using the comparison principle again, one can prove that
$$c_1(1+|x|)^{-\frac{N-1}{2}} e^{-\sqrt{\lambda_1}|x|}\leq u_1(x)\;\hbox{for all}\;x\in \mathbb{R}^N\;\hbox{and some}\;c_1>0,$$
and hence \eqref{eq:20210329-e1} holds.
Noting that for any $\nu_1>\lambda_1$, it is trivial to  find some $C_{\nu_1}>0$ such that
$$C_{\nu_1} e^{-\sqrt{\nu_1}|x|}\leq c_1(1+|x|)^{-\frac{N-1}{2}} e^{-\sqrt{\lambda_1}|x|}\;\hbox{in}\;\mathbb{R}^N.$$
Hence, for any $0<\nu_2<\lambda_1<\nu_1$, we can find some $C_{\nu_1}, C_{\nu_2}>0$ such that
$$C_{\nu_1} e^{-\sqrt{\nu_1}|x|}\leq u_1(x)\leq C_{\nu_2} e^{-\sqrt{\nu_2}|x|}\;\hbox{in}\;\mathbb{R}^N.$$
We note that $C_{\nu_2}\geq C_{\nu_1}$ is trivial by putting $x=0$. Hence, \eqref{eq:20210328-estimation1} holds and we finish the proof of $(i)$.

(iii-2) We note that \eqref{eq:20210329-e2} is trivial by the proof of (iii-1). Next, we focus on the proof of \eqref{eq:20210328-estimation2}. We divide the proof  into different cases.

{\bf Case for $r_2\geq 2$:} In this case, we have that $\bar{\lambda}_2=\lambda_2$.
Firstly, we take some $R_1>0$ large enough such that
$$u_1(x)+u_2(x)\leq 1, \ \  \partial_2G(u_1,u_2)\leq C_1(u_{2}^{q_1-1}+u_{2}^{r_2-1}[u_{1}+u_2]^{r_1})\;\hbox{for}\;|x|\geq R_1.$$
It follows that
\begin{equation}
-\Delta u_2+[\lambda_2+V_2(x)-C_1u_{2}^{q_1-2} -C_1u_{2}^{r_2-2}[u_{1}+u_2]^{r_1}]u_2\leq 0\;\hbox{for}\;|x|\geq R_1.
\end{equation}
By $r_2\geq 2$, we have that
$$\lim_{|x|\rightarrow \infty}\left[\lambda_2+V_2(x)-C_1u_{2}(x)^{q_1-2} -C_1u_{2}(x)^{r_2-2}[u_{1}(x)+u_2(x)]^{r_1}\right]=\lambda_2.$$
So for any $0<\nu_2<\lambda_2$, by the comparison principle, we can find some $C_{\nu_2}>0$ such that
$$u_2(x)\leq C_{\nu_2}e^{-\sqrt{\nu_2}|x|}\;\hbox{in}~\mathbb{R}^N.$$
And we note that the lower bound can be proved as (iii-1). Hence, \eqref{eq:20210328-estimation2} holds for $r_2\geq 2$.

{\bf Case for $1<r_2<2$:} We firstly prove the upper bound of $u_2(x)$ in \eqref{eq:20210328-estimation2}.
Noting that we can find some $C>0$ such that
$$\partial_2G(u_1(x),u_2(x))\leq C\left\{u_{2}(x)^{q_1-1}+u_{2}(x)^{r_2-1}[u_{1}(x)+u_2(x)]^{r_1}\right\}\;\hbox{in}\;\mathbb{R}^N.$$
Put
$f_0(x):=\partial_2G(u_1(x),u_2(x)).$
Observing that for any $\lambda\in (0,\lambda_1)$, there exists some $C_\lambda$ such that
$u_1(x)+u_2(x)\leq C_\lambda e^{-\sqrt{\lambda}|x|}\;\hbox{in}\;\mathbb{R}^N.$
So we see that
$$f_0(x)\leq C_\lambda (e^{-(q_1-1)\sqrt{\lambda}|x|}+e^{-(r_1+r_2-1)\sqrt{\lambda}|x|}).$$
Since both $q_1-1>1$ and $r_1+r_2-1>1$, we can take $\varepsilon\in (0,\frac{\lambda_1}{2})$ small enough, such that
$$\sqrt{\xi_{1,\varepsilon}}:=\min\{(q_1-1)\sqrt{\lambda_1-\varepsilon}, (r_1+r_2-1)\sqrt{\lambda_1-\varepsilon}\}>\sqrt{\lambda_1}.$$
Noting that
\begin{equation}\label{eq:20210328-we1}
-\Delta e^{-\sqrt{\nu}|x|} +(\lambda_2+V_2(x))e^{-\sqrt{\nu}|x|}=\left(\lambda_2+V_2(x)-\nu-\frac{N-1}{|x|}\sqrt{\nu}\right)e^{-\sqrt{\nu}|x|},
\end{equation}
for any $\nu\in (0, \min\{\lambda_2, \xi_{1,\varepsilon}\})$, one can find some $R_\nu>0$ large enough such that
\begin{equation}\label{eq:20210328-we2}
\left(\lambda_2+V_2(x)-\nu-\frac{N-1}{|x|}\sqrt{\nu}\right)e^{-\sqrt{\nu}|x|}\geq f_0(x)\;\hbox{for}\;|x|\geq R_\nu.
\end{equation}
Then by the comparison principle again, one can prove that for any $\nu\in (0, \min\{\lambda_2, \xi_{1,\varepsilon}\})$, there exists some $C_{\nu}>0$ such that
\begin{equation}\label{eq:20210328-we3}
u_2(x)\leq C_\nu e^{-\sqrt{\nu}|x|}\;\hbox{in}\;\mathbb{R}^N.
\end{equation}
In particular, if $ \xi_{1,\varepsilon}<\lambda_2$, one can improve the statement that
\begin{equation}\label{eq:20210328-we4}
u_2(x)\leq C_\varepsilon e^{-\sqrt{\xi_{1,\varepsilon}}|x|}\;\hbox{in}\;\mathbb{R}^N.
\end{equation}

So if $\xi_{1,\varepsilon}\geq \lambda_2$, the upper bound of $u_2(x)$ in \eqref{eq:20210328-estimation2} holds.
If $\xi_{1,\varepsilon}<\lambda_2$, noting that $\xi_{1,\varepsilon}>\lambda_1$, so by such argument, we can improve the decay of $u_2$ from $e^{-\sqrt{\lambda}|x|}$ with $\lambda\in (0,\lambda_1)$ to $e^{-\sqrt{\lambda}|x|}$ with $\lambda\in (0,\xi_{1,\varepsilon})$.
For each $\varepsilon\in (0,\frac{\lambda_1}{2})$, we set
\begin{equation}\label{eq:diedai-1}
\sqrt{\xi_{n+1,\varepsilon}}:=\min\left\{(q_1-1)\sqrt{\xi_{n,\varepsilon}}, r_1\sqrt{\lambda_1-\varepsilon}+(r_2-1)\sqrt{\xi_{n,\varepsilon}}\right\}.
\end{equation}
Suppose that $\xi_{n,\varepsilon}<\lambda_2$, then we have
$f_0(x)\leq C_{n+1,\varepsilon} e^{-\sqrt{\xi_{n+1,\varepsilon}}|x|}\;\hbox{in}\;\mathbb{R}^N.$
Apply a similar argument as \eqref{eq:20210328-we1}-\eqref{eq:20210328-we4},
one can prove that for any $\nu\in (0, \min\{\lambda_2, \xi_{n+1,\varepsilon}\})$,
\begin{equation}\label{eq:20210328-we5}
u_2(x)\leq C_\nu e^{-\sqrt{\nu}|x|}\;\hbox{in}\;\mathbb{R}^N.
\end{equation}
In particular, if $\xi_{n+1,\varepsilon}<\lambda_2$,
\begin{equation}\label{eq:20210328-we5}
u_2(x)\leq C_\varepsilon e^{-\sqrt{\xi_{n+1,\varepsilon}}|x|}\;\hbox{in}\;\mathbb{R}^N.
\end{equation}

By solving the equation
\begin{equation}
r_1\sqrt{\lambda_1}+(r_2-1)\sqrt{t}=\sqrt{t},
\end{equation}
we have that
$t=\frac{r_1^2}{(2-r_2)^2}\lambda_1.$
So if
$\lambda_3:=\frac{r_1^2}{(2-r_2)^2}\lambda_1>\lambda_2,$
we can take $\varepsilon_0>0$ small enough such that
$\frac{r_1^2}{(2-r_2)^2}(\lambda_1-\varepsilon_0)>\lambda_2.$
Since $q_1-1>1$, we have that
\begin{equation}\label{eq:diedai-2}
\sqrt{s}<(q_1-1)\sqrt{s}\;\hbox{for any}\;s>0
\end{equation}
and
\begin{equation}\label{eq:diedai-3}
\sqrt{s}<r_1\sqrt{\lambda_1-\varepsilon_0}+(r_2-1)\sqrt{s}\;
\hbox{iff}~s<\frac{r_1^2}{(2-r_2)^2}(\lambda_1-\varepsilon_0).
\end{equation}
So if $\xi_{n,\varepsilon_0}\leq \lambda_2$ for all $n\in \mathbb{N}$, then
by \eqref{eq:diedai-1}, \eqref{eq:diedai-2} and \eqref{eq:diedai-3}, we see that
$\{\xi_{n,\varepsilon_0}\}$ is an increasing sequence with upper bound $\lambda_2$.
Let $n\rightarrow \infty$, we have that
$$\sqrt{\xi_{\infty,\varepsilon_0}}=r_1\sqrt{\lambda_1-\varepsilon_0}+(r_2-1)\sqrt{\xi_{\infty,\varepsilon_0}},$$
which implies $\xi_{\infty,\varepsilon_0}=\frac{r_1^2}{(2-r_2)^2}(\lambda_1-\varepsilon_0)\leq \lambda_2$, a contradiction to the choice of $\varepsilon_0$.
Hence, there must exists some $n_0\in \mathbb{N}$ such that $\xi_{n_0,\varepsilon_0}>\lambda_2$.
Therefore, for the case of $1<r_2<2$ with $\frac{r_1^2}{(2-r_2)^2}\lambda_1>\lambda_2$, we have $\bar{\lambda}_2=\lambda_2$.
 And for any $0<\nu_2<\lambda_2$, we can find some $C_{\nu_2}>0$ such that
$u_2(x)\leq C_{\nu_2}e^{-\sqrt{\nu_2}|x|}\;\hbox{in}~\mathbb{R}^N.$
And we note that the lower bound can be proved as (iii-1). Hence, \eqref{eq:20210328-estimation2} holds for $1<r_2<2$ with $\frac{r_1^2}{(2-r_2)^2}\lambda_1>\lambda_2$.

Finally, we focus on the case of $1<r_2<2$ with $\frac{r_1^2}{(2-r_2)^2}\lambda_1\leq \lambda_2$. In this case, we have $\lambda_1<\lambda_2$ and $\bar{\lambda}_2=\frac{r_1^2}{(2-r_2)^2}\lambda_1$.
We take $\xi_1$ satisfying that
$\lambda_1<\xi_1<\frac{r_1^2}{(2-r_2)^2}\lambda_1\leq \lambda_2$
and fix an $\varepsilon_1\in (0,\frac{\lambda_1}{2})$ small enough such that
\begin{equation}
\lambda_1<\xi_1<\frac{r_1^2}{(2-r_2)^2}(\lambda_1-\varepsilon_1).
\end{equation}
For any $\nu\in (0,\bar{\lambda}_2)$, we can take $\varepsilon_0\in (0,\varepsilon_1)$ small enough such that
$$\nu<\frac{r_1^2}{(2-r_2)^2}(\lambda_1-\varepsilon_0).$$
We argue as above and define $\{\xi_{n,\varepsilon_0}\}$ inductively by \eqref{eq:diedai-1} with $\xi_{1,\varepsilon_0}:=\xi_1-\varepsilon_0$.
We {\bf claim} that there exists some $n_0\in \mathbb{N}$ such that
$\min\{\xi_{n_0,\varepsilon_0},\lambda_2\}=\xi_{n_0,\varepsilon_0}>\nu.$
Then we can find some $C_\nu>0$ such that
$u_2(x)\leq C_\nu e^{-\sqrt{\nu}|x|}\;\hbox{in}\;\mathbb{R}^N$, which is the result we are looking for. If the claim does not hold, then one can check that $\{\xi_{n,\varepsilon_0}\}_{n=1}^{\infty}$ is an increasing sequence such that
$$\xi_{n,\varepsilon_0}\uparrow \xi_{\infty,\varepsilon_0}:=\frac{r_1^2}{(2-r_2)^2}(\lambda_1-\varepsilon_0)>\nu,$$
a contradiction to the choice of $\varepsilon_0$.

If $\frac{r_1^2}{(2-r_2)^2}\lambda_1\geq  \lambda_2$, we have $\bar{\lambda}_2=\lambda_2$, the lower bound can be proved as (iii-1). So we suppose that
$\bar{\lambda}_2:=\frac{r_1^2}{(2-r_2)^2}\lambda_1< \lambda_2$.
We can prove the lower bound by a similar inductively way. Firstly, we note that for any $\nu>\lambda_2$, as in (iii-1), one can find some $C_\nu>0$ such that
\begin{equation}
u_2(x)\geq C_\nu e^{-\sqrt{\nu}|x|}\;\hbox{in}\;\mathbb{R}^N.
\end{equation}
Hence, we only need to consider $\bar{\lambda}_2<\nu\leq \lambda_2$ in the following.
Put
$$\sqrt{\eta_{1}}:=r_1\sqrt{\lambda_1}+(r_2-1)\sqrt{\lambda_2}<\sqrt{\lambda_2}.$$
Then by \eqref{eq:partial-G-1}, we have
\begin{equation}
f_0(x)\geq C_2(u_{2}^{q_1-1}+u_{2}^{r_2-1}u_{1}^{r_1})\geq \tilde{c}_1(1+|x|)^{-\alpha_1} e^{-\sqrt{\eta_{1}}|x|}\;\hbox{in}\;\mathbb{R}^N
\end{equation}
with
$\alpha_1:=(r_1+r_2-1)\frac{N-1}{2}>0.$ Here we use the result in (iii-1) since $u_1\not\equiv 0$.
Recalling the Green function of $G_{\lambda_2}$,
 by
$$-\Delta u_2+\lambda_2 u_2=-V_2(x)u_2+f_0(x)\;\hbox{in}\;\mathbb{R}^N,$$
and $V_2(x)\leq 0$, we have that
\begin{align*}
u_2(x)=&G_{\lambda_2}\star [-V_2(x)u_2+f_0(x)]\\
=&\int_{\mathbb{R}^N}G_{\lambda_2}(x-y)[-V_2(y)u_2(y)+f_0(y)]dy\\
\geq&\int_{\mathbb{R}^N}G_{\lambda_2}(x-y)f_0(y)dy\\
\geq&\tilde{c}_1\int_{\mathbb{R}^N}G_{\lambda_2}(x-y) (1+|y|)^{-\alpha_1}e^{-\sqrt{\eta_{1}}|y|}dy
\end{align*}
Then  by Corollary \ref{cro:20210329-zcro1}, we have that
\begin{align*}
&\lim_{r\rightarrow \infty}(1+r)^{\alpha_1}e^{\sqrt{\eta_{1}} r}\int_{\mathbb{R}^N} G_{\lambda_2}(r\omega -y) (1+|y|)^{-\alpha_1}e^{-\sqrt{\eta_{1}}|y|}dy\\
=&\int_{\mathbb{R}^N}G_{\lambda_2}(y) e^{\sqrt{\eta_{1}}\omega\cdot y}dy>0\;\hbox{uniformly with respect to}~\omega\in \mathcal{S}^{N-1}.
\end{align*}
Hence, we can find some $c_1>0$ such that
\begin{equation}\label{eq:20210328-xbe1}
u_2(x)\geq c_1 (1+|x|)^{-\alpha_1}e^{-\sqrt{\eta_{1}}|x|}\;\hbox{in}\;\mathbb{R}^N.
\end{equation}
If $\nu> \eta_{1}$, by \eqref{eq:20210328-xbe1}, we can find some $C_\nu>0$ such that
$u_2(x)\geq C_\nu e^{-\sqrt{\nu}|x|}\;\hbox{in}\;\mathbb{R}^N.$
If $\eta_{n}\geq \nu$,  inductively, we define
\begin{equation}\label{eq:20210328-xbe2}
\sqrt{\eta_{n+1}}:=r_1\sqrt{\lambda_1}+(r_2-1)\sqrt{\eta_{n}}.
\end{equation}
We can find some $c_n>0, \alpha_n>0$ such that
\begin{equation}
u_2(x)\geq c_n (1+|x|)^{-\alpha_n}e^{-\sqrt{\eta_{n}}|x|}\;\hbox{in}\;\mathbb{R}^N.
\end{equation}
If there exists some $n_0\in \mathbb{N}$ such that $\eta_{n_0}<\nu $, our assertion holds. Suppose that $\eta_{n}\geq  \nu$ for all $n\in \mathbb{N}$,  we have
$\eta_n> \frac{r_1^2}{(2-r_2)^2}\lambda_1$ for all $n\in \mathbb{N}$.
Observing that
\begin{equation}\label{eq:20210328-whhe1}
r_1\sqrt{\lambda_1}+(r_2-1)\sqrt{t}<\sqrt{t}\;\hbox{iff} ~t>\frac{r_1^2}{(2-r_2)^2} \lambda_1,
\end{equation}
 we see that $\{\eta_{n}\}_{n=1}^{\infty}$ is an decreasing sequence. Thus, $\eta_{n}\leq \eta_1<\lambda_2$ holds for all $n\in \mathbb{N}$. Hence, Corollary \ref{cro:20210329-zcro1} is applied for every step in the induction. Let $n\rightarrow \infty$, we obtain that
$$\frac{r_1^2}{(2-r_2)^2}\lambda_1=\eta_{\infty}:=\lim_{n\rightarrow \infty}\eta_{n}\geq \nu,$$
a contradiction to $\nu>\bar{\lambda}_2$. We complete the proof.


\s{Appendix B}
\renewcommand{\theequation}{B.\arabic{equation}}
In this appendix, we give the proof of Lemma \ref{lemma:intersection-estimation}.

(i) If $u_1\equiv 0$ or $\omega_1\equiv 0$, then we have $\sigma_{1,R}\equiv 0$ and $\tau_{1,R}\equiv 1$. So the assertion of (i) holds automatically. Next, we assume that $u_1>0$ and $\omega_1>0$ in $\mathbb{R}^N$. By \cite[Lemma 3.4-(i)]{Ikoma2021} and the formula \eqref{eq:20210401-be2} that
\begin{align*}
&-\frac{\mu_i}{p_i}\tau_{1,R}^{p_i}\int_{\mathbb{R}^N}(u_1+\omega_{1,R})^{p_i}\mathrm{d} x
\leq -\frac{\mu_i}{p_i}\left(1-\frac{p_i}{a_1+b_1} \sigma_{1,R}\right)\\
& \ \ \ \ \ \ \ \ \times \int_{\mathbb{R}^N}\left[u_{1}^{p_i}+\omega_{1,R}^{p_i}+p_iu_{1}^{p_i-1}\omega_{1,R}+p_iu_{1}\omega_{1,R}^{p_i-1}-C_{p_i} u_{1}^{\frac{p_i}{2}} \omega_{1,R}^{\frac{p_i}{2}}\right]\mathrm{d} x+O(\sigma_{1,R}^{2}).
\end{align*}
By Proposition \ref{prop:sharp-estimation}, for any $\varepsilon_0>0$, there exists some $C_{\varepsilon_0}>0$ such that
$$u_1(x)+\omega_1(x)\leq C_{\varepsilon_0} e^{-\sqrt{\lambda_1-\varepsilon_0}|x|}\;\hbox{in}\;\mathbb{R}^N.$$
We can find some positive numbers $\varepsilon_is>0$ such that
\begin{equation}\label{eq:20210401-e7}
u_{1}^{\frac{p_i}{2}}(x) +\omega_{1}^{\frac{p_i}{2}}(x)\leq C_{\varepsilon_1} e^{-\sqrt{\frac{p_i^2}{4}\lambda_1-\varepsilon_1}|x|}
\end{equation}
and
\begin{equation}\label{eq:20210408-e1}
u_{1}^{p_i-1}(x)\leq C_{\varepsilon_2} e^{-\sqrt{(p_i-1)^2\lambda_1-\varepsilon_2}|x|}, \omega_{1}^{p_i-1}(x)\leq C_{\varepsilon_3} e^{-\sqrt{(p_i-1)^2\lambda_1-\varepsilon_3}|x|}.
\end{equation}
Then by Lemma  \ref{lemma:Ikoma-Miyamoto}, we have that
\begin{equation}\label{eq:20210408-e2}
\begin{cases}
\int_{\mathbb{R}^N}u_{1}^{p_i-1}\omega_{1,R}\mathrm{d} x\leq Ce^{-\sqrt{\lambda_1-\varepsilon_0}R},\\
\int_{\mathbb{R}^N}u_{1}\omega_{1,R}^{p_i-1}\mathrm{d} x\leq Ce^{-\sqrt{\lambda_1-\varepsilon_0}R},\\
\int_{\mathbb{R}^N}u_{1}^{\frac{p_i}{2}} \omega_{1,R}^{\frac{p_i}{2}}\mathrm{d} x\leq Ce^{-\sqrt{\frac{p_i^2}{4}\lambda_1-\frac{\varepsilon_1}{2}}R}.
\end{cases}
\end{equation}
It follows from  \eqref{20210329-we1} that
$c(1+R)^{-\frac{N-1}{2}} e^{-\sqrt{\lambda_1}R}\leq \sigma_{1,R}$.
Thanks to $p_i>2$,  for any $\theta_i(p_i)\in (1,\frac{p_i}{2})\cap (1,2)$, we can take suitable $\varepsilon_0,\varepsilon_1,\varepsilon_2,\varepsilon_3$ and conclude that
\begin{equation}\label{eq:20210401-e8}
\sigma_{1,R} \int_{\mathbb{R}^N}\left[u_{1}^{p_i-1}\omega_{1,R}+u_{1}\omega_{1,R}^{p_i-1}\right] \mathrm{d} x +
\int_{\mathbb{R}^N}u_{1}^{\frac{p_i}{2}} \omega_{1,R}^{\frac{p_i}{2}}\mathrm{d} x=o(\sigma_{1,R}^{\theta_i(p_i)})~\hbox{as}~R\rightarrow +\infty.
\end{equation}
Hence, the assertion of (i) holds.

(ii) We only need to prove the case of $u_2>0, \omega_2>0$.

 Suppose that $u_1\neq 0, \omega_1\neq 0$. Recalling the definition of $\bar{\lambda}_2$, by Proposition \ref{prop:sharp-estimation}, for any $\varepsilon>0$ small, there exists some $C_\varepsilon>0$ such that
$$u_2(x)+\omega_2(x)\leq C_\varepsilon e^{-\sqrt{\bar{\lambda}_2-\varepsilon}|x|}\;\hbox{in}\;\mathbb{R}^N.$$
Then thanks to  $\frac{q_j}{2}>1$ and $q_j-1>1$, we can similarly find some $\varepsilon>0$ such that
\begin{equation}\label{eq:20210401-e9}
u_{2}^{\frac{q_j}{2}}(x) +\omega_{2}^{\frac{q_j}{2}}(x)\leq C_\varepsilon e^{-\frac{q_j}{2}\sqrt{\bar{\lambda}_2-\varepsilon}|x|}
\end{equation}
\begin{equation}\label{eq:20210408-e3}
u_{2}^{q_j-1}(x)\leq C_\varepsilon e^{-(q_j-1)\sqrt{\bar{\lambda}_2-\varepsilon}|x|}, \omega_{2}^{q_j-1}(x)\leq C_\varepsilon e^{-(q_j-1)\sqrt{\bar{\lambda}_2-\varepsilon}|x|}.
\end{equation}
Then by Lemma  \ref{lemma:Ikoma-Miyamoto} again, we have that
\begin{equation}\label{eq:20210408-e4}
\begin{cases}
\int_{\mathbb{R}^N}u_{2}^{q_j-1}\omega_{2,R}\mathrm{d} x\leq Ce^{-\sqrt{\bar{\lambda}_2-\varepsilon}R},\\
\int_{\mathbb{R}^N}u_{2}\omega_{2,R}^{q_j-1}\mathrm{d} x\leq Ce^{-\sqrt{\bar{\lambda}_2-\varepsilon}R},\\
\int_{\mathbb{R}^N}u_{2}^{\frac{q_j}{2}} \omega_{2,R}^{\frac{q_j}{2}}\mathrm{d} x\leq Ce^{-\frac{q_j}{2}\sqrt{\bar{\lambda}_2-2\varepsilon}R}.
\end{cases}
\end{equation}
On the other hand, by Corollary \ref{cro:20210329-crow1}-(iii), we have that
$\tilde{C}_\varepsilon e^{-\sqrt{\bar{\lambda}_2+\frac{\varepsilon}{4}}R} \leq \sigma_{2,R}$. So for any $\theta(q_j)\in (1,\frac{q_j}{2})\cap(1,2)$, we can take $\varepsilon>0$ small enough and conclude that
\begin{equation}\label{eq:20210401-e10}
\sigma_{2,R} \int_{\mathbb{R}^N}\left[u_{2}^{q_j-1}\omega_{2,R}+u_{2}\omega_{2,R}^{q_j-1}\right] \mathrm{d} x +
\int_{\mathbb{R}^N}u_{2}^{\frac{q_j}{2}} \omega_{2,R}^{\frac{q_j}{2}}\mathrm{d} x=o(\sigma_{2,R}^{\theta_j(q_j)}).
\end{equation}
Hence, similar to (i), we can prove  the desired conclusion.

For the case of $u_1=\omega_1=0$, replace $\bar{\lambda}_2$ by $\lambda_2$, the argument above is valid.
Next, we consider the case of $u_1=0, \omega_1\neq 0$, and the case of $u_1\neq 0, \omega_1=0$.

(iii) If $0=\sigma_{1,R}<\sigma_{2,R}$, remark that either $u_1\equiv 0$ or $\omega_{1,R}\equiv 0$. Then  by \eqref{eq:20210409-xe1} and \eqref{eq:20210401-be2}, we have that
{\allowdisplaybreaks
\begin{align*}
&-\beta_k\tau_{1,R}^{r_{1,k}}\tau_{2,R}^{r_{2,k}}\int_{\mathbb{R}^N}(u_1+\omega_{1,R})^{r_{1,k}} (u_2+\omega_{2,R})^{r_{2,k}}\mathrm{d} x\\
\leq&-\beta_k  \left(1-\frac{r_{2,k}}{a_2+b_2} \sigma_{2,R}+O(\sigma_{2,R}^{2})\right)\\
&\times \int_{\mathbb{R}^N}\left[u_{1}^{r_{1,k}}u_{2}^{r_{2,k}}+r_{2,k}u_{1}^{r_{1,k}}u_{2}^{r_{2,k}-1}\omega_{2,R}
+r_{2,k}\omega_{1,R}^{r_{1,k}}\omega_{2,R}^{r_{2,k}-1}u_2+\omega_{1,R}^{r_{1,k}}\omega_{2,R}^{r_{2,k}}\right]\mathrm{d} x.
\end{align*}}
Noting that if $u_1\equiv 0$, we have
$$\int_{\mathbb{R}^N}r_{2,k}u_{1}^{r_{1,k}}u_{2}^{r_{2,k}-1}\omega_{2,R}\mathrm{d} x=0.$$
And if $u_1>0$, we can find some $0<\varepsilon<\frac{\bar{\lambda}_2}{2}$ and $C>0$ such that
$$u_{1}^{r_{1,k}}u_{2}^{r_{2,k}-1}(x)\leq C e^{-\sqrt{\varepsilon}|x|}, \omega_2(x)\leq C e^{-\sqrt{\bar{\lambda}_2-\varepsilon}|x|}.$$
It follows Lemma \ref{lemma:Ikoma-Miyamoto} that
$$\int_{R^N}u_{1}^{r_{1,k}}u_{2}^{r_{2,k}-1}\omega_{2,R}\mathrm{d} x\leq C e^{-\sqrt{\varepsilon}R}.$$
So we may find some $\theta\in (1,2)$ such that
\begin{equation}\label{eq:20210409-xe2}
\sigma_{2,R}\int_{R^N}u_{1}^{r_{1,k}}u_{2}^{r_{2,k}-1}\omega_{2,R}\mathrm{d} x=o(\sigma_{2,R}^{\theta}).
\end{equation}
Similarly, we can also prove that
\begin{equation}\label{eq:20240912-1014}
\sigma_{2,R}\int_{R^N}\omega_{1,R}^{r_{1,k}}\omega_{2,R}^{r_{2,k}-1}u_2\mathrm{d} x=o(\sigma_{2,R}^{\theta}).
\end{equation}

And thus
{\allowdisplaybreaks
\begin{align*}
&-\beta_k\tau_{1,R}^{r_{1,k}}\tau_{2,R}^{r_{2,k}}\int_{\mathbb{R}^N}(u_1+\omega_{1,R})^{r_{1,k}} (u_2+\omega_{2,R})^{r_{2,k}}\mathrm{d} x\\
\leq&-\beta_k \int_{\mathbb{R}^N}\left[u_{1}^{r_{1,k}}u_{2}^{r_{2,k}}+r_{2,k}u_{1}^{r_{1,k}}u_{2}^{r_{2,k}-1}\omega_{2,R}
+r_{2,k}\omega_{1,R}^{r_{1,k}}\omega_{2,R}^{r_{2,k}-1}u_2+\omega_{1,R}^{r_{1,k}}\omega_{2,R}^{r_{2,k}}\right]\mathrm{d} x\\
&+\frac{\beta_k r_{2,k}}{a_2+b_2} \sigma_{2,R} \int_{\mathbb{R}^N}\left[u_{1}^{r_{1,k}}u_{2}^{r_{2,k}}+\omega_{1,R}^{r_{1,k}}\omega_{2,R}^{r_{2,k}}\right]\mathrm{d} x+o(\sigma_{2,R}^{\theta}).
\end{align*}}

If $\sigma_{1,R}>0$, by Corollary \ref{cro:20210409-c1}, we have that
{\allowdisplaybreaks
\begin{align*}
&\int_{\mathbb{R}^N} (u_1+\omega_{1,R})^{r_{1,k}} (u_2+\omega_{2,R})^{r_{2,k}}\mathrm{d} x\\
\geq& \int_{\mathbb{R}^N}\left[ u_{1}^{r_{1,k}}u_{2}^{r_{2,k}}+r_{1,k}u_{1}^{r_{1,k}-1}u_{2}^{r_{2,k}}\omega_{1,R}
+r_{2,k}u_{1}^{r_{1,k}}u_{2}^{r_{2,k}-1}\omega_{2,R} +\omega_{1,R}^{r_{1,k}}\omega_{2,R}^{r_{2,k}}\right]\mathrm{d} x\\
&-C_{r_{1,k},r_{2,k},\eta}\int_{\mathbb{R}^N}\left[u_{1}^{1+\eta}\omega_{1,R}^{r_{1,k}-1-\eta}\omega_{2,R}^{r_{2,k}}
+u_{2}^{1+\eta}\omega_{1,R}^{r_{1,k}}\omega_{2,R}^{r_{2,k}-1-\eta}\right]\mathrm{d} x.
\end{align*}}
Basing on Proposition \ref{prop:sharp-estimation}-(iii), thanks to $1+\eta>1$ and $r_{1,k}-1-\eta+r_{2,k}>1$, we can find some $\eta_1>0$ such that
$$u_{1}^{1+\eta}(x)+\omega_{1}^{r_{1,k}-1-\eta}\omega_{2}^{r_{2,k}}(x)\leq Ce^{-\sqrt{\lambda_1+\eta_1}|x|}.$$
Then by Lemma \ref{lemma:Ikoma-Miyamoto} and \eqref{20210329-we1}, we may find some $\tilde{\theta}_1\in (1,2)$ such that
$$\int_{\mathbb{R}^N}u_{1}^{1+\eta}\omega_{1,R}^{r_{1,k}-1-\eta}\omega_{2,R}^{r_{2,k}}\mathrm{d} x=o(\sigma_{1,R}^{\tilde{\theta}_1})~\hbox{as}~R\rightarrow \infty.$$
Similarly, we can also find some $\tilde{\theta}_2\in (1,2)$ such that
$$\int_{\mathbb{R}^N}u_{2}^{1+\eta}\omega_{1,R}^{r_{1,k}}\omega_{2,R}^{r_{2,k}-1-\eta}\mathrm{d} x=o(\sigma_{1,R}^{\tilde{\theta}_2})~\hbox{as}~R\rightarrow \infty.$$
On the other hand, similar to the proof of \eqref{eq:20210409-xe2}, combing with Corollary \ref{cro:20210329-crow1}-(i), we may find some $\tilde{\theta}_3\in (1,2)$ such that
\begin{equation}\label{eq:20210409-we1}
\sigma_{\iota,R} \int_{\mathbb{R}^N}\left[r_{1,k}u_{1}^{r_{1,k}-1}u_{2}^{r_{2,k}}\omega_{1,R}
+r_{2,k}u_{1}^{r_{1,k}}u_{2}^{r_{2,k}-1}\omega_{2,R}\right]\mathrm{d} x=o(\sigma_{1,R}^{\tilde{\theta}_3}), \quad \iota=1,2.
\end{equation}
Hence, there exists some $\theta\in (1,2)$ such that
{\allowdisplaybreaks
\begin{align*}
&-\beta_k\tau_{1,R}^{r_{1,k}}\tau_{2,R}^{r_{2,k}}\int_{\mathbb{R}^N}(u_1+\omega_{1,R})^{r_{1,k}} (u_2+\omega_{2,R})^{r_{2,k}}\mathrm{d} x\\
\leq&-\beta_k \left(1-\frac{r_{1,k}}{a_1+b_1} \sigma_{1,R}+O(\sigma_{1,R}^{2})\right) \left(1-\frac{r_{2,k}}{a_2+b_2} \sigma_{2,R}+O(\sigma_{2,R}^{2})\right)\\
&\times \left\{\int_{\mathbb{R}^N}\left[ u_{1}^{r_{1,k}}u_{2}^{r_{2,k}}+r_{1,k}u_{1}^{r_{1,k}-1}u_{2}^{r_{2,k}}\omega_{1,R}
+r_{2,k}u_{1}^{r_{1,k}}u_{2}^{r_{2,k}-1}\omega_{2,R} +\omega_{1,R}^{r_{1,k}}\omega_{2,R}^{r_{2,k}}\right]\mathrm{d} x\right.\\
&\left.\quad-C_{r_{1,k},r_{2,k},\eta}\int_{\mathbb{R}^N}\left[u_{1}^{1+\eta}\omega_{1,R}^{r_{1,k}-1-\eta}\omega_{2,R}^{r_{2,k}}
+u_{2}^{1+\eta}\omega_{1,R}^{r_{1,k}}\omega_{2,R}^{r_{2,k}-1-\eta}\right]\mathrm{d} x\right\}\\
=&-\beta_k \int_{\mathbb{R}^N}\left[ u_{1}^{r_{1,k}}u_{2}^{r_{2,k}}+r_{1,k}u_{1}^{r_{1,k}-1}u_{2}^{r_{2,k}}\omega_{1,R}
+r_{2,k}u_{1}^{r_{1,k}}u_{2}^{r_{2,k}-1}\omega_{2,R} +\omega_{1,R}^{r_{1,k}}\omega_{2,R}^{r_{2,k}}\right]\mathrm{d} x\\
&+\frac{\beta_k r_{1,k}}{a_1+b_1}\sigma_{1,R}\int_{\mathbb{R}^N}\left[u_{1}^{r_{1,k}}u_{2}^{r_{2,k}}
+\omega_{1,R}^{r_{1,k}}\omega_{2,R}^{r_{2,k}}\right]\mathrm{d} x\\
&+\frac{\beta_k r_{2,k}}{a_2+b_2}\sigma_{2,R}\int_{\mathbb{R}^N}\left[u_{1}^{r_{1,k}}u_{2}^{r_{2,k}}
+\omega_{1,R}^{r_{1,k}}\omega_{2,R}^{r_{2,k}}\right]\mathrm{d} x+o(\sigma_{1,R}^{\theta}).
\end{align*}}


\s{Appendix C}
\renewcommand{\theequation}{C.\arabic{equation}}
In this appendix, we give the proof of Proposition \ref{prop:binding-inequality}.

For $\iota=1,2$, put
$a_\iota:=\|u_\iota\|_2^2, b_{\iota}:=\|\omega_\iota\|_2^2$.
Let
$\omega_{\iota, R}:=\omega_{\iota}(\cdot-R{\bf e_1}) $
and
$\sigma_{\iota,R}:=\int_{\mathbb{R}^N}u_\iota \omega_{\iota, R} \mathrm{d} x.$
Let
$\displaystyle
\tau_{\iota,R}:=\sqrt{\frac{a_\iota+b_\iota}{a_\iota+b_\iota+2\sigma_{\iota,R}}}.
$
For $\iota=1,2$, we have
{\allowdisplaybreaks
\begin{align*}
&\frac{\tau_{\iota,R}^{2}}{2}\|\nabla  (u_\iota+\omega_{\iota,R})\|_2^2-\frac{\tau_{\iota,R}^{2}}{2}\|u_\iota+ \omega_{\iota,R}\|_{V_\iota}^{2}\\
=&\frac{1}{2}\left(1-\frac{2\sigma_{\iota,R}}{a_{\iota}+b_{\iota}}+O(\sigma_{\iota,R}^{2})\right)
\Big\{\left(\|\nabla u_\iota\|_2^2+\|\nabla \omega_\iota\|_2^2+2\langle \nabla u_\iota, \nabla \omega_{\iota,R}\rangle_{L^2}\right)\\
&\quad \quad \quad \quad \quad \quad \quad \quad \quad \quad \quad \quad -\left(\|u_\iota\|_{V_\iota}^{2}+\|\omega_{\iota,R}\|_{V_\iota}^{2}+2\langle u_\iota,\omega_{\iota,R}\rangle_{V_\iota}\right)\Big\}\\
=&\frac{1}{2}\Big\{\left(\|\nabla u_\iota\|_2^2+\|\nabla \omega_\iota\|_2^2+2\langle \nabla u_\iota, \nabla \omega_{\iota,R}\rangle_{L^2}\right) -\left(\|u_\iota\|_{V_\iota}^{2}+\|\omega_{\iota,R}\|_{V_\iota}^{2}+2\langle u_\iota,\omega_{\iota,R}\rangle_{V_\iota}\right)\Big\}\\
&-\frac{\sigma_{\iota,R}}{a_{\iota}+b_{\iota}}\Big\{\left(\|\nabla u_\iota\|_2^2+\|\nabla \omega_\iota\|_2^2\right) -\left(\|u_\iota\|_{V_\iota}^{2}+\|\omega_{\iota,R}\|_{V_\iota}^{2}\right)\Big\}+O(\sigma_{\iota,R}^{2})\\
&+O\left(\sigma_{\iota,R}\left(\langle \nabla u_\iota, \nabla \omega_{\iota,R}\rangle_{L^2}-\langle u_\iota,\omega_{\iota,R}\rangle_{V_\iota}\right)\right)\\
=&\frac{1}{2}\Big\{\left(\|\nabla u_\iota\|_2^2+\|\nabla \omega_\iota\|_2^2+2\langle \nabla u_\iota, \nabla \omega_{\iota,R}\rangle_{L^2}\right) -\left(\|u_\iota\|_{V_\iota}^{2}+2\langle u_\iota,\omega_{\iota,R}\rangle_{V_\iota}\right)\Big\}\\
&-\frac{\sigma_{\iota,R}}{a_{\iota}+b_{\iota}}\Big\{\left(\|\nabla u_\iota\|_2^2+\|\nabla \omega_\iota\|_2^2\right) -\|u_\iota\|_{V_\iota}^{2}\Big\}+O(\sigma_{\iota,R}^{2})\\
&+O\left(\sigma_{\iota,R}\left(\langle \nabla u_\iota, \nabla \omega_{\iota,R}\rangle_{L^2}-\langle u_\iota,\omega_{\iota,R}\rangle_{V_\iota}\right)\right)-\left(\frac{1}{2}-\frac{\sigma_{\iota,R}}{a_{\iota}+b_{\iota}}\right)\|\omega_{\iota,R}\|_{V_\iota}^{2}.
\end{align*}}

Recalling that $(u_1,u_2)$ is a solution to \eqref{eq:system-20210329-1}, we have
\begin{equation}\label{eq:20210331-wbue1}
\begin{cases}
\langle \nabla u_1, \nabla \omega_{1,R}\rangle_{L^2}-\langle u_1,\omega_{1,R}\rangle_{V_1}=-\lambda_1 \sigma_{1,R}+\int_{\mathbb{R}^N}\partial_1G(u_1,u_2)\omega_{1,R}\mathrm{d} x,\\
\langle \nabla u_2, \nabla \omega_{2,R}\rangle_{L^2}-\langle u_2,\omega_{2,R}\rangle_{V_2}=-\lambda_2 \sigma_{2,R}+\int_{\mathbb{R}^N}\partial_2G(u_1,u_2)\omega_{2,R}\mathrm{d} x.
\end{cases}
\end{equation}
When $0<\lambda_1\leq \lambda_2$, due to
$p_i-1>1, q_j-1>1, r_{1,k}+r_{2,k}-1>1,$
apply a similar argument of Lemma \ref{lemma:intersection-estimation}, we can find some $\theta\in (1,2)$ such that if $\sigma_{1,R}> 0$,
\begin{align*}
&\sigma_{\iota,R} \left(\langle \nabla u_\iota, \nabla \omega_{\iota,R}\rangle_{L^2}-\langle u_\iota,\omega_{\iota,R}\rangle_{V_\iota}\right)\\
=&-\lambda_{\iota}\sigma_{\iota,R}^{2}+\sigma_{\iota,R}\int_{\mathbb{R}^N}\partial_\iota G(u_1,u_2)\omega_{\iota,R}\mathrm{d} x=o(\sigma_{1,R}^{\theta}), \iota=1,2,
\end{align*}
and if $\sigma_{1,R}=0<\sigma_{2,R}$,
\begin{align*}
\sigma_{2,R} \big(\langle \nabla u_2, \nabla \omega_{2,R}\rangle_{L^2}-\langle u_2,\omega_{2,R}\rangle_{V_2}\big)=o(\sigma_{2,R}^{\theta}).
\end{align*}

On the other hand , by $\sigma_{\iota,R}\rightarrow 0$ as $R\rightarrow \infty$, we consider $R>0$ large enough such that
$$\frac{1}{2}-\frac{\sigma_{\iota,R}}{a_{\iota}+b_{\iota}}>0, \iota=1,2.$$
Hence, if $\sigma_{1,R}\neq 0$, then for $\iota=1,2$, we have
\begin{align*}
&\frac{\tau_{\iota,R}^{2}}{2}\|\nabla  (u_\iota+\omega_{\iota,R})\|_2^2-\frac{\tau_{\iota,R}^{2}}{2}\|u_\iota+ \omega_{\iota,R}\|_{V_\iota}^{2}\\
\leq&\frac{1}{2}\Big\{\|\nabla u_\iota\|_2^2+\|\nabla \omega_\iota\|_2^2 -\|u_\iota\|_{V_\iota}^{2}\Big\}-\lambda_\iota \sigma_{\iota,R}+\int_{\mathbb{R}^N}\partial_\iota G(u_1,u_2)\omega_{\iota,R}\mathrm{d} x\\
&-\frac{\sigma_{\iota,R}}{a_{\iota}+b_{\iota}}\Big\{\left(\|\nabla u_\iota\|_2^2+\|\nabla \omega_\iota\|_2^2\right) -\|u_\iota\|_{V_\iota}^{2}\Big\}+o(\sigma_{1,R}^{\theta}).
\end{align*}
And if $\sigma_{1,R}=0$, we have that
\begin{align*}
&\frac{\tau_{2,R}^{2}}{2}\|\nabla  (u_2+\omega_{2,R})\|_2^2-\frac{\tau_{2,R}^{2}}{2}\|u_2+ \omega_{2,R}\|_{V_2}^{2}\\
\leq&\frac{1}{2}\Big\{\|\nabla u_2\|_2^2+\|\nabla \omega_2\|_2^2 -\|u_2\|_{V_2}^{2}\Big\}-\lambda_2 \sigma_{2,R}+\int_{\mathbb{R}^N}\partial_2 G(u_1,u_2)\omega_{2,R}\mathrm{d} x\\
&-\frac{\sigma_{2,R}}{a_{2}+b_{2}}\Big\{\left(\|\nabla u_2\|_2^2+\|\nabla \omega_2\|_2^2\right) -\|u_2\|_{V_2}^{2}\Big\}+o(\sigma_{2,R}^{\theta}),
\end{align*}
and
\begin{align*}
&\frac{\tau_{1,R}^{2}}{2}\|\nabla  (u_1+\omega_{1,R})\|_2^2-\frac{\tau_{1,R}^{2}}{2}\|u_1+ \omega_{1,R}\|_{V_1}^{2}\\
=&\frac{1}{2}\Big\{\|\nabla  (u_1+\omega_{1,R})\|_2^2-\|u_1+ \omega_{1,R}\|_{V_1}^{2}\Big\}\\
=&\begin{cases}
\frac{1}{2}\left[\|\nabla u_1\|_2^2-\|u_1\|_{V_1}^{2}\right]\quad &\hbox{if}~\omega_1=0,\\
\frac{1}{2}\left[\|\nabla \omega_1\|_2^2-\|\omega_{1,R}\|_{V_1}^{2}\right]\quad &\hbox{if}~u_1=0
\end{cases}\\
\leq&\frac{1}{2}\Big\{\|\nabla u_1\|_2^2+\|\nabla \omega_1\|_2^2 -\|u_1\|_{V_1}^{2}\Big\}.
\end{align*}

In the following, we shall prove that there exists some $R>0$ such that
\begin{equation}\label{eq:20210410-e2}
J(\tau_{1,R}(u_1+\omega_{1,R}), \tau_{2,R}(u_2+\omega_{2,R}))<J[u_1,u_2]+I[\omega_1,\omega_2].
\end{equation}
We divide the proof into three cases.

{\bf Case 1: $\sigma_{1,R}=0, \sigma_{2,R}=0$}.
 In this case, we have that $\tau_{1,R}\equiv \tau_{2,R}\equiv 1, \forall R>0$. Remark that either $u_1= \omega_2=0$ or $u_2=\omega_1=0$.
We only prove the case of $u_1=\omega_2=0$.
By
$$G(\omega_{1,R},u_2)>G(\omega_{1,R},0)+G(0,u_2)=G(\omega_{1,R},\omega_{2,R})+G(u_1,u_2),$$
we have that
{\allowdisplaybreaks
\begin{align*}
&J(\tau_{1,R}(u_1+\omega_{1,R}), \tau_{2,R}(u_2+\omega_{2,R}))=J[\omega_{1,R},u_2]\\
=&\frac{1}{2}\left[\|\nabla \omega_{1,R}\|_2^2-\|\omega_{1,R}\|_{V_1}^{2}\right]
+\frac{1}{2}\left[\|\nabla u_2\|_2^2-\|u_2\|_{V_2}^{2}\right]-\int_{\mathbb{R}^N}G(\omega_{1,R},u_2)\mathrm{d} x\\
<&\frac{1}{2}\left[\|\nabla \omega_{1,R}\|_2^2-\|\omega_{1,R}\|_{V_1}^{2}\right]-\int_{\mathbb{R}^N}G(\omega_{1,R},0)\mathrm{d} x\\
&+\frac{1}{2}\left[\|\nabla u_2\|_2^2-\|u_2\|_{V_2}^{2}\right]-\int_{\mathbb{R}^N}G(0,u_2)\mathrm{d} x\\
<&\frac{1}{2}\left[\|\nabla u_2\|_2^2-\|u_2\|_{V_2}^{2}\right]-\int_{\mathbb{R}^N}G(0,u_2)\mathrm{d} x\\
&+\frac{1}{2}\|\nabla \omega_{1,R}\|_2^2-\int_{\mathbb{R}^N}G(\omega_{1,R},0)\mathrm{d} x\\
=&J[0,u_2]+I[\omega_1,0].
\end{align*}}
The assertion \eqref{eq:20210410-e2} holds.

{\bf Case 2: $\sigma_{1,R}=0$ and $\sigma_{2,R}>0$}. In this case, we have that
{\allowdisplaybreaks
\begin{align*}
&\sum_{\iota=1}^{2}\frac{\tau_{\iota,R}^{2}}{2}\|\nabla  (u_\iota+\omega_{\iota,R})\|_2^2-\frac{\tau_{\iota,R}^{2}}{2}\|u_\iota+ \omega_{\iota,R}\|_{V_\iota}^{2}\\
\leq&\frac{1}{2}\sum_{\iota=1}^{2}\left[\|\nabla u_\iota\|_2^2+\|\nabla \omega_\iota\|_2^2 -\|u_\iota\|_{V_\iota}^{2}\right]-\lambda_2 \sigma_{2,R}+\int_{\mathbb{R}^N}\partial_2 G(u_1,u_2)\omega_{2,R}\mathrm{d} x\\
&-\frac{\sigma_{2,R}}{a_{2}+b_{2}}\left[\left(\|\nabla u_2\|_2^2+\|\nabla \omega_2\|_2^2\right) -\|u_2\|_{V_2}^{2}\right]+o(\sigma_{2,R}^{\theta}).
\end{align*}}
Then combing with Lemma \ref{lemma:intersection-estimation}, we can find some common $\theta\in (1,2)$ such that
{\allowdisplaybreaks
\begin{align*}
&J[\tau_{1,R}(u_1+\omega_{1,R}), \tau_{2,R}(u_2+\omega_{2,R})]=J[u_1+\omega_{1,R}, \tau_{2,R}(u_2+\omega_{2,R})]\\
=&\sum_{\iota=1}^{2}\frac{\tau_{\iota,R}^{2}}{2}\|\nabla  (u_\iota+\omega_{\iota,R})\|_2^2-\frac{\tau_{\iota,R}^{2}}{2}\|u_\iota+ \omega_{\iota,R}\|_{V_\iota}^{2}-\int_{\mathbb{R}^N}G\big(u_1+\omega_{1,R}, \tau_{2,R}(u_2+\omega_{2,R})\big)\mathrm{d} x\\
\leq&\frac{1}{2}\sum_{\iota=1}^{2}\left[\|\nabla u_\iota\|_2^2+\|\nabla \omega_\iota\|_2^2 -\|u_\iota\|_{V_\iota}^{2}\right]-\lambda_2 \sigma_{2,R}+\int_{\mathbb{R}^N}\partial_2 G(u_1,u_2)\omega_{2,R}\mathrm{d} x\\
&-\frac{\sigma_{2,R}}{a_{2}+b_{2}}\left[\left(\|\nabla u_2\|_2^2+\|\nabla \omega_2\|_2^2\right) -\|u_2\|_{V_2}^{2}\right]+o(\sigma_{2,R}^{\theta})\\
&-\sum_{i=1}^{\ell}\frac{\mu_i}{p_i}\left[\int_{\mathbb{R}^N}u_{1}^{p_i} +\omega_{1,R}^{p_i}\mathrm{d} x\right]\\
&-\sum_{j=1}^{m}\frac{\nu_j}{q_j}\left[\int_{\mathbb{R}^N}u_{2}^{q_j} +\omega_{2,R}^{q_j}\mathrm{d} x\right]
-\sum_{j=1}^{m}\nu_j\int_{\mathbb{R}^N} \left[u_{2}^{q_j-1}\omega_{2,R}+u_{2}\omega_{2,R}^{q_j-1}\right]\mathrm{d} x\\
&+\sum_{j=1}^{m}\frac{\nu_j\sigma_{2,R}}{a_2+b_2} \left[\int_{\mathbb{R}^N}u_{2}^{q_j} +\omega_{2,R}^{q_j}\mathrm{d} x\right]+o(\sigma_{2,R}^{\theta})\\
&-\sum_{k=1}^{n}\beta_k \int_{\mathbb{R}^N}\left[ \omega_{1,R}^{r_{1,k}}\omega_{2,R}^{r_{2,k}}+r_{2,k}\omega_{1,R}^{r_{1,k}}\omega_{2,R}^{r_{2,k}-1}u_2+ u_{1}^{r_{1,k}}u_{2}^{r_{2,k}}+r_{2,k}u_{1}^{r_{1,k}}u_{2}^{r_{2,k}-2}\omega_{2,R}\right]\mathrm{d} x\\
&+\sum_{k=1}^{n}\frac{\beta_k r_{2,k}}{a_2+b_2}\sigma_{2,R}\int_{\mathbb{R}^N}\left[u_{1}^{r_{1,k}}u_{2}^{r_{2,k}}
+\omega_{1,R}^{r_{1,k}}\omega_{2,R}^{r_{2,k}}\right]\mathrm{d} x+o(\sigma_{2,R}^{\theta}).
\end{align*}}
Noting that
{\allowdisplaybreaks
\begin{align*}
&\frac{1}{2}\sum_{\iota=1}^{2}\left[\|\nabla u_\iota\|_2^2+\|\nabla \omega_\iota\|_2^2 -\|u_\iota\|_{V_\iota}^{2}\right]
-\sum_{i=1}^{\ell}\frac{\mu_i}{p_i}\left[\int_{\mathbb{R}^N}u_{1}^{p_i} +\omega_{1,R}^{p_i}\mathrm{d} x\right]\\
&-\sum_{j=1}^{m}\frac{\nu_j}{q_j}\left[\int_{\mathbb{R}^N}u_{2}^{q_j} +\omega_{2,R}^{q_j}\mathrm{d} x\right]
-\sum_{k=1}^{n}\beta_k \int_{\mathbb{R}^N}\left[ u_{1}^{r_{1,k}}u_{2}^{r_{2,k}}
 +\omega_{1,R}^{r_{1,k}}\omega_{2,R}^{r_{2,k}}\right]\mathrm{d} x\\
 =&J[u_1,u_2]+I[\omega_1,\omega_2],
\end{align*}}
{\allowdisplaybreaks
\begin{align*}
&-\lambda_2\sigma_{2,R}-\frac{\sigma_{2,R}}{a_{2}+b_{2}}\Big\{\left(\|\nabla u_2\|_2^2+\|\nabla \omega_2\|_2^2\right) -\|u_2\|_{V_2}^{2}\Big\}\\
&+\sum_{j=1}^{m}\frac{\nu_j\sigma_{2,R}}{a_2+b_2} \left[\int_{\mathbb{R}^N}u_{2}^{q_j} +\omega_{2,R}^{q_j}\mathrm{d} x\right]\\
&+\sum_{k=1}^{n}\frac{\beta_k r_{2,k}}{a_2+b_2}\sigma_{2,R}\int_{\mathbb{R}^N}\left[u_{1}^{r_{1,k}}u_{2}^{r_{2,k}}
+\omega_{1,R}^{r_{1,k}}\omega_{2,R}^{r_{2,k}}\right]\mathrm{d} x\\
=&-\lambda_2\sigma_{2,R}-\frac{\sigma_{2,R}}{a_{2}+b_{2}}\Big\{\left(\|\nabla u_2\|_2^2+\|\nabla \omega_2\|_2^2\right) -\|u_2\|_{V_2}^{2}\Big\}\\
&+\frac{\sigma_{2,R}}{a_{2}+b_{2}}\int_{\mathbb{R}^N}\left[\partial_2G(u_1,u_2)u_2 +\partial_2G(\omega_{1,R},\omega_{2,R})\omega_{2,R}\right]\mathrm{d} x\\
=&-\lambda_2\sigma_{2,R}+\frac{\sigma_{2,R}}{a_{2}+b_{2}}\int_{\mathbb{R}^N}\left[\partial_2G(u_1,u_2)u_2 +\partial_2G(\omega_1,\omega_2)\omega_2\right]\mathrm{d} x\\
&-\frac{\sigma_{2,R}}{a_{2}+b_{2}}\Big\{-\lambda_2a_2+\int_{\mathbb{R}^N}\partial_2G(u_1,u_2)u_2\mathrm{d} x
-\lambda_2b_2+\int_{\mathbb{R}^N}\partial_2G(\omega_{1},\omega_{2})\omega_{2}\mathrm{d} x\Big\}\\
=&0
\end{align*}}
and
{\allowdisplaybreaks
\begin{align*}
&\int_{\mathbb{R}^N}\partial_2 G(u_1,u_2)\omega_{2,R}\mathrm{d} x-\sum_{j=1}^{m}\nu_j\int_{\mathbb{R}^N} \left[u_{2}^{q_j-1}\omega_{2,R}+u_{2}\omega_{2,R}^{q_j-1}\right]\mathrm{d} x\\
&-\sum_{k=1}^{n}\beta_k\int_{\mathbb{R}^N} \left[r_{2,k}u_{1}^{r_{1,k}}u_{2}^{r_{2,k}-1}\omega_{2,R}+r_{2,k}\omega_{1,R}^{r_{1,k}}\omega_{2,R}^{r_{2,k}-1}u_2 \right] \mathrm{d} x\\
=&-\sum_{j=1}^{m}\nu_j\int_{\mathbb{R}^N} u_{2}\omega_{2,R}^{q_j-1}\mathrm{d} x-\sum_{k=1}^{n}\beta_k\int_{\mathbb{R}^N}r_{2,k}\omega_{1,R}^{r_{1,k}}\omega_{2,R}^{r_{2,k}-1}u_2 \mathrm{d} x,
\end{align*}}
we obtain that
\begin{equation}\label{eq:20240911-1149}
\begin{aligned}
&J[\tau_{1,R}(u_1+\omega_{1,R}), \tau_{2,R}(u_2+\omega_{2,R})]-J[u_1,u_2]-I[\omega_1,\omega_2]\\
&\displaystyle \leq -\sum_{j=1}^{m}\nu_j\int_{\mathbb{R}^N} u_{2}\omega_{2,R}^{q_j-1}\mathrm{d} x-\sum_{k=1}^{n}\beta_k\int_{\mathbb{R}^N}r_{2,k}\omega_{1,R}^{r_{1,k}}\omega_{2,R}^{r_{2,k}-1}u_2 \mathrm{d} x+o(\sigma_{2,R}^{\theta}).
\end{aligned}
\end{equation}
So, if $\omega_1\equiv 0$, we have that
\begin{equation}\label{eq:20240911-1149}
\begin{aligned}
&J[\tau_{1,R}(u_1+\omega_{1,R}), \tau_{2,R}(u_2+\omega_{2,R})]-J[u_1,u_2]-I[\omega_1,\omega_2]\\
&\displaystyle \leq -\sum_{j=1}^{m}\nu_j\int_{\mathbb{R}^N} u_{2}\omega_{2,R}^{q_j-1}\mathrm{d} x+o(\sigma_{2,R}^{\theta}).
\end{aligned}
\end{equation}
While if $u_1\equiv 0$, we deduce further that
\begin{equation}\label{eq:20210410-wle1}
\begin{array}{ll}
&\displaystyle J[\tau_{1,R}(u_1+\omega_{1,R}), \tau_{2,R}(u_2+\omega_{2,R})]-J[u_1,u_2]-I[\omega_1,\omega_2]\\
&\displaystyle \leq -\sum_{j=1}^{m}\nu_j\int_{\mathbb{R}^N} u_{2}\omega_{2,R}^{q_j-1}\mathrm{d} x-\sum_{k=1}^{n}\beta_k\int_{\mathbb{R}^N}r_{2,k}\omega_{1,R}^{r_{1,k}}\omega_{2,R}^{r_{2,k}-1}u_2 \mathrm{d} x+o(\sigma_{2,R}^{\theta})\\
&=-\langle -\Delta \omega_{2,R}+\lambda_2 \omega_{2,R}, u_2\rangle +o(\sigma_{2,R}^{\theta})\\
&=-\langle -\Delta u_2+\lambda_2 u_2, \omega_{2,R}\rangle +o(\sigma_{2,R}^{\theta})\\
&=-\left[-\int_{\mathbb{R}^N}V_2(x)u_2\omega_{2,R}\mathrm{d} x+\sum_{j=1}^{m}\nu_j\int_{\mathbb{R}^N} u_{2}^{q_j-1}\omega_{2,R}\right]+o(\sigma_{2,R}^{\theta})\\
&\leq-\sum_{j=1}^{m}\nu_j\int_{\mathbb{R}^N} u_{2}^{q_j-1}\omega_{2,R}+o(\sigma_{2,R}^{\theta}).
\end{array}
\end{equation}
Recalling Corollary \ref{cro:20210329-crow1} that for any $\nu_1<\bar{\lambda}_2$, we have that
$\sigma_{2,R}\leq C_{\nu_1} e^{-\sqrt{\nu_1}R}.$
So by $\theta>1$, we can find some suitable $\bar{\eta}_1>0$ such that
\begin{equation}\label{eq:20210410-wle2}
\sigma_{2,R}^{\theta}\leq C e^{-\sqrt{\bar{\lambda}_2+\bar{\eta}_1}R}.
\end{equation}
On the other hand,  by Proposition \ref{prop:sharp-estimation} and Lemma \ref{lemma:Ikoma-Miyamoto}, if $u_1\equiv 0$, for any $\sqrt{\nu_2}>\min\{(q_1-1)\sqrt{\lambda_2}, \sqrt{\bar{\lambda}_2}\}=\sqrt{\bar{\lambda}_2}$, there exists $C_{\nu_2}>0$ such that
\begin{equation}\label{eq:20210410-wle3}
\int_{\mathbb{R}^N}u_{2}^{q_1-1} \omega_{2,R}\mathrm{d} x\geq C_{\nu_2} e^{-\sqrt{\nu_2}R}.
\end{equation}
In particular, we can take $\nu_2=\bar{\lambda}_2+\frac{\bar{\eta}_1}{2}\in (\bar{\lambda}_2, \bar{\lambda}_2+\bar{\eta}_1)$. Then one can see that
\begin{equation}
\sigma_{2,R}^{\theta} \left(\int_{\mathbb{R}^N}u_2 \omega_{2,R}^{q_1-1}\mathrm{d} x\right)^{-1}\rightarrow 0\;\hbox{as}\;R\rightarrow \infty,
\end{equation}
which implies that
$$-\sum_{j=1}^{m}\nu_j\int_{\mathbb{R}^N} u_{2}^{q_j-1}\omega_{2,R}\mathrm{d} x+o(\sigma_{2,R}^{\theta})<0\;\hbox{for large $R$}.$$
Hence, for $R$ large enough, by \eqref{eq:20210410-wle1} we have that
$$J[\tau_{1,R}(u_1+\omega_{1,R}), \tau_{2,R}(u_2+\omega_{2,R})]<J[u_1,u_2]+I[\omega_1,\omega_2].$$
The assertion \eqref{eq:20210410-e2} also holds.

Similarly, if $\omega_1\equiv 0$, by \eqref{eq:20240911-1149}, we can also deduce \eqref{eq:20210410-e2}.

{\bf Case 3: $\sigma_{1,R}>0$}.
In this case, we have that
{\allowdisplaybreaks
\begin{align*}
&\sum_{\iota=1}^{2}\left[\frac{\tau_{\iota,R}^{2}}{2}\|\nabla  (u_\iota+\omega_{\iota,R})\|_2^2-\frac{\tau_{\iota,R}^{2}}{2}\|u_\iota+ \omega_{\iota,R}\|_{V_\iota}^{2}\right]\\
\leq&\sum_{\iota=1}^{2}\left\{\frac{1}{2}\left[\|\nabla u_\iota\|_2^2+\|\nabla \omega_\iota\|_2^2 -\|u_\iota\|_{V_\iota}^{2}\right]-\lambda_\iota \sigma_{\iota,R}+\int_{\mathbb{R}^N}\partial_\iota G(u_1,u_2)\omega_{\iota,R}\mathrm{d} x\right.\\
&\left.\quad \quad -\frac{\sigma_{\iota,R}}{a_{\iota}+b_{\iota}}\left[\left(\|\nabla u_\iota\|_2^2+\|\nabla \omega_\iota\|_2^2\right) -\|u_\iota\|_{V_\iota}^{2}\right]\right\}+o(\sigma_{1,R}^{\theta}).
\end{align*}}
Then combining with Lemma \ref{lemma:intersection-estimation}, we can find some common $\theta\in (1,2)$ such that
{\allowdisplaybreaks
\begin{align*}
&J[\tau_{1,R}(u_1+\omega_{1,R}), \tau_{2,R}(u_2+\omega_{2,R})]\\
=&\sum_{\iota=1}^{2}\left[\frac{\tau_{\iota,R}^{2}}{2}\|\nabla  (u_\iota+\omega_{\iota,R})\|_2^2-\frac{\tau_{\iota,R}^{2}}{2}\|u_\iota+ \omega_{\iota,R}\|_{V_\iota}^{2}\right]-\int_{\mathbb{R}^N}G\big(\tau_{1,R}(u_1+\omega_{1,R}), \tau_{2,R}(u_2+\omega_{2,R})\big)\mathrm{d} x\\
\leq&\sum_{\iota=1}^{2}\left\{\frac{1}{2}\left[\|\nabla u_\iota\|_2^2+\|\nabla \omega_\iota\|_2^2 -\|u_\iota\|_{V_\iota}^{2}\right]-\lambda_\iota \sigma_{\iota,R}+\int_{\mathbb{R}^N}\partial_\iota G(u_1,u_2)\omega_{\iota,R}\mathrm{d} x\right.\\
&\left.\quad \quad -\frac{\sigma_{\iota,R}}{a_{\iota}+b_{\iota}}\left[\left(\|\nabla u_\iota\|_2^2+\|\nabla \omega_\iota\|_2^2\right) -\|u_\iota\|_{V_\iota}^{2}\right]\right\}+o(\sigma_{1,R}^{\theta})\\
&-\sum_{i=1}^{\ell}\frac{\mu_i}{p_i}\left[\int_{\mathbb{R}^N}u_{1}^{p_i} +\omega_{1,R}^{p_i}\mathrm{d} x\right]
-\sum_{i=1}^{\ell}\mu_i\int_{\mathbb{R}^N} \left[u_{1}^{p_i-1}\omega_{1,R}+u_{1}\omega_{1,R}^{p_i-1}\right]\mathrm{d} x\\
&+\sum_{i=1}^{\ell}\frac{\mu_i\sigma_{1,R}}{a_1+b_1} \left[\int_{\mathbb{R}^N}u_{1}^{p_i} +\omega_{1,R}^{p_i}\mathrm{d} x\right]+o(\sigma_{1,R}^{\theta})\\
&-\sum_{j=1}^{m}\frac{\nu_j}{q_j}\left[\int_{\mathbb{R}^N}u_{2}^{q_j} +\omega_{2,R}^{q_j}\mathrm{d} x\right]
-\sum_{j=1}^{m}\nu_j\int_{\mathbb{R}^N} \left[u_{2}^{q_j-1}\omega_{2,R}+u_{2}\omega_{2,R}^{q_j-1}\right]\mathrm{d} x\\
&+\sum_{j=1}^{m}\frac{\nu_j\sigma_{2,R}}{a_2+b_2} \left[\int_{\mathbb{R}^N}u_{2}^{q_j} +\omega_{2,R}^{q_j}\mathrm{d} x\right]+o(\sigma_{2,R}^{\theta})\\
&-\sum_{k=1}^{n}\beta_k \int_{\mathbb{R}^N}\left[ u_{1}^{r_{1,k}}u_{2}^{r_{2,k}}+r_{1,k}u_{1}^{r_{1,k}-1}u_{2}^{r_{2,k}}\omega_{1,R}
+r_{2,k}u_{1}^{r_{1,k}}u_{2}^{r_{2,k}-1}\omega_{2,R} +\omega_{1,R}^{r_{1,k}}\omega_{2,R}^{r_{2,k}}\right]\mathrm{d} x\\
&+\sum_{k=1}^{n}\frac{\beta_k r_{1,k}}{a_1+b_1}\sigma_{1,R}\int_{\mathbb{R}^N}\left[u_{1}^{r_{1,k}}u_{2}^{r_{2,k}}
+\omega_{1,R}^{r_{1,k}}\omega_{2,R}^{r_{2,k}}\right]\mathrm{d} x\\
&+\sum_{k=1}^{n}\frac{\beta_k r_{2,k}}{a_2+b_2}\sigma_{2,R}\int_{\mathbb{R}^N}\left[u_{1}^{r_{1,k}}u_{2}^{r_{2,k}}
+\omega_{1,R}^{r_{1,k}}\omega_{2,R}^{r_{2,k}}\right]\mathrm{d} x+o(\sigma_{1,R}^{\theta}).
\end{align*}}
Noting that firstly we have
{\allowdisplaybreaks
\begin{align*}
&\sum_{\iota=1}^{2}\frac{1}{2}\left[\|\nabla u_\iota\|_2^2+\|\nabla \omega_\iota\|_2^2 -\|u_\iota\|_{V_\iota}^{2}\right]
-\sum_{i=1}^{\ell}\frac{\mu_i}{p_i}\left[\int_{\mathbb{R}^N}u_{1}^{p_i} +\omega_{1,R}^{p_i}\mathrm{d} x\right]\\
&-\sum_{j=1}^{m}\frac{\nu_j}{q_j}\left[\int_{\mathbb{R}^N}u_{2}^{q_j} +\omega_{2,R}^{q_j}\mathrm{d} x\right]
-\sum_{k=1}^{n}\beta_k \int_{\mathbb{R}^N}\left[ u_{1}^{r_{1,k}}u_{2}^{r_{2,k}} +\omega_{1,R}^{r_{1,k}}\omega_{2,R}^{r_{2,k}}\right]\mathrm{d} x\\
=&J[u_1,u_2]+I[\omega_1,\omega_2].
\end{align*}}
Secondly, by
{\allowdisplaybreaks
\begin{align*}
&\sum_{\iota=1}^{2}\left\{ -\lambda_\iota \sigma_{\iota,R}-\frac{\sigma_{\iota,R}}{a_{\iota}+b_{\iota}}\left[\left(\|\nabla u_\iota\|_2^2+\|\nabla \omega_\iota\|_2^2\right)-\|u_\iota\|_{V_{\iota}}^{2}\right]\right\}\\
=&\sum_{\iota=1}^{2}\left\{ -\lambda_\iota \sigma_{\iota,R}-\frac{\sigma_{\iota,R}}{a_{\iota}+b_{\iota}}\left[-\lambda_\iota a_\iota+\int_{\mathbb{R}^N}\partial_\iota G(u_1,u_2)u_\iota \mathrm{d} x -\lambda_\iota b_\iota +\int_{\mathbb{R}^N}\partial_\iota G(\omega_1,\omega_2)\omega_\iota \mathrm{d} x\right]\right\}\\
=&\sum_{\iota=1}^{2}-\frac{\sigma_{\iota,R}}{a_{\iota}+b_{\iota}} \left[\int_{\mathbb{R}^N}\partial_\iota G(u_1,u_2)u_\iota \mathrm{d} x+\int_{\mathbb{R}^N}\partial_\iota G(\omega_1,\omega_2)\omega_\iota \mathrm{d} x\right],
\end{align*}}
we have that
{\allowdisplaybreaks
\begin{align*}
&\sum_{\iota=1}^{2}\left\{ -\lambda_\iota \sigma_{\iota,R}-\frac{\sigma_{\iota,R}}{a_{\iota}+b_{\iota}}\left[\left(\|\nabla u_\iota\|_2^2+\|\nabla \omega_\iota\|_2^2\right)-\|u_\iota\|_{V_{\iota}}^{2}\right]\right\}\\
&+\sum_{i=1}^{\ell}\frac{\mu_i\sigma_{1,R}}{a_1+b_1} \left[\int_{\mathbb{R}^N}u_{1}^{p_i} +\omega_{1,R}^{p_i}\mathrm{d} x\right]
+\sum_{j=1}^{m}\frac{\nu_j\sigma_{2,R}}{a_2+b_2} \left[\int_{\mathbb{R}^N}u_{2}^{q_j} +\omega_{2,R}^{q_j}\mathrm{d} x\right]\\
&+\sum_{k=1}^{n}\frac{\beta_k r_{1,k}}{a_1+b_1}\sigma_{1,R}\int_{\mathbb{R}^N}\left[u_{1}^{r_{1,k}}u_{2}^{r_{2,k}}
+\omega_{1,R}^{r_{1,k}}\omega_{2,R}^{r_{2,k}}\right]\mathrm{d} x\\
&+\sum_{k=1}^{n}\frac{\beta_k r_{2,k}}{a_2+b_2}\sigma_{2,R}\int_{\mathbb{R}^N}\left[u_{1}^{r_{1,k}}u_{2}^{r_{2,k}}
+\omega_{1,R}^{r_{1,k}}\omega_{2,R}^{r_{2,k}}\right]\mathrm{d} x
=0.
\end{align*}}
Thirdly,
{\allowdisplaybreaks
\begin{align*}
&\sum_{\iota=1}^{2}\int_{\mathbb{R}^N}\partial_\iota G(u_1,u_2)\omega_{\iota,R}\mathrm{d} x
-\sum_{i=1}^{\ell}\mu_i\int_{\mathbb{R}^N} \left[u_{1}^{p_i-1}\omega_{1,R}+u_{1}\omega_{1,R}^{p_i-1}\right]\mathrm{d} x\\
&-\sum_{j=1}^{m}\nu_j\int_{\mathbb{R}^N} \left[u_{2}^{q_j-1}\omega_{2,R}+u_{2}\omega_{2,R}^{q_j-1}\right]\mathrm{d} x\\
&-\sum_{k=1}^{n}\beta_k \int_{\mathbb{R}^N}\left[ r_{1,k}u_{1}^{r_{1,k}-1}u_{2}^{r_{2,k}}\omega_{1,R}
+r_{2,k}u_{1}^{r_{1,k}}u_{2}^{r_{2,k}-1}\omega_{2,R} \right]\mathrm{d} x\\
=&-\sum_{i=1}^{\ell}\mu_i\int_{\mathbb{R}^N}u_{1}\omega_{1,R}^{p_i-1}\mathrm{d} x
-\sum_{j=1}^{m}\nu_j\int_{\mathbb{R}^N}u_{2}\omega_{2,R}^{q_j-1}\mathrm{d} x.
\end{align*}}
Hence, we have that
\begin{align}\label{eq:20210411-xe1}
&J[\tau_{1,R}(u_1+\omega_{1,R}), \tau_{2,R}(u_2+\omega_{2,R})]
-J[u_1,u_2]-I[\omega_1,\omega_2]\nonumber\\
\leq &-\sum_{i=1}^{\ell}\mu_i\int_{\mathbb{R}^N}u_{1}\omega_{1,R}^{p_i-1}\mathrm{d} x
-\sum_{j=1}^{m}\nu_j\int_{\mathbb{R}^N}u_{2}\omega_{2,R}^{q_j-1}\mathrm{d} x+o(\sigma_{1,R}^{\theta}),
\end{align}
here we use the fact $o(\sigma_{2,R}^{\theta})=o(\sigma_{1,R}^{\theta})$ due to Corollary \ref{cro:20210329-crow1}.
 Then applying a similar argument as the Case 2, we can prove that
 \begin{equation}\label{eq:20210411-xe2}
 \sigma_{1,R}^{\theta} \left(\int_{\mathbb{R}^N}u_1 \omega_{1,R}^{p_1-1}\mathrm{d} x\right)^{-1}\rightarrow 0\;\hbox{as}\;R\rightarrow \infty.
 \end{equation}
Hence, by \eqref{eq:20210411-xe1} and \eqref{eq:20210411-xe2}, for $R>0$ large enough, we have that
$$J[\tau_{1,R}(u_1+\omega_{1,R}), \tau_{2,R}(u_2+\omega_{2,R})]<J[u_1,u_2]+I[\omega_1,\omega_2].$$
The assertion \eqref{eq:20210410-e2} also holds.

\vskip 0.2in
\noindent
{\bf Declarations}\\
{\bf Conflict of interest:} On behalf of all authors, the corresponding author states that there is no conflict of interest.

\noindent
{\bf Ethical Statement:}
The manuscript has not been previously published, is not currently submitted for review to any other journal, and will not be submitted elsewhere before a decision is made by your journal.

\noindent
{\bf Data Availability Statements}:
The datasets generated during and/or analysed during the current study are available from the author on reasonable request.

\end{document}